\documentclass{amsart}
\usepackage{amssymb}
\usepackage{bbm}
\usepackage{amsfonts}
\usepackage{mathrsfs}

\newtheorem{theorem}{Theorem}[section]
\newtheorem{lemma}[theorem]{Lemma}

\theoremstyle{definition}
\newtheorem{definition}[theorem]{Definition}
\newtheorem{example}[theorem]{Example}

\newtheorem{cor}[theorem]{Corollary}

\theoremstyle{remark}
\newtheorem{remark}[theorem]{Remark}

\numberwithin{equation}{section}

\def\and{\cap}
\def\u{\cup}
\def\bref#1{(\ref{#1})}
\def\proof{{\noindent\em Proof:} }

\DeclareFontFamily{U}{fsy}{} \DeclareFontShape{U}{fsy}{m}{n}{<->s*[.
9]psyr}{} \DeclareSymbolFont{der@m}{U}{fsy}{m}{n}
\DeclareMathSymbol{\diff}{\mathord}{der@m}{182}

\def\X{{\mathbb{X}}}
\def\Y{{\mathbb{Y}}}
\def\U{{\mathbb{U}}}
\def\V{{\mathbb{V}}}

\def\VB{{\mathbf{V}}}

\def\A{{\mathcal A}}

\def\C{{\mathcal C}}

\def\ff{{\mathcal F}}
\def\CI{{\mathcal{I}}}
\def\gg{{\mathcal G}}
\def\ee{{\mathcal E}}

\def\CV{{\mathbb{CV}}}

\def\P{{\mathbb P}}

\def\I{{\mathbb{I}}}    

\def\bu{{\mathbf{u}}}
\def\bv{{\mathbf{v}}}
\def\ba{{\mathbf{a}}}

\def\dZero{\mathbb{V}}
\def\dzero{\mathbb{V}}

\def\sat{\hbox{\rm{sat}}}
\def\asat{\hbox{\rm{asat}}}

\def\max{\hbox{\rm{max}}}

\def\gcd{\hbox{\rm{gcd}}}
\def\deg{\hbox{\rm{deg}}}

\def\init{\hbox{\rm{I}}}
\def\ord{\hbox{\rm{ord}}}

\def\H{\hbox{\rm{H}}}
\def\lead{\hbox{\rm{ld}}}
\def\sep{\hbox{\rm{S}}}

\def\dim{\hbox{\rm{dim}}}

\def\mod{\hbox{\rm{mod}}}

\def\ord{\hbox{\rm{ord}}}
\def\rk{\hbox{\rm{rk}}}

\def\den{\hbox{\rm{den}}}
\def\wdeg{\hbox{\rm{wdeg}}}

\def\CI{{\mathcal{I}}}

\def\Q{{\mathbb Q}}

\def\trdeg{\hbox{\rm{tr.deg}}}
\def\dtrdeg{\hbox{\rm{d.tr.deg}}}

\def\and{\cap}

\newcounter{bean}
\def\bl{\begin{list}{Step \arabic{bean}}{\usecounter{bean}}\labelwidth=34pt}
\def\el{\end{list}}

\def\deg{{\rm deg}}

\def\init{{\rm I}}

\def\normalization1{{\rm normalization1}}
\def\normalization{{\rm normalization}}

\def\irrfactor1{{\rm irrfactor1}}
\def\irrfactor{{\rm irrfactor}}

\def\sat{{\rm sat}}

\def\rank{{\rm rk}}

\begin{document}

\title[Differential Chow Form]{
Intersection Theory in \\Differential Algebraic Geometry:\\  Generic
Intersections  and \\ the Differential Chow Form}

\author{Xiao-Shan Gao, Wei Li, Chun-Ming Yuan}
\address{KLMM, Academy of Mathematics and Systems Science\\
Chinese Academy of Sciences, Beijing 100190, China}
\email{xgao@mmrc.iss.ac.cn, liwei@mmrc.iss.ac.cn, cmyuan@mmrc.iss.ac.cn}
%
\thanks{Partially supported by a National Key Basic Research Project of China (2011CB302400) and  by a grant from NSFC (60821002). }



\date{August 15, 2010.}
%


\begin{abstract}
In this paper, an intersection theory for generic differential
polynomials is presented. The intersection of an irreducible
differential variety of dimension $d$ and order $h$ with a generic
differential hypersurface of order $s$ is shown to be an irreducible
variety of dimension $d-1$ and order $h+s$. As a consequence, the
dimension conjecture for generic differential polynomials is proved.
Based on the intersection theory, the Chow form for an irreducible
differential variety is defined and most of the properties of the
Chow form in the algebraic case are established for its differential
counterpart. Furthermore, the generalized differential Chow form is
defined and its properties are proved. As an application of the
generalized differential Chow form, the differential resultant of
$n+1$ generic differential polynomials in $n$ variables is defined
and properties similar to that of the Macaulay resultant for
multivariate polynomials are proved.
\end{abstract}

\maketitle

\tableofcontents

\section{Introduction} \label{sec-i}

Differential algebra or differential algebraic geometry founded by
Ritt and Kolchin aims to study algebraic differential equations in a
similar way that polynomial equations are studied in commutative
algebra or algebraic geometry \cite{ritt,kol}. Therefore, the basic
concepts of commutative differential algebra are based on those of
commutative algebra.
%
%
%
An excellent survey on this subject can be found in \cite{bc1}.

It is known that, for many properties in algebraic geometry, their
differential counterparts are much more difficult to prove and some
of them are still open problems. For instance, many of the 16
questions proposed by Ritt in his classic book {\em Differential
Algebra} \cite[p.177]{ritt} are still not solved.
%
%
In this paper, two naturally connected problems in differential
algebraic geometry are studied: the differential dimension
conjecture for generic differential polynomials and the differential
Chow form.


The first part of the paper is concerned with the differential
dimension conjecture  which is one of the problems proposed by Ritt:
%
%
   Let $F_{1},    \ldots,    F_{r}$ be differential polynomials in
   $\mathcal {F}\{y_{1},$ $\ldots,  y_{n}\}$ with $r<n$, where $\mathcal {F}$
   is a differential field. If the
   differential variety  of the system $\{F_{1}, \ldots,  F_{r}\}$ is nonempty,
   then   each  of its components is of dimension at least $n-r$ \cite[p.178]{ritt}.

Ritt proved that the conjecture is correct when $r=1$, that is, any
component of a differential polynomial equation in $\mathcal
{F}\{y_{1},$ $\ldots,  y_{n}\}$ is of dimension $n-1$
\cite[p.57]{ritt}. The general differential dimension conjecture is
still open.
In \cite{cohn}, it is shown that the differential dimension
conjecture is closely related with Jacobi's bound for the order of
differential polynomial systems, which is another  well-known
conjecture in differential algebra.

In this paper, we consider the dimension and order for the
intersection of an irreducible differential variety with generic
differential hypersurfaces.
A differential polynomial $f$ is said to be generic of order $s$ and
degree $m$, if $f$ contains all the monomials with degree less than
or equal to $m$ in $y_{1},\ldots,y_{n}$ and their derivatives of
order up to $s$, and the coefficients of $f$ are differential
indeterminates.
%
%
A generic differential hypersurface is the set of solutions of a
generic differential polynomial.
We show that for generic differential hypersurfaces,
we can determine the dimension and order of their intersection with an irreducible
differential variety explicitly. More precisely, we will prove

\begin{theorem}\label{th-main1}
Let $\mathcal {I}$ be a prime differential polynomial ideal in
$\mathcal {F}\{y_{1},    \ldots,    y_{n}\}$ with dimension $d$ and
order $h$ and $f$ a generic differential polynomial with order $s$
and degree greater than zero. If $d>0$, then $\mathcal
{I}_1=[\mathcal {I}, f]$ is a prime differential polynomial ideal in
$\mathcal {F}\langle\bu_f \rangle\{y_{1}, \ldots, y_{n}\}$ with
dimension $d-1$ and order $h+s$, where $\bu_f$ is the set of
coefficients of $f$. And if $d=0$, $\mathcal {I}_1$ is the unit
ideal in $\mathcal {F}\langle\bu_f \rangle\{y_{1}, \ldots, y_{n}\}$.
\end{theorem}

As a direct consequence of this result, we show that the dimension conjecture
is valid for a system of generic differential polynomials. Furthermore, the order of the system
is also given explicitly.

Another purpose of studying the intersection of an irreducible
differential variety with generic differential hypersurfaces is to
establish the theory of the differential Chow form, which is the
concern of the second part of the paper consisting of Sections 4 to
6.

%

The Chow form, also known as the Cayley form or the Cayley-Chow
form, is a basic concept in algebraic geometry \cite{van,hodge}.
More recently, the Chow form also becomes a powerful tool in elimination theory.
This is not surprising, since the Chow form is a resultant in certain sense.
The Chow form was used as a tool to obtain deep results in
transcendental number theory by Nesterenko  \cite{nes1}
and  Philippon  \cite{ph1}.
Brownawell made a major breakthrough in elimination theory by
developing new properties of the Chow form and proving an effective
version of the Nullstellensatz with optimal bounds
\cite{brownawell}.
Gel'fand et al and Sturmfels started the sparse elimination theory
which is to study the Chow form and the resultant associated with
sparse polynomials \cite{gelfand,sturmfels}.
Eisenbud et al proposed a new expression for the Chow form via
exterior algebra and used it to give explicit formulas in many new cases \cite{eisenbud}.
Jeronimo et al gave a bounded probabilistic algorithm which
can be used to compute the Chow form, whose complexity
is polynomial in the size and the geometric degree of the input
equation system \cite{Complexitychowform}.
Other properties of the Chow form can be found in
\cite{Pedersen,Plumer,Seidenberg}.
Given the fact that the Chow form plays an important role in both
theoretic and algorithmic aspects of algebraic geometry and has
applications in many fields, it is worthwhile to develop the theory
of the differential Chow form and hope that it will play a similar
role as its algebraic counterpart.

Let $V$ be an irreducible differential variety of dimension $d$ in
an $n$-dimensional differential affine space and
$$\P_{i} = u_{i0} + u_{i1} y_1 + \cdots + u_{in}y_n \, (i=0,\ldots,d)$$
$d+1$ generic hyperplanes in differential variables
$y_1,\ldots,y_n$, where $u_{ij}\,(i=0,  \ldots,  d;$ $j=0, \ldots,
n)$ are differential indeterminates.
The differential Chow form of $V$ is roughly defined to be the
elimination differential polynomial in $u_{ij}$ for the intersection
of $V$ with $\P_i=0 (i=0,\ldots,d)$.
More intuitively, the differential Chow form of $V$ can be roughly
considered as the condition on the coefficients of $\P_i$ such that
these $d+1$ hyperplanes will meet $V$. We will show that most of the
properties of the Chow form in the algebraic case presented in
\cite{hodge,van} can be generalized to the differential case.
Precisely, we will prove

\begin{theorem} \label{th-main2}
Let $V$ be an irreducible differential variety with dimension $d$
and order $h$ over a differential field $\mathcal{F}$ and
$F(\bu_{0},\bu_{1},\ldots,\bu_{d})\in\mathcal{F}\{\bu_0,\bu_1,\ldots,\bu_d\}$ the Chow form of
$V$ where $\bu_{i}=(u_{i0},u_{i1},\ldots,u_{in})\,$ $(i=0,1,\ldots,d)$.
Then $F(\bu_{0},\bu_{1},\ldots,\bu_{d})$ has the
following properties:

1.
 $F(\bu_{0},\bu_{1},\ldots,\bu_{d})$ is differentially
 homogenous of the same degree in each  set  $\bu_{i}$
 and
 $\ord(F,u_{ij})=h$ for all $u_{ij}$ occurring in $F$.

2. 
%
$F(\bu_{0},\bu_{1},\ldots,\bu_{d})$  can be factored uniquely into
the following form
\begin{eqnarray}F(\bu_{0},\bu_{1},\ldots,\bu_{d})&=&A(\bu_{0},\bu_{1},\ldots,\bu_{d})
 \prod^g_{\tau=1}(u_{00}^{(h)}+\sum_{\rho=1}^n u_{0\rho}^{(h)}\xi_{\tau \rho}+t_{\tau}) \nonumber \\&=&A(\bu_{0},\bu_{1},\ldots,\bu_{d})
 \prod^g_{\tau=1}(u_{00}+\sum_{\rho=1}^n u_{0\rho}\xi_{\tau \rho})^{(h)} \nonumber
\end{eqnarray}
where $g=\deg(F,u_{00}^{(h)})$,  $\xi_{\tau \rho}$ are in a
differential extension field $\ff_\tau$ of $\mathcal{F}$,and
$(u_{00}+\sum_{\rho=1}^n u_{0\rho}\xi_{\tau \rho})^{(h)}$ is the
$h$-th derivative of $(u_{00}+\sum_{\rho=1}^n u_{0\rho}\xi_{\tau
\rho})$. The first $``="$ is obtained by factoring
$F(\bu_{0},\bu_{1},\ldots,\bu_{d})$ as an algebraic polynomial in
the variables $u_{00}^{(h)},u_{01}^{(h)},\ldots,u_{0n}^{(h)}$, while
the second $``="$ is a differential expression to be explained in
Section 4.4 of this paper.
%

3. $\Xi_\tau=(\xi_{\tau1},\ldots,\xi_{\tau n})\,(\tau=1,\ldots,g)$
are generic points of  $V$. And they are the only elements of $V$
lying on the differential hyperplanes
$\P_{\sigma}=0\,(\sigma=1,\ldots,d)$ as well as on the algebraic
hyperplanes $^a\P_{0}^{(l)}=0\,(l=0,\ldots,h-1)$.

4. Suppose that $\bu_i(i=0,\ldots,d)$ are differentially specialized
to sets $\bv_i$ of specific elements in  $\mathcal{F}$ and
$\overline{\P}_{i}\,(i=0,\ldots,d)$ are obtained by substituting
$\bu_i$ by $\bv_i$ in $\P_i$.
If\, $\overline{\P}_i=0(i=0,\ldots,d)$ meet $V$, then
$F(\bv_{0},\ldots,\bv_{d})$ $=0$. Furthermore, if
$F(\bv_{0},\ldots,\bv_{d})=0$ and $S_{F}(\bv_{0},\ldots,\bv_{d})\neq
0$, then the $d+1$ hyperplanes $\overline{\P}_{i}=0$
$(i=0,\ldots,d)$ meet $V$, where $S_{F}=\frac{\partial F}{\partial
u_{00}^{(h)}}$.
\end{theorem}

The number $g$ in the above theorem is called the leading
differential degree of $V$. From the third statement of the theorem,
we see that $V$ intersects with $\P_{\sigma}=0\,(\sigma=1,\ldots,d)$
and  $^a\P_{0}^{(l)}=0\,(l=0,\ldots,h-1)$ in exactly $g$ points.
So the leading differential degree satisfies similar properties with
the degree for an algebraic variety.

Furthermore, we prove that the four conditions given in Theorem
\ref{th-main2} are also the sufficient conditions for a differential
polynomial $F(\bu_{0},\bu_{1},\ldots,\bu_{d})$ to be the Chow form
for an order un-mixed differential variety, or a differential
algebraic cycle.
As a consequence of this result, we define the Chow quasi-variety
for a class of differential algebraic cycles in the sense that each
point in the Chow quasi-variety represents a differential algebraic
cycle $V$ in that class via the Chow form of $V$. These are clearly
generalizations of the algebraic Chow variety and algebraic cycles
\cite{gelfand,hodge}.

Note that both the differential Chow form and the generators of the
differential Chow quasi-variety are proven to be differentially
homogenous. Further developments of these concepts may need the
knowledge of differential projective space \cite{kol51}, which has
been ignored in certain degree in the past.

In \cite{ph1}, Philippon considered the intersection of a variety of
dimension $d$ with $d+1$ homogeneous polynomials with generic coefficients
and developed the theory for an elimination form which can be regarded as
a type of generalized Chow form.
In \cite{bost}, Bost, Gillet, and Soul\'e further generalized
the concept to generalized Chow divisors of cycles
and estimated their heights.
In this paper, we will introduce the generalized differential Chow form
which is roughly defined to be the elimination differential polynomial
obtained by intersecting an irreducible differential variety $V$ of dimension $d$
with $d+1$ generic differential hypersurfaces. We show that the generalized differential Chow
form satisfies similar properties to that given in Theorem
\ref{th-main2}.

As an application of the generalized differential Chow form,
we can define the differential resultant.
The differential resultant for two nonlinear differential
polynomials in one variable was studied by Ritt in \cite[p.47]{ritt0}.
General differential resultants were defined by
Carra' Ferro \cite{dres1,dres2} using Macaulay's definition of
algebraic resultant of polynomials.
But, the treatment in \cite{dres1} is not complete. For instance,
the differential resultant for two generic differential polynomials
with degree two and order one in one variable is always zero if
using the definition in  \cite{dres1}.
%
%
%
Differential resultants for linear ordinary differential
polynomials were studied by Rueda and Sendra in \cite{lres1}.
In this paper, a rigorous definition for the differential resultant
of $n+1$ generic differential polynomials in $n$ variables is given
as the generalized differential Chow form of the prime differential
ideal $\CI=[0]$. In this way, we obtain the following properties for
differential resultants, which are similar to that of the Macaulay
resultant for $n+1$ algebraic polynomials in $n$ variables.

\begin{theorem} \label{th-main3}
Let $\P_i$$(i=0,\ldots,n)$ be generic differential polynomials  in
$n$ differential variables $y_1,\ldots,y_n$ with orders $s_{i}$,
degrees $m_i$, and degree zero terms $u_{i0}$ respectively. Let
$R(\bu_{0},\bu_{1},\ldots,\bu_{n})$ be the differential resultant of
$\P_0,\ldots,\P_n$, where $\bu_i$ is the set of coefficients of
$\P_i$.
 Then

a)  $R(\bu_{0},\bu_{1},\ldots,\bu_{n})$ is differentially
homogeneous in each $\bu_i$ and is of order $h_i=s-s_i$ in $\bu_i$
$(i=0,\ldots,n)$ with $s=\sum_{l=0}^n s_l$.

b) There exist $\xi_{\tau \rho}(\rho=1,\ldots,n)$ in the
differential extension fields $\ff_\tau(\tau=1,\ldots,t_0)$ of
$\mathcal{F}$ such that
\begin{eqnarray*}
R(\bu_{0},\bu_{1},\ldots,\bu_{n})=A(\bu_{0},\bu_{1},\ldots,\bu_{n})\prod^{t_0}_{\tau=1}
 \P_0(\xi_{\tau 1},\ldots,\xi_{\tau n})^{(h_0)}
\end{eqnarray*}
where $A(\bu_{0},\bu_{1},\ldots,\bu_{n})$ is a differential
polynomial in $\bu_i$, $t_0=\deg(R,u_{00}^{(h_0)})$, $\P_0(\xi_{\tau
1},$ $\ldots,$ $\xi_{\tau n})^{(h_0)}$ is the $h_0$-th derivative
of $\P_0(\xi_{\tau 1},$ $\ldots,\xi_{\tau n})$, and
$(\xi_{\tau1},\ldots,$ $\xi_{\tau n})\,(\tau=1,\ldots,t_0)$ are
certain generic points of the zero dimensional prime differential
ideal $[\P_1,\ldots,\P_n]$.

c) The differential resultant can be written as a linear combination
of $\P_i$ and their derivatives up to the order $s-s_i\, (i=0,\ldots,n)$.
Precisely, we have
 $$R(\bu_{0},\bu_{1},\ldots,\bu_{n}) = \sum_{i=0}^n\sum_{j=0}^{s-s_i}
 h_{ij} \P_i^{(j)}.$$
In the above expression,  $h_{ij}\in \mathcal
{F}\langle\bu\rangle[y_1,\ldots,y_n,\ldots,y_{1}^{(s)},\ldots,y_{n}^{(s)}]$
have degrees at most $(sn+n)^2 D^{sn+n}+D(sn+n)$, where
$\bu=\cup_{i=0}^n\bu_{i}\setminus\{u_{00},\ldots,u_{n0}\}$,
$y_i^{(j)}$ is the $j$-th derivative of $y_i$, and
$D=\max\{m_0,m_1,\ldots,m_n\}$.

d) Suppose that $\bu_i(i=0,\ldots,n)$ are differentially specialized
to sets $\bv_i$ of specific elements in  $\mathcal{F}$ and
$\overline{\P}_{i}\,(i=0,\ldots,n)$ are obtained by substituting
$\bu_i$ by $\bv_i$ in $\P_i$.
If $\overline{\P}_{i}=0(i=0,\ldots,n)$ have a common solution, then
$R(\bv_{0},\ldots,\bv_{n})= 0$. On the other hand, if
$R(\bv_{0},\ldots,\bv_{n})= 0$ and $S_{R}(\bv_{0},\ldots,\bv_{n})\ne
0$, then $\overline{\P}_{i}=0 (i=0,\ldots,n)$ have a common
solution, where $S_{R}=\frac{\partial R}{\partial u_{00}^{(h_0)}}$.


\end{theorem}

%

As a prerequisite result, we prove a useful property of differential
specializations, which roughly asserts that if a set of differential
polynomial functions in a set of indeterminates  are differentially
dependent, then they are still differentially dependent when the
indeterminates are specialized to any concrete values. This property
plays a key role throughout this paper.
The algebraic version of this result is also a key result in algebraic
elimination theory (\cite[p.168]{hodge1}, \cite[p.161]{wu}).

It is not straightforward to extend the intersection theory for generic polynomials and the theory of Chow forms
from the algebraic case to the differential case.
Due to the complicated structure of differential polynomials,
most proofs in the algebraic case cannot
be directly used in the differential case.
In particular, we need to consider the orders
of differential polynomials, which is not an issue in the algebraic case.
For instance, the second property of the differential Chow form in Theorem \ref{th-main2}
has a different form as its algebraic counterpart.

One of the main tools used in the paper is the theory of
characteristic sets developed by Ritt \cite[p.47]{ritt}. The
algorithmic character of Ritt's work on differential algebra is
mainly due to the usage of characteristic sets. Properties of
characteristic sets  proved more recently in
\cite{boulier2010,ardm1,hubert,drpe,wu} will be also used in this
paper.

The rest of this paper is organized as follows.
In Section 2, we will present the notations and preliminary results
used in this paper.
%
In Section 3, the intersection theory for generic differential
polynomials is given and Theorem \ref{th-main1} is proved.
In Section 4, the Chow form for an irreducible differential variety is defined
and its properties will be proved. Basically, we will prove Theorem \ref{th-main2}.
In Section 5, necessary and sufficient conditions for a differential
polynomial to be the Chow form of an differential algebraic cycle is
given and the Chow quasi-variety for a class of differential
algebraic cycles is defined.
In Section 6, we present the theory of the generalized differential Chow form and
the differential resultant. Theorem \ref{th-main3} will be proved.
In Section 7, we present the conclusion and propose several problems
for further study.

\section{Preliminaries} \label{sec-p}

 In this section,   some basic notations and preliminary results in
 differential algebra will be given. For more details about differential algebra,
 please refer to \cite{ritt,kol,klmp,bc1,sit}.
\subsection{Differential polynomial algebra and Kolchin topology}
Let $\mathcal {F}$ be a fixed ordinary differential field of
characteristic zero, with a derivation  $\delta$. An element
$c\in\ff$ such that $\delta( c)=0$ is called a constant of $\ff.$ In
this paper, unless otherwise indicated, $\delta$ is kept fixed
during any discussion and we use primes and exponents $(i)$ to
indicate derivatives under $\delta$. Let $\Theta$ denote the free
commutative semigroup with unit (written multiplicatively) generated
by $\delta$.
Throughout the paper, we shall often use the prefix ``$\delta$-" as
a synonym of  ``differential" or ``differentially".

A typical example of differential field is $\mathbb{Q}(t)$ which is
the field of rational functions in variable $t$ with
$\delta=\frac{d}{dt}$.

%
%
Let $\mathcal{G}$ be a $\delta$-extension field of $\mathcal {F}$
and $S$  a subset of  $\mathcal{G}$. We will    denote respectively
by  $\mathcal {F}[S]$, $\mathcal {F}(S)$, $\mathcal {F}\{S\}$, and
$\mathcal {F}\langle S\rangle$ the smallest subring, the smallest
subfield, the smallest $\delta$-subring, and the smallest
$\delta$-subfield of $\mathcal {G}$ containing $\mathcal {F}$ and
$S$.  If we denote $\Theta(S)$ to be the    smallest subset of
$\mathcal {G}$ containing $S$ and stable under $\delta$,    we have
$\mathcal {F}\{S\}=\mathcal {F}[\Theta(S)]$ and $\mathcal {F}\langle
S\rangle=\mathcal {F}(\Theta(S))$. A $\delta$-extension field $\gg$
of $\ff$ is said to be finitely  generated if $\gg$ has a finite
subset $S$ such that $\gg=\ff\langle S\rangle$.

A subset $\Sigma$ of a $\delta$-extension field $\mathcal {G}$ of
$\mathcal {F}$ is said to be {\em $\delta$-dependent} over $\mathcal
{F}$
if the set $(\theta\alpha)_{\theta \in \Theta, \alpha\in\Sigma}$ is
algebraically dependent over $\mathcal {F}$, and is said to be {\em
$\delta$-independent} over $\mathcal {F}$, or to be a family of {\em
$\delta$-indeterminates}  over $\mathcal {F}$ (abbr.
$\delta$-$\ff$-indeterminates) in the contrary case.
In the case $\Sigma$ consists of one element $\alpha$, we say that
$\alpha$ is $\delta$-algebraic or $\delta$-transcendental over
$\mathcal {F}$ respectively. A maximal subset $\Omega$ of $\mathcal
{G}$ which is $\delta$-independent over $\mathcal {F}$ is said to be
a $\delta$-transcendence basis of $\mathcal {G}$ over $\mathcal
{F}$. We use $\dtrdeg \,\mathcal {G}/\mathcal {F}$ (see
\cite[p.105-109]{kol}) to denote the {\em $\delta$-transcendence
degree} of $\mathcal {G}$ over $\mathcal {F}$, which is the cardinal
number of $\Omega$. Considering $\mathcal {F}$ and $\mathcal {G}$ as
ordinary algebraic fields, we denote the algebraic transcendence
degree of $\mathcal {G}$ over $\mathcal {F}$ by $\trdeg\,\mathcal
{G}/\mathcal {F}$. If $S$ is a set of $\delta$-indeterminates over
$\ff$, then we call $\mathcal {F}\langle S\rangle$ a {\em pure
$\delta$-extension} of $\ff$.

A homomorphism $\varphi$ from a  $\delta$-ring $(\mathcal
{R},\delta)$ to a $\delta_1$-ring $(\mathcal {S},\delta_1)$ is a
{\em differential homomorphism} if $\varphi\circ
\delta=\delta_1\circ \varphi$. If $\mathcal {R}_0$ is a common
$\delta$-subring of $\mathcal {R}$ and $\mathcal {S}$ and the
homomorphism $\varphi$ leaves every element of $\mathcal {R}_0$
invariant, it is said to be over $\mathcal {R}_0$ and denoted by
$\delta$-$\mathcal {R}_0$-homomorphism. If, in addition $\mathcal
{R}$ is an integral domain and $\mathcal {S}$ is a $\delta$-field,
$\varphi$ is called a {\em $\delta$-specialization} of $\mathcal
{R}$ into $\mathcal {S}$.

A $\delta$-extension field  $\mathcal {E}$ of $\ff$ is called a {\em
universal $\delta$-extension field}, if for any finitely generated
$\delta$-extension field $\ff_1$ of $\ff$ in $\mathcal {E}$ and any
finitely generated $\delta$-extension field $\ff_2$ of $\ff_1$ not
necessarily in $\mathcal {E}$, $\ff_2$ can be embedded in $\mathcal
{E}$ over $\ff_1$, i.e. there exists a $\delta$-extension field of
$\ff_1$ in $\mathcal {E}$ that is $\delta$-isomorphic to $\ff_2$
over $\ff_1$. Such a $\delta$-universal extension field of $\ff$
always exists (\cite[Theorem 2, p. 134]{kol}). By definition, any
finitely generated $\delta$-extension field of $\ff$ can be embedded
over $\ff$ into $\mathcal {E}$, and $\mathcal {E}$ is a universal
$\delta$-extension field of every finitely generated
$\delta$-extension field of $\ff$. In particular, for any natural
number $n$, we can find in $\ee$ a subset of $\ee$ of cardinality
$n$ whose elements are $\delta$-independent over $\ff.$ Throughout
the present paper, $\ee$ stands for a fixed universal
$\delta$-extension field of $\ff$.

Now suppose $\Y=\{y_{1}, y_{2}, \ldots, y_{n}\}$ is a set of
$\delta$-indeterminates over $\ee$. For any $y\in\Y$, denote
$\delta^ky$ by $y^{(k)}.$ The elements of $\mathcal
{F}\{\Y\}=\mathcal {F}[y_j^{(k)}:j=1,\ldots,n;k\in \mathbb{N}]$ are
called {\em $\delta$-polynomials} over $\ff$ in $\Y$,
 and $\mathcal {F}\{\Y\}$ itself is called the {\em $\delta$-polynomial ring } over $\ff$ in $\Y$.
A $\delta$-polynomial ideal $\mathcal {I}$ in $\mathcal {F}\{\Y\}$
is an ordinary algebraic ideal which is closed under derivation,
 i.e. $\delta(\mathcal {I})\subset\mathcal {I}$.
In this paper, by $\delta$-ideals  we mean $\delta$-polynomial
ideals, and by a $\delta$-$\ff$-ideal we mean an ideal in $\mathcal
{F}\{\Y\}$. And a prime (resp. radical) $\delta$-ideal is a
$\delta$-ideal which is prime (resp. radical) as an ordinary
algebraic polynomial ideal.
For convenience, a prime $\delta$-ideal is assumed not to be the
unit ideal in this paper.

By a {\em $\delta$-affine space} we mean any one of the sets
$\ee^n\,(n\in \mathbb{N}).$ An element $\eta=(\eta_1,\ldots,\eta_n)$
of $\ee^n$ will be called a point. Let $\Sigma$ be a subset of
$\delta$-polynomials in $\mathcal {F}\{\Y\}$. A point
$\eta=(\eta_{1},\ldots,\eta_{n}) \in \mathcal {E}^n$ is called a
$\delta$-zero of $\Sigma$ if $f(\eta)=0$ for any $f \in\Sigma$. The
set of $\delta$-zeros of $\Sigma$ is denoted by $\dzero(\Sigma)$,
which is called a {\em $\delta$-variety} defined over $\ff$ (abbr.
$\delta$-$\ff$-variety). And for any $D\in \ff\{\Y\}$,
$\dzero(\Sigma/D) = \dzero(\Sigma)\setminus\dzero(D)$ is called a
{\em $\delta$-quasi-variety}. By convenience, we also call
$\cup_{i=1}^m \dzero(\Sigma_i/D_i)$ a $\delta$-quasi-variety, where
$\Sigma_i$ and $D_i$ are $\delta$-polynomial sets and
$\delta$-polynomials respectively. The $\delta$-varieties in $\ee^n$
(resp. the $\delta$-$\ff$-varieties in $\ee^n$ ) are the closed sets
in a topology called the {\em Kolchin topology} (resp. the Kolchin
$\ff$-topology).

For a $\delta$-$\ff$-variety $V$,  we denote $\mathbb{I}(V)$ to be
the set of all $\delta$-polynomials in $\ff\{\Y\}$ that vanish at
every point of $V$. Clearly, $\mathbb{I}(V)$ is a radical
$\delta$-ideal in $\ff\{\Y\}$. And there exists a bijective
correspondence between Kolchin $\ff$-closed sets and radical
$\delta$-ideals in $\ff\{\Y\}$. That is, for any $\delta$-$\ff$
variety $V$, $\dzero(\mathbb{I}(V))=V$ and for any radical
$\delta$-ideal $\mathcal{I}$ in $\ff\{\Y\}$,
$\mathbb{I}(\dzero(\mathcal {I}))=\mathcal {I}$.

Similarly as in algebraic geometry, an $\ff$-irreducible
$\delta$-variety can be defined. And there is a bijective
correspondence between $\ff$-irreducible $\delta$-varieties and
prime $\delta$-ideals in $\ff\{\Y\}$.
A point $\eta\in\V(\CI)$ is called a {\em generic point} of a prime
$\delta$-ideal $\CI\subset\ff\{\Y\}$ if for any $\delta$-polynomial
$p\in\ff\{\Y\}$ we have $p(\eta)=0 \Leftrightarrow p\in\CI$.
%
%
It is well known that \cite[p.27]{ritt}
\begin{lemma}\label{lm-gp}
A non-unit $\delta$-ideal is prime if and only if it has a generic
point.
\end{lemma}
And by the definition of universal $\delta$-fields, a prime
$\delta$-ideal over any finitely generated $\delta$-extension field
of $\ff$ has a generic point in $\ee^n$.
A point $\eta\in V$ is called a {\em generic point} of an
irreducible $\delta$-variety $V$ if $\eta$ is a generic point for
$\I(V)$.

Notice that irreducibility  depends on the base field over which the
$\delta$-polynomials are defined. For instance, $(\delta y_1)^2 - t$
is an irreducible $\delta$-polynomial in $\Q(t)\{y_1\}$ where
$\delta = {\partial \over \partial t}$; but it can be factored in
$\Q(\sqrt{t})\{y_1\}$.
Thus, to emphasize the base field, in the rest of the paper, we will
use $\mathcal {I}=[\Sigma]\subset \mathcal {G}\{\Y\}$ to denote the
$\delta$-ideal generated by $\Sigma$ in $\mathcal {G}\{\Y\}$, and
use $\mathcal {G}_1\cdot\mathcal {I}$ to denote the $\delta$-ideal
generated by $\CI$ in $\mathcal {G}_1\{\Y\}$ where $\mathcal
{G}_1\subset \ee$ is a $\delta$-extension field of $\mathcal {G}$.

\subsection{Characteristic sets of a  differential polynomial set}

     Let $f$ be a $\delta$-polynomial in $\ff\{\Y\}$.
     We define the order of $f$ w.r.t. $y_i$ to be the greatest number $k$ such that $y_{i}^{(k)}$
appears effectively in $f$, which is denoted by $\ord(f,y_{i})$. And
if $y_{i}$ does not appear in $f$, then we set
$\ord(f,y_{i})=-\infty$.
     The {\em order} of $f$ is defined to be $\max_{i}\,\ord(f,y_{i})$, that is,
     $\ord(f)=\max_{i}\,\ord(f,y_{i})$.

     A {\em ranking} $\mathscr{R}$ is a total order over $\Theta (\Y)$, which is compatible with
     the derivations over the alphabet:

     1) $\delta \theta y_{j} >\theta y_{j}$ for all derivatives $\theta y_{j}\in\Theta (\Y)$.

     2) $\theta_{1} y_{i} >\theta_{2} y_{j}$ $\Longrightarrow$ $\delta\theta_{1} y_{i} >\delta\theta_{2} y_{j}$
for $\theta_{1} y_{i}, \theta_{2} y_{j}\in \Theta (\Y)$.

     By convention, $1<\theta y_{j}$ for all $\theta y_{j}\in \Theta (\Y)$.

    Two important kinds of rankings are the following:

    1) {\em Elimination ranking}: \ $y_{i} > y_{j}$ $\Longrightarrow$ $\delta^{k}  y_{i} >\delta^{l}
    y_{j}$ for any $k,    l\geq 0$.

    2) {\em Orderly ranking}: \
     $k>l$  $\Longrightarrow$  $\delta^{k}  y_{i} >\delta^{l}
    y_{j}$,    for any $i,    j \in \{1,    2,    \ldots,    n\}$.

    Let $p$ be a $\delta$-polynomial in
    $\mathcal {F}\{\Y\}$ and $\mathscr{R}$  a ranking
    endowed on it.  The greatest derivative w.r.t.  $\mathscr{R}$ which  appears effectively in $p$ is called the {\em leader} of $p$,
    which will be denoted by $u_{p}$ or $\lead(p)$.  The
    two conditions mentioned above imply that the leader of $\theta
    p$ is $\theta u_{p}$ for $\theta\in\Theta$.  Let the degree of $p$  in $u_{p}$ be $d$.     We rewrite
    $p$
    as an algebraic polynomial in $u_{p}$.  Then
     $$p=I_{d} u_{p}^{d}+I_{d-1}u_{p}^{d-1}+\cdots+I_{0}.$$
    We call $I_{d}$ the {\em initial} of $p$ and denote it by $\init_{p}$.  The
    partial derivative of $p$ w.r.t. $u_{p}$ is called the {\em separant} of $p$,
    which will be denoted by $\sep_{p}$.  Clearly,    $\sep_{p}$
is the initial of any   proper derivative of $p$.
The {\em rank} of $p$ is $u_{p}^{d}$, and we denote it by $\rk(p)$.
For any two
    $\delta$-polynomials $p$, $q$ in $\mathcal {F}\{\Y\}\backslash \mathcal {F}$,
    $p$  is said to be of {\em lower rank} than  $q$ if either  $u_{p}<u_{q}$
    or  $u_{p}=u_{q}=u$ and $\deg(p,u)<\deg(q,u)$.
    By convention, any element of $\mathcal {F}$ is of lower rank than elements of $\mathcal{F}\{\Y\}\backslash\mathcal
    {F}$.  We denote $p \preceq
    q$ if and only if either $p$ is of lower rank than $q$ or they have
 the same rank.  Clearly,    $\preceq$ is a totally ordering of
$\mathcal{F}\{\Y\}$.

    Let $p$ and $q$ be two $\delta$-polynomials and $u_{p}^{d}$  the rank of $p$.  $q$ is said to be
    {\em partially reduced} w.r.t. $p$ if no proper derivatives of $u_{p}$ appear
    in $q$.  $q$ is said to be {\em reduced} w.r.t. $p$ if $q$ is partially reduced
    w.r.t. $p$ and $\deg(q,u_{p})<d$.  Let $\mathcal {A}$ be a set of
    $\delta$-polynomials.    $\mathcal {A}$ is said to be an
    {\em auto-reduced set} if each $\delta$-polynomial of $\mathcal {A}$ is reduced
    w.r.t.  any other element of $\mathcal{A}$.
   Every auto-reduced set is  finite.

Let $\mathcal {A}$ be an auto-reduced set. We
    denote $\H_{\mathcal {A}}$ to be  the set of all the initials and
    separants of $\mathcal {A}$ and $\H_{\mathcal {A}}^\infty$ to be the minimal
    multiplicative set containing $\H_{\mathcal {A}}$.
    The {\em saturation ideal} of $\A$ is defined to be
    $$\sat(\A)=[\mathcal
   {A}]:H_{\mathcal {A}}^\infty = \{p\in\ff\{\Y\}\big| \exists h\in H_{\mathcal
{A}}^\infty, \,{\text s. t. }\, hp\in[A]\}.$$

    Let $\mathcal {A}=A_{1},A_{2},\ldots,A_{s}$ and $\mathcal
    {B}=B_{1},B_{2},\ldots,B_{l}$  be two auto-reduced
    sets with the $A_{i}$, $B_{j}$ arranged in nondecreasing ordering.
    $\mathcal {A}$ is said to be of {\em lower rank} than $\mathcal {B}$, if
    either 1)\ there is some $k$ ($\leq  \min\{s,l\}$) such that for each
    $i<k$, $A_{i}$ has the same rank as $B_{i}$, and $A_{k}\prec B_{k}$
    or     2)\ $s>l$ and for each $i\in \{1, 2, \ldots, l\}$, $A_{i}$ has the same rank as
    $B_{i}$.
     It is easy to see that the above definition introduces really a
     partial ordering among all auto-reduced sets.
Any sequence of auto-reduced sets steadily decreasing in ordering
$\mathcal {A}_{1}\succ \mathcal {A}_{2}\succ \cdots\mathcal
{A}_{k}\succ \cdots$ is necessarily finite.

Let $\mathcal {A}=A_{1},A_{2},\ldots,A_{t}$ be an auto-reduced
   set with $\sep_{i}$ and $\init_{i}$ as the separant and initial of $A_{i}$, and $f$ any $\delta$-polynomial. Then there exists an
   algorithm, called Ritt's algorithm of reduction, which reduces
   $f$ w.r.t. $\mathcal {A}$ to a  $\delta$-polynomial $r$ that is
   reduced w.r.t. $\mathcal {A}$, satisfying the relation
   $$\prod_{i=1}^t\sep_{i}^{d_{i}}\init_{i}^{e_{i}} \cdot f \equiv
   r, \mod \,[\mathcal {A}],$$ for nonnegative integers
   $d_{i},e_{i}\, (i=1,2,\ldots,t)$. We call $r$ the {\em pseudo remainder} of $f$
w.r.t. $\A$.

 An auto-reduced set $\mathcal {C}$ contained in a $\delta$-polynomial set
 $\mathcal {S}$ is said to be a {\em characteristic set} of $\mathcal {S}$,
 if  $\mathcal {S}$ does not contain any nonzero element reduced w.r.t.
$\mathcal {C}$. All the characteristic  sets of $\mathcal {S}$ have
the same and minimal rank among all auto-reduced sets contained in
$\mathcal {S}$. A characteristic set $\mathcal{C}$ of a
$\delta$-ideal $\mathcal {J}$ reduces to zero all elements of
$\mathcal {J}$. If the $\delta$-ideal is prime, $\mathcal {C}$
reduces to zero only the elements of $\mathcal {J}$ and we have
$\mathcal {J}=\sat(\C)$ (\cite[Lemma 2, p.167]{kol}).

In polynomial algebra, let $\mathcal {A}=A_{1},A_{2},\ldots,A_{t}$
be an algebraic auto-reduced set arranged in nondecreasing order.
$\mathcal{A}$  is called an {\em irreducible  auto-reduced set} if
for any $1\leq i\leq t$, there can not exist any relation of the
form
 $$T_iA_i=B_iC_i,\, \mod \,(A_1,\ldots,A_{i-1})$$
where $B_i,C_i$ are polynomials with the same leader as $A_i$, $T_i$
is a polynomial with lower leader than $A_i$, and $B_i, C_i, T_i$
are reduced w.r.t. $A_1,\ldots,A_{i-1}$ (\cite{wu2}). Equivalently,
an algebraic auto-reduced set $\mathcal {A}$ is irreducible if and
only if  there exist no polynomials $P$ and $Q$ which are reduced
w.r.t $\A$ and $PQ\in\asat(\A)=(\mathcal
{A}):I_\mathcal{A}^\infty$, where $I_\mathcal{A}^\infty$ stands for
the set of all products of powers of $\init_{A_i}$.

In ordinary differential algebra, we can define an auto-reduced set
to be {\em irreducible} if when considered as an algebraic
auto-reduced set in the underlying polynomial ring, it is
irreducible. We have (\cite[p.107]{ritt})

\begin{theorem}\label{th-sat} Let $\mathcal {A}$  be an auto-reduced set. Then
a necessary and sufficient condition for $\mathcal {A}$ to be a
characteristic set of a prime $\delta$-ideal is that $\mathcal {A}$
is irreducible. Moreover, in the case $\mathcal {A}$ is irreducible,
$\sat(\mathcal {A})$=$[\mathcal {A}]:\hbox{\rm{H}}_{\mathcal
{A}}^{\infty}$ is prime with $\mathcal {A}$ being a characteristic
set of it.
\end{theorem}

\begin{remark} \label{re-1}
A set of $\delta$-polynomials $\A=\{A_1\ldots,A_p\}$ is called a
{\em $\delta$-chain} if the following conditions are satisfied,

1) the leaders of $A_i$ are $\delta$-auto-reduced,

2) each $A_i$ is partially reduced w.r.t. all the others,

3) no initial of an element of $\A$ is reduced to zero by $\A$.

Similar properties to auto-reduced sets can be developed for
$\delta$-chains \cite{hubert}. In particular, we can define a
$\delta$-characteristic set of a $\delta$-ideal $\mathcal{I}$ to be
a $\delta$-chain contained in $\mathcal{I}$ of minimal rank among
all the $\delta$-chains contained in $\mathcal {I}$.
So, in this paper we will not distinguish auto-reduced sets and
$\delta$-chains. Note that we can also use the weak $\delta$-chains
introduced in \cite{ardm1}.
\end{remark}

\subsection{Dimension and order of a prime differential polynomial ideal}
\label{sec-ord} Let $\mathcal {I}$ be a prime $\delta$-ideal in
$\ff\{\Y\}$ and $\xi=(\xi_{1},\ldots,\xi_{n})$ a generic point of
$\mathcal {I}$ \cite[p.19]{kol}.
The {\em dimension} of $\mathcal {I}$ or  $\dzero(\mathcal {I})$ is
defined to be the $\delta$-transcendence degree of $\mathcal
{F}\langle \xi_{1},\ldots,\xi_{n}\rangle$ over $\mathcal {F}$, that
is, $\dim(\mathcal {I})=\dtrdeg\, \mathcal {F}\langle
\xi_{1},\ldots,\xi_{n}\rangle/\mathcal {F}$.

In \cite{ritt}, Ritt gave another definition of the dimension of
$\mathcal{I}$. An independent set modulo $\mathcal {I}$ is defined
to be a variable set $\U\subset\Y$ such that $\mathcal {I}\cap
\mathcal {F}\{\U\}=\{0\}$, and in this case $\U$ is also said to be
$\delta$-independent modulo $\mathcal {I}$. And a {\em parametric
set} of $\mathcal {I}$ is a maximal independent set modulo $\mathcal
{I}$. Then Ritt defined the dimension of $\mathcal{I}$ to be the
cardinal number of its parametric set. Clearly, the two definitions
are equivalent.
%
%

\begin{definition}\cite{kol3}
Let  $\mathcal {I}$ be a prime $\delta$-ideal of $\mathcal
{F}\{\Y\}$ with a generic point $\eta=(\eta_{1},\ldots,\eta_{n})$.
Then there exists a unique numerical polynomial $\omega_{\mathcal
{I}}(t)$  such that $\omega_{\mathcal {I}}(t)=\omega_{\eta/\ff}(t)=\trdeg\,\mathcal
{F}(\eta_{i}^{(j)}:i=1,\ldots,n;j\leq t)/\mathcal {F} $ for all
sufficiently large $t \in \mathbb{N}$. $\omega_{\mathcal {I}}(t)$ is
called the {\em $\delta$-dimension polynomial} of $\mathcal {I}$.
\end{definition}

We now define the order of a prime $\delta$-ideal $\CI$, which is
also related with the characteristic set of $\CI$.
\begin{definition} \label{def-order-char}  For an auto-reduced set $\mathcal
{A}=A_{1}, A_{2}, \ldots, A_{t}$,  with
$\lead(A_{i})=y_{c_i}^{(o_{i})}$, the order of $\mathcal {A}$ is
defined to be $\ord(\A) = \sum\limits_{i=1}^{t} o_{i}$, and the set
$\Y\backslash\{y_{c_1},\ldots,y_{c_t}\}$ is called a parametric set
of $\mathcal {A}$.
\end{definition}

\begin{theorem} \cite[Theorem 13]{sadik} \label{def-order-sadik}
Let $\mathcal {I}$ be a prime $\delta$-ideal of dimension $d$. Then
the $\delta$-dimension polynomial has the form $\omega_{\mathcal
{I}}(t)=d(t+1)+h$, where $h$ is defined to be the order of $\mathcal
{I}$ or $\dzero(\mathcal {I})$, that is, $\ord(\mathcal {I})= h$.
Let $\mathcal{A}$ be a characteristic set of $\mathcal {I}$ under
any orderly ranking. Then, $\ord(\mathcal {I})= \ord(\A)$.
\end{theorem}

In \cite{ritt}, Ritt introduced the concept of relative order for a
prime $\delta$-ideal w.r.t. a particular parametric set.

\begin{definition}\label{def-rord}
  Let $\mathcal {I}$ be a prime $\delta$-ideal of $\mathcal
  {F}\{\Y\},    \mathcal{A}$  a characteristic
  set of $\mathcal {I}$ w.r.t. any elimination ranking,    and
  $\{u_{1},    \ldots,    u_{d}\}\subset\Y$   the parametric set of $\mathcal
  {A}$. The {\em relative order} of $\mathcal {I}$  w.r.t.
  $\{u_{1}, \ldots, u_{d}\}$, denoted by $\ord_{u_{1},\ldots,u_{d}}\mathcal {I}$,
  is defined to be $\ord(\A)$.
\end{definition}

The relative order of a prime $\delta$-ideal $\mathcal{I}$ can be
computed from its generic points as shown by the following result
(\cite{kol5}).
\begin{cor} \cite{kol5}\label{cor-order1}
Let $\mathcal{I}$ be a prime $\delta$-ideal in $\mathcal {F}\{\Y\}$
with a generic point $(\xi_1,\ldots,\xi_n)$. If $\{y_1,\ldots,y_d\}$
is a parametric set of $\mathcal{I}$, then
$\ord_{y_1,\ldots,y_d}(\mathcal{I}) =
\trdeg\,\mathcal{F}\langle\xi_1,$
$\ldots,\xi_d\rangle\langle\xi_{d+1},$ $\ldots,\xi_n\rangle/\mathcal
{F}\langle\xi_1,\ldots,$ $\xi_d\rangle$.
\end{cor}

Ritt's definition of relative order is based on the elimination
ranking. Hubert proved that all characteristic sets of $\mathcal
{I}$ admitting the same parametric set have the same order
\cite{hubert}.

\begin{theorem} \cite{hubert} \label{relativehubert}
Let $\mathcal {A}$ be a characteristic set of a prime $\delta$-ideal
$\mathcal {I}$ in $\mathcal {F}\{\Y\}$ endowed with any ranking. The
parametric set $\U$ of $\mathcal {A}$ is a maximal independent set
modulo $\mathcal {I}$. Its cardinality gives the dimension of
$\mathcal {I}$. Furthermore, the order of $\mathcal{I}$ relative to
$\U$ is the order of $\mathcal{A}$.
\end{theorem}

\begin{cor} \label{cor-0 order} Let  $\mathcal {I}$ be a prime $\delta$-ideal
with dimension zero,    and $\mathcal {A}$ a characteristic set of
$\mathcal{I}$ w.r.t. any ranking $\mathscr{R}$. Then $\ord
(\mathcal{I})=\ord (\mathcal {A})$.
\end{cor}

The following result gives the relation between the order and
relative order for a prime $\delta$-ideal.
\begin{theorem} \label{th-maxorder}
Let $\mathcal {I}$ be a prime $\delta$-ideal in $\mathcal
{F}\{\Y\}$. Then $\ord(\mathcal{I})$ is the maximum of all the
relative  orders of $\mathcal {I}$, that is, $\ord(\mathcal
{I})=\max_{\U }\, \ord_{\U}(\mathcal {I})$, where $\U$ is any
parametric set of $\mathcal {I}$.
\end{theorem}
\proof Let $\mathcal {C}$ be a characteristic set of $\mathcal {I}$
w.r.t. some orderly ranking.
Firstly, we claim that any relative order of $\mathcal {I}$ is less
than or equal to $\ord(\mathcal {C})$. Let
$\U=\{u_{1},\ldots,u_{q}\}$ be any parametric set of $\mathcal {I}$,
$\{y_{1},\ldots,y_{p}\} (p+q=n)$ the set of the remaining variables,
and $\mathcal {B}$ any characteristic set of $\mathcal {I}$ w.r.t.
the elimination ranking $u_{1} \prec \ldots \prec u_{q} \prec y_{1}
\prec \ldots \prec y_{p}$. By Theorem \ref{relativehubert}, it
suffices to prove
 $\ord_{\U}(\mathcal {I}) \leq \ord(\mathcal {C})$.

Let $\eta=(\overline{u_{1}},    \ldots,    \overline{u_{q}},
\overline{y_{1}}, \ldots,    \overline{y_{p}})$ be a generic point
of $\mathcal {I}$. Then for sufficiently large  $t$,    the
$\delta$-dimension polynomial of $\mathcal {I}$ is
 \begin{eqnarray}
&& \omega_{\mathcal {I}} (t) \nonumber \\
 & = & \omega_{\eta / \mathcal {F}} (t) \nonumber \\
                & = &  \trdeg\, \mathcal {F}(\delta ^{s} \overline{u_{i}},    \delta ^{k} \overline{y_{j}}: s,    k \leq t;i=1,\ldots,q;j=1,\ldots,p)/\mathcal
 {F} \nonumber \\
                & = &  \trdeg\, \mathcal {F}(\delta ^{s} \overline{u_{i}}:    s \leq t)(\delta
^{k} \overline{y_{j}}:    k \leq t)/\mathcal {F}(\delta ^{s}
\overline{u_{i}}: s \leq t)+\trdeg\, \mathcal {F}(\delta ^{s}
\overline{u_{i}}:    s \leq
t)/\mathcal {F} \nonumber \\
                & = & q(t+1)+\trdeg\,
\mathcal {F}(\delta ^{s} \overline{u_{i}}:   s \leq t)(\delta ^{k}
\overline{y_{j}}:   k \leq t)/\mathcal {F}(\delta ^{s}
\overline{u_{i}}:   s \leq  t) \nonumber
\end{eqnarray}

Since $\omega_{\mathcal {I}} (t)=q(t+1)+\ord(\mathcal {C})$,
$\trdeg\, \mathcal {F}(\delta ^{s} \overline{u_{i}}:   s \leq
t)(\delta ^{k} \overline{y_{j}}:   k \leq t)/\mathcal {F}(\delta
^{s} \overline{u_{i}}: s \leq  t)=\ord(\mathcal {C})$. By Corollary
\ref{cor-order1}, we have
 \begin{eqnarray}
 && \ord_{\U}(\mathcal {I})  \nonumber \\
 & = &  \trdeg\,  \mathcal {F}\langle \overline{u_{1}},    \ldots,    \overline{u_{q}}\rangle (\delta ^{k} \overline{y_{l}}:  k \geq 0)/\mathcal{F}\langle \overline{u_{1}},    \ldots,    \overline{u_{q}} \rangle  \nonumber \\
                & = &  \trdeg\, \mathcal {F}\langle \overline{u_{1}},    \ldots,    \overline{u_{q}}\rangle (\delta ^{k} \overline{y_{l}}:    k \leq t)/\mathcal{F}\langle \overline{u_{1}},    \ldots,    \overline{u_{q}} \rangle\quad (\hbox{for }t \geq \ord(\mathcal {B})) \nonumber \\
                & \leq &  \trdeg\, \mathcal {F}(\delta ^{s} \overline{u_{i}}:    s \leq t)(\delta ^{k}
\overline{y_{j}}:    k \leq t)/\mathcal {F}(\delta ^{s}
\overline{u_{i}}:   s \leq
 t) \nonumber \\
                & = &  \ord(\mathcal {C}). \nonumber
\end{eqnarray}
Thus, the claim is proved.

Now, let $\U^*$ be the parametric set of $\mathcal {C}$. Then, by
Theorem~\ref{relativehubert},
$\ord(\mathcal{I})=\ord(C)=\ord_{\U^*}(\mathcal{I})$.
That is, for any parametric set $\U$ of $\mathcal{I}$, we have
$\ord_\U(\mathcal{I}) \le \ord(\mathcal{I})$ and there exists one
parametric set $\U^*$ of $\mathcal{I}$ such that
$\ord_{\U^*}(\mathcal{I}) = \ord(\mathcal{I})$.
As a consequence, $\ord(\mathcal {I}) =\max_{\U} \ord_{\U}(\mathcal
{I})$. \qed

The following well known result about the adjoining indeterminates
to the base field will be used in this paper \cite{kol,ritt}.
\begin{lemma}\label{le-fieldextension}
Let $\U=\{u_1,\ldots,u_r\}\subset \ee$ be a set of
$\delta$-$\ff$-indeterminates,
 $\mathcal {I}_0$ a prime $\delta$-ideal of dimension $d$
and order $h$ in $\mathcal {F}\{\Y\}$, and $\mathcal {I}=\mathcal
{F}\langle \U\rangle\cdot\mathcal {I}_0$ the $\delta$-ideal
generated by $\mathcal {I}_0$ in $\mathcal {F}\langle
\U\rangle\{\Y\}$. Then $\mathcal {I}$ is a prime $\delta$-ideal in
$\mathcal {F}\langle \U\rangle\{\Y\}$ of dimension $d$ and order $h$
and $\mathcal {I}\cap\mathcal {F}\{\Y\}=\mathcal {I}_0$.
\end{lemma}
%

\subsection{A property on differential specialization } \label{sec-spec}
The following lemma is a key result in algebraic elimination theory,
which is used to develop the theory of Chow form (\cite[p.
168-169]{hodge1}, \cite[p. 161]{wu}). The result is originally given
for homogenous polynomials. We will show that it also holds for
non-homogenous polynomials.

\begin{lemma}\label{lem-algspe}
Let $P_{i} \in \mathcal{F}[\U, \Y]\,(i=1, \ldots, m)$ be polynomials
in the independent indeterminates $\U=(u_{1},  \ldots, u_{r})$ and
$\Y=(y_{1}, \ldots,  y_{n})$. Let $\overline{\Y}=(\overline{y}_{1},
\ldots, \overline{y}_{n})$, where $\overline{y}_{i}$ are elements of
some extension field of $\mathcal{F}$ free from
$\ff(\U)$\footnote{By saying $\overline{\Y}$ free from
$\ff(\U)$\,(resp. $\delta$-free from $\ff\langle\U\rangle$), we mean
that $\U$ is a set of indeterminates over
$\ff(\overline{\Y})$\,(resp. $\delta$-indeterminates over
$\ff\langle\overline{\Y}\rangle$).}. If $P_{i}(\U,
\overline{\Y})\,(i=1, \ldots, m)$ are algebraically dependent over
$\mathcal{F}(\U)$, then for any specialization $\U$ to
$\overline{\U}=(\overline{u}_1,\ldots,\overline{u}_r)\in\mathcal{F}^r$,
$P_{i}(\overline{\U}, \overline{\Y})\,(i=1,    \ldots, m)$ are
algebraically dependent over $\mathcal{F}$.
\end{lemma}
\proof We sketch the proof for the case $r=1$. Since $P_{i}(u_1,
\overline{\Y})\,(i=1, \ldots,m)$  are algebraically dependent over
$\mathcal{F}(u_1)$, there exists some nonzero polynomial
$f\in\ff(u_1)[z_1,\ldots,z_m]$ which vanishes for $z_i=P_{i}(u_1,
\overline{\Y}).$
 By clearing of fractions when necessary, we suppose $f\in\ff[u_1,z_1,\ldots,z_m]$.
Now specialize $u_1$ to $\overline{u}_1$ in $f$, then we have
$f(\overline{u}_1;P_{1}(\overline{u}_1,
\overline{\Y}),\ldots,P_{m}(\overline{u}_1, \overline{\Y}))=0.$ If
$f(\overline{u}_1;z_1,\ldots,z_m)\neq0$, it follows that
$P_{i}(\overline{u}_1, \overline{\Y})\,(i=1,    \ldots, m)$ are
algebraically dependent over $\mathcal{F}$. If
$f(\overline{u}_1;z_1,\ldots,z_m)=0$, then
$f(u_1;z_1,\ldots,z_m)=(u_1-\overline{u}_1)^lf_1$ where
$f_1(\overline{u}_1;z_1,\ldots,z_m)\neq0$. Since
$f(\overline{u}_1;P_{1}(\overline{u}_1,
\overline{\Y}),\ldots,P_{m}(\overline{u}_1, \overline{\Y}))=0$, we
have $f_1(\overline{u}_1;P_{1}(\overline{u}_1,
\overline{\Y}),\ldots,P_{m}(\overline{u}_1, \overline{\Y}))=0$.
Thus, it follows that $P_{i}(\overline{\U}, \overline{\Y})\,(i=1,
\ldots, m)$ are algebraically dependent over $\mathcal{F}$. \qed

To generalize the above  result to the differential case,  we need
the following lemma \cite[p35]{ritt}.

%

\begin{lemma} \label{le-notannul}  Suppose $\mathcal {F}$ contains at least one nonconstant
element. If $G \in \mathcal {F}\{u\}$ is a nonzero
$\delta$-polynomial with order $r$,    then for any nonconstant
$\eta \in \mathcal {F}$, there exists an element
$c_{0}+c_{1}\eta+c_{2}\eta^2+\cdots+c_{r}\eta^{r}$ which does not
annul $G$,    where $c_{0},    \ldots,    c_{r}$ are constants in
$\mathcal {F}$.
\end{lemma}

By {\em arbitrary constants} $a_1,\ldots,a_s$ over $\ff$, we mean
$a_1,\ldots,a_s$ are constants in $\ee$, which are algebraically
independent over $\ff$. As a consequence of Lemma \ref{le-notannul},
we have

\begin{cor} \label{cor-notannul}
Suppose $\mathcal {F}$ contains at least one nonconstant element. If
$G \in \mathcal {F}\{u\}$ is a nonzero $\delta$-polynomial with
order $r$,   then for any nonconstant $\eta \in \mathcal {F}$ and
arbitrary constants $a_0,\ldots,a_r$ over $\ff$,
$a_{0}+a_{1}\eta+a_{2}\eta^2+\cdots+a_{r}\eta^{r}$ does not annul
$G$.
\end{cor}
\proof Suppose the contrary. Let $a_0,\ldots,a_r$ be arbitrary
constants over $\ff$ in $\ee$ such that
$G(a_{0}+a_{1}\eta+\cdots+a_{r}\eta^{r})=0$, where $\eta$ is a
nonconstant in $\ff$. Since $a_0,\ldots,a_r$ are constants in $\ee$
which are algebraically independent over $\ff$,
$g(a_0,\ldots,a_r)=G(a_{0}+a_{1}\eta+\cdots+a_{r}\eta^{r})$ is a
polynomial in $\ff[a_0,\ldots,a_r]$. Now by the hypothesis, $g$ is a
zero polynomial. Thus, for any constants
$\overline{a_i}\in\ff\,(i=0,\ldots,r)$,
$g(\overline{a_0},\ldots,\overline{a_r})=G(\overline{a_0}+\overline{a_{1}}\eta+\cdots+\overline{a_{r}}\eta^{r})=0$,
which contradicts  Lemma \ref{le-notannul}.\qed

\vskip5pt Now we prove the following result, which is crucial
throughout the paper.
\begin{theorem}\label{th-specil}
Let $\{u_{1}, \ldots,  u_{r}\}\subset \ee$ be a set of
$\delta$-$\ff$-indeterminates,    and $P_{i}(\U, \Y)$ $\in \mathcal
{F}\{\U, \Y\}$ $(i=1, \ldots, m)$  $\delta$-polynomials in the
$\delta$-$\ff$-indeterminates $\U=(u_{1},\ldots,u_{r})$ and
$\Y=(y_{1}, \ldots, y_{n})$. Let $\overline{\Y}=(\overline{y}_{1},
\overline{y}_{2}, \ldots, \overline{y}_{n})$,  where
$\overline{y}_{i}\in \ee$ are $\delta$-free from $\ff\langle \U
\rangle$. If $P_{i}(\U, \overline{\Y})$ $(i=1, \ldots, m)$ are
$\delta$-dependent over $\mathcal {F}\langle \U \rangle$, then for
any specialization $\U$ to $\overline{\U}$ in $\mathcal {F}$,
$P_{i}(\overline{\U}, \overline{\Y}) \, (i=1, \ldots,  m)$ are
$\delta$-$\mathcal {F}$-dependent.
\end{theorem}
\proof  It suffices to prove the case $r=1$. Denote $u_1$ by $u$.
Firstly, we suppose $\mathcal {F}$ contains at least one nonconstant
element.

Since $P_{i}(u, \overline{\Y}) \, (i=1,    \ldots,    m)$ are
$\delta$-dependent over $\mathcal {F}\langle u \rangle$, there
exists a nonzero $G(z_{1},    \ldots,    z_{m}) \in \mathcal
{F}\langle u \rangle \{z_{1},    \ldots,    z_{m}\}$ such that
$G(P_{1}(u, \overline{\Y}), \ldots,    P_{m}(u, \overline{\Y}))$
$=0$. We can take $G \in \mathcal {F} \{u, z_{1}, \ldots, z_{m}\}$
by clearing denominators when necessary.

Since $G(u, z_{1}, \ldots, z_{m}) \neq 0$,  by
Corollary~\ref{cor-notannul}, for any nonconstant $\eta \in \mathcal
{F}$ and arbitrary constants $ c_{0}, \ldots,c_{s}$ $(s=\ord(G, u))$
over $\mathcal {F}\langle \overline{\Y}\rangle$,
$G(u^*,z_1,\ldots,z_m)\neq 0$ where
$u^*=\sum_{i=0}^{s}c_{i}\eta^{i}$.  Since $G(u,P_{1}(u,
\overline{\Y}), \ldots, P_{m}(u, \overline{\Y}))=0$ and $u$ is a
$\delta$-indeterminate over $\ff\langle \overline{\Y}\rangle$, when
$u$ is specialized to $u^*$, we have $G(u^*, P_{1}(u^*,
\overline{\Y}), \ldots,$ $P_{m}(u^*, \overline{\Y})) =0$. Regarding
$G$ as an algebraic polynomial over $\ff$  in $c_{i} \, (i=0,
\ldots, s)$ and $z_{i}^{(j)}\,(i=1,\ldots,m;j\geq0)$ which appear
effectively, we have
$$G(c_{0},\ldots, c_{s}, \ldots, z_{i}^{(j)},
\ldots)\neq 0$$ while $$G(c_{0},    \ldots, c_{s}, \ldots,
(P_{i}(u^*, \overline{\Y}))^{(j)},   \ldots)= 0.$$ So $P_{i}(u^*,
\overline{\Y}))^{(j)}\,(i=1,\ldots,m;j\geq0)$ are algebraically
dependent over $\mathcal {F}(c_0,\ldots,c_s)$, by
Lemma~\ref{lem-algspe}, when the $c_{i}$ are specialized to
constants $\overline{c}_{i}$ in $\mathcal {F}$, the corresponding
$P_{i}(\overline{u^*}, \overline{\Y})^{(j)}\,(i=1,\ldots,m)$ are
algebraically dependent over $\mathcal {F}$, where
$\overline{u^*}=\sum_{i=0}^{s}\overline{c}_{i}\eta^{i}$. That is,
$P_{i}(\overline{u^*}, \overline{\Y}) \,(i=1, \ldots, m)$ are
$\delta$-dependent over $\mathcal {F}$.

To complete the proof, if $\overline{u}$ is a nonconstant,    as
above we take $\eta=\overline{u}$,    and specialize
$c_{1}\rightarrow 1$ and other $c_{i}$ to zero; else we take $\eta$
as an arbitrary nonconstant and specialize $c_{0} \rightarrow
\overline{u}$ and other $c_{i}$ to zero. Then in either case, $u$ is
specialized to $\overline{u}$, and we have completed the proof in
the case that $\mathcal {F}$ contains at least one nonconstant
element.

 If $\mathcal {F}$ consists of constant elements, take $v\in\ee$ to be a
 $\delta$-indeterminate over $\mathcal {F}\langle \U,\overline{\Y}\rangle$. Now we consider in
 the $\delta$-field $\mathcal {F}\langle v \rangle$.    Following the first
 case,   for any specialization $\U$ to $\overline{\U}\subset \mathcal {F}$, we can show that $P_{i}(\overline{\U},  \overline{\Y}) \, (i=1,    \ldots,    m)$ are $\delta$-dependent over $\mathcal {F}\langle v\rangle$. Since $v$ is a
$\delta$-indeterminate over $\mathcal {F}\langle
\overline{\Y}\rangle$, $P_{i}(\overline{\U}, \overline{\Y})\,(i=1,
\ldots, m)$ are $\delta$-dependent over $\mathcal {F}$. \qed

From the proof above, we can obtain the following result easily:
 \begin{cor} \label{corspecial}
 Let $\{u_{1}, \ldots, u_{r}\}\subset\ee$ be a set of $\delta$-$\ff$-indeterminates, and $P_{i}(\U, \Y)$ $\in \mathcal {F}\{\U,
\Y\}\,(i=1, \ldots, m)$  $\delta$-polynomials in the
$\delta$-indeterminates $\U=(u_{1},\ldots,u_{r})$ and $\Y=(y_{1},
\ldots, y_{n})$.  Let $\overline{\Y}=(\overline{y}_{1},
\overline{y}_{2}, \ldots, \overline{y}_{n})$, where
$\overline{y}_{i}\in \mathcal {E}$ are free from $\ff\langle
\U\rangle$. If the set $(P_{i}(\U,
\overline{\Y}))^{(\sigma_{ij})}\,(i=1, \ldots, m; j=1,\ldots,n_{i})$
are algebraically dependent over $\mathcal {F}\langle \U \rangle$,
then for  any specialization $\U$ to $\overline{\U}$ in $\mathcal
{F}$, $(P_{i}(\overline{\U}, \overline{\Y}))^{(\sigma_{ij})}(i=1,
\ldots, m; j=1,\ldots,n_{i})$ are algebraically dependent over
$\mathcal {F}$.
\end{cor}

Now we give an example to illustrate the proof of
Theorem~\ref{th-specil}.

\begin{example}
In this example $r=m=n=1$ and $\ff=\mathbb{Q}(t)$ with
$\delta=\frac{d}{dt}$. Let $P(u,y)=uy$. Suppose that $\xi$ is a
generic point of $\mathcal
{I}=\sat(y'^2-2y)\subset\mathbb{Q}(t)\{y\}$ that is $\delta$-free
from $\mathbb{Q}(t)\langle u\rangle$. Let
$G(u;z)=u'^2z''^2+2u(u'^2-uu'')z''-2u'u''z'z''+u''^2z'^2+z'(-8u'^3+6uu'u'')-3u^2u'^2+2u^3u''\in\mathbb{Q}(t)\{u,z\}$.
Let $\overline{y}=\xi$. It is easy to verify that $G(u;z)\neq0$ and
$G(u;P(u,\overline{y}))=0$. That is, $P(u,\overline{y})$ is
$\delta$-algebraic over $\mathbb{Q}(t)\langle u\rangle$. We now
proceed to show that when $u$ is specialized to $\overline{u}=1$,
$P(\overline{u},\overline{y})$ is $\delta$-algebraic over
$\mathbb{Q}(t)$.

Clearly, $G(\overline{u};z)=0$. In the algebraic case, following the
proof of Lemma~\ref{lem-algspe}, $u-\overline{u}$ must be a factor
of $G(u;z)$. Removing all factors of the form $u-\overline{u}$, we
obtain a new $\delta$-polynomial $G_1(u;z)$ which still satisfies
$G_1(u;P(u,\overline{y}))=0$ and $G_1(\overline{u};z)\ne0$. But, in
the differential case $G(u;z)$  may not contain any factor involving
$u-\overline{u}$ and its derivatives.

Follow the steps in Theorem~\ref{th-specil}, let $c_0,c_1,c_2$ be
arbitrary constants over $\mathbb{Q}(t)\langle \xi\rangle$ and
$\eta=t$. Denote $u^*=c_0+c_1t+c_2t^2$. Then $G(u^*;z)=-40 z'c_2^3
t^3-6c_0c_1^3 t-12c_0 c_2^3 t^4-14 c_1^3t^3 c_2-27 c_1^2 t^4 c_2^2+4
z'^2 c_2^2- 8 z'c_1^3-3 c_0^2 c_1^2-3 c_1^4 t^2-8 c_2^4 t^6+4 c_0^3
c_2+z''^2 c_1^2-36 z' c_1^2 c_2 t- 60 z' c_1 c_2^2 t^2+12 z'  c_0
c_1 c_2+24 z'  c_0 c_2^2 t-18 c_0 c_1^2 t^2 c_2-24 c_0 c_1 t^3
c_2^2+ 4 z''^2 c_1 c_2 t+6 z''c_1^2 t^2 c_2+8 z'' c_1 t^3 c_2^2-4
z'z'' c_1 c_2-8 z'  z'' c_2^2 t-24 c_1 t^5 c_2^3+2 z'' c_0 c_1^2-4
z'' c_0^2 c_2+ 4 z''^2 c_2^2 t^2+2 z'' c_1^3 t+4 z'' c_2^3 t^4\neq0$
and $G(u^*;P(u^*,\overline{y}))=0.$ Regard $G(u^*;z)$ as an
algebraic polynomial in $\mathbb{Q}(t)[c_0,c_1,c_2;z,z',z'']$,
denoted by $g(c_0,c_1,c_2;z,z',z'')$. Then
$P(u^*,\overline{y}),P(u^*,\overline{y})',P(u^*,\overline{y})''$ are
algebraically dependent over $\mathbb{Q}(t)(c_0,c_1,c_2)$. So the
problem is converted to an algebraic one, and we can use
Lemma~\ref{lem-algspe} to solve it.

To be more precise, firstly, specialize $c_2$ to 0, we obtain
$g(c_0,c_1,0;z,z',z'')=z''^2c_1^2+2z''c_0c_1^2+2z''c_1^3t-6c_0c_1^3t-8z'c_1^3-3c_0^2c_1^2-3c_1^4t^2\neq0$
and
$g(c_0,c_1,0;P(c_0+c_1t,\overline{y}),P(c_0+c_1t,\overline{y})',P(c_0+c_1t,\overline{y})'')=0$.
Then specialize $c_1$ to 0, we obtain $g(c_0,0,0;z,z',z'')=0$ and
$g(c_0,c_1,0;z,z',z'')=c_1^2g_1$ where
$g_1=z''^2+2z''c_0+2z''c_1t-6c_0c_1t-8z'c_1-3c_0^2-3c_1^2t^2$.
Clearly,
$g_1(c_0,c_1;P(c_0+c_1t,\overline{y}),P(c_0+c_1t,\overline{y})',P(c_0+c_1t,\overline{y})'')=0.$
Specialize $c_1$ to 0 in $g_1$,
$\overline{g}_1=z''^2+2z''c_0-3c_0^2$ while
$\overline{g}_1(c_0;P(c_0,\overline{y}),P(c_0,\overline{y})',P(c_0,\overline{y})'')=0$.
Now specialize $c_0$ to $\overline{u}=1$ in $\overline{g}_1$, we
obtain $g_2(z,z',z'')=z''^2+2z''-3$. Clearly,
$g_2(\xi,\xi',\xi'')=0$ where $\xi=P(\overline{u},\overline{y})$.
Thus, $P(\overline{u},\overline{y})$ is $\delta$-algebraic over
$\mathbb{Q}(t)$. Also,
$P(\overline{u},\overline{y}),P(\overline{u},\overline{y})'$ and
$P(\overline{u},\overline{y})''$ are algebraically dependent over
$\mathbb{Q}(t)$.
\end{example}

\section{Intersection theory for generic differential polynomials}
\label{sec-dim} In this section, we will develop an intersection
theory for generic $\delta$-polynomials by proving Theorem
\ref{th-main1}. As a consequence, the dimension conjecture is shown
to be true for generic $\delta$-polynomials. These results will also
be used in Sections 4 and 6 to determine the order  of the Chow
form.

\subsection{Generic dimension theorem}
In this section, we will show that the dimension conjecture is valid
for certain generic $\delta$-polynomials.
To prove the dimension conjecture in the general case, one simple idea
is to generalize the following theorem (\cite[p.43]{kol}) in algebra to the differential case.
\begin{theorem}\label{th-00}
Let $\CI$ be a prime ideal of dimension $d>0$  and $f \in
\mathcal{F}[\Y]$. If $(\CI,f) \neq (1)$. Then every prime component
of $\V(\CI,f)$ has dimension not less than $d-1$. Moreover, if $f$
is not in $\CI$,  then every prime component of $\V(\CI,f)$ has
dimension $d-1$.
\end{theorem}
Unfortunately, in the differential case,    the above theorem does
not hold. Ritt gave the following counter example.
\begin{example}\label{ex-dim1} \cite[p.133]{ritt}
$p=y_1^{5}-y_2^{5}+y_3(y_1y_2'-y_2y_1')^2 \in \mathcal
{F}\{y_1,y_2,y_3\}$ and $f=y_3$, where $\mathcal {F}$ is the  field
of complex numbers.  Then $\sat(F)$ is a prime $\delta$-ideal of
dimension two. But, $\sqrt{[\sat(p), f]} =[y_1,y_2,y_3]$ which is a
prime $\delta$-ideal of dimension zero.
%
\end{example}

It could also happen that when adding a $\delta$-polynomial to a
prime $\delta$-ideal which it does not belong to, the dimension is
still the same.
\begin{example}\label{ex-dim2}
Let $p=y_1'y_2''-y_1''y_2'$. Then $\sat(p)=[p]:y_1'^\infty$ is a
prime $\delta$-ideal of dimension one. It is clear that
$y_2'\not\in\sat(p)$ and $[\sat(p),y_2'] = [y_2']$ is still a prime
$\delta$-ideal of dimension one.
\end{example}

In this section, we will prove that Theorem \ref{th-00} is valid for
certain generic $\delta$-polynomials, which will lead to the
solution to the dimension conjecture in these generic cases.

 \begin{definition} \label{def-genericpol}
Let $\mathbbm{m}_{s,r}$ be the set of all $\delta$-monomials in
$\ff\{\Y\}$ of  order $\leq s$ and degree $\leq r$. Let
$\mathbb{U}=\{u_m\}_{m\in \mathbbm{m}_{s,r}}$ be a set of elements
of $\ee$ that are $\delta$-$\ff$-indeterminates. Then,
$$f=\sum_{m\in \mathbbm{m}_{s,r}}u_m m$$ is called a {\em generic $\delta$-polynomial} of order $s$ and degree $r$.
If $s=0$, then $f$ is also called a {\em generic polynomial} of
degree $r$.
A {\em generic $\delta$-hypersurface} is the set of zeros of a
generic $\delta$-polynomial.
\end{definition}
We use $\bu_f$ to denote the set of coefficients of a generic
$\delta$-polynomial $f$ and
 \begin{equation}\label{eq-buf}\widetilde{\bu}_f=\bu_f\setminus\{u_0\}\end{equation}
where $u_0$ is the degree zero term of $f$.
By saying that a point $\eta\in\ee^n$ is free from the pure
$\delta$-extension field $\ff\langle\bu_f\rangle$ over $\ff$, we
mean $\bu_f$ are $\delta$-$\ff\langle\eta\rangle$-indeterminates.

Throughout the paper, a generic $\delta$-polynomial is assumed to be
of degree greater than zero.

\begin{lemma}\label{lm-pu}
Let $\mathcal {I}$ be a prime $\delta$-ideal in $\mathcal {F}\{\Y\}$
with dimension $d$ and $f$ a generic $\delta$-polynomial. Then
$\mathcal{I}_0=[\mathcal {I}, f]$ is a prime $\delta$-ideal in
$\mathcal {F}\langle \widetilde{\bu}_f\rangle\{\Y,u_0\}$ with
dimension $d$, where $\widetilde{\bu}_f$ is defined in
\bref{eq-buf}. Furthermore, $\mathcal{I}_0\cap\mathcal {F}\langle
\widetilde{\bu}_f\rangle\{u_0\} =\{0\}$ if and only if $d>0$.
\end{lemma}
\proof Let $\xi=(\xi_{1}, \ldots, \xi_{n})$ be a generic point of
$\mathcal {I}$ over $\ff$ that is free from $\ff\langle
\bu_f\rangle$ and $f=u_{0}+f_0$ where $f_0=\sum\limits_{\deg(m)\geq
1} u_{m}m$.
We claim that $(\xi_{1}, \ldots,  \xi_{n},-f_0(\xi))$ is a generic
point of $\mathcal {I}_{0}$ over $\mathcal {F}\langle
\widetilde{\bu}_f\rangle$. Thus, it follows that $\mathcal {I}_{0}$
is a prime $\delta$-ideal by Lemma \ref{lm-gp}.

Clearly, $(\xi_{1}, \ldots,  \xi_{n},-f_0(\xi))$ is a zero of
$\CI_0$. Let $g$ be any $\delta$-polynomial in $\mathcal {F}\langle
\widetilde{\bu}_f\rangle\{y_{1}, \ldots, y_{n},u_0\}$ which vanishes
at $(\xi_{1}, \ldots, \xi_{n},-f_0(\xi))$. Regarding $u_0$ as the
leader of $f$, suppose the  pseudo remainder of $g$ w.r.t. $f$ is
$g_1$, then we have
 $$g\equiv g_1,\,\mod\,[f]$$
where $g_1\in \mathcal {F}\langle \widetilde{\bu}_f\rangle\{y_{1},
\ldots, y_{n}\}.$ From the above expression, $g_1(\xi)=0$. Since
$\xi$ is also a generic point of
$\ff\langle\widetilde{\bu}_f\rangle\cdot \CI$ over
$\ff\langle\widetilde{\bu}_f\rangle$, $g_1\in
\ff\langle\widetilde{\bu}_f\rangle\cdot \CI$. Thus, $g\in \CI_0$ and
it follows that $(\xi_{1}, \ldots,  \xi_{n},-f_0(\xi))$ is a generic
point of $\mathcal {I}_{0}$ over $\mathcal {F}\langle
\widetilde{\bu}_f\rangle$ and $\mathcal {I}_{0}$ is  prime.

By Lemma \ref{le-fieldextension},
 \begin{eqnarray}
\dim \mathcal {I}_{0} &=& \dtrdeg \,\mathcal {F}\langle
\widetilde{\bu}_f\rangle \langle \xi_{1}, \ldots, \xi_{n}, -f_0(\xi)
\rangle /\mathcal {F}\langle \widetilde{\bu}_f\rangle \nonumber \\
&=& \dtrdeg\, \mathcal {F}\langle \widetilde{\bu}_f\rangle \langle
\xi_{1},    \ldots, \xi_{n}\rangle /\mathcal {F}\langle
\widetilde{\bu}_f\rangle \nonumber \\ &=& \dtrdeg\,\mathcal {F}
\langle \xi_{1},    \ldots,    \xi_{n}\rangle /\mathcal {F}  = d.
\nonumber
\end{eqnarray}

Now consider the second part of the lemma. If $d=0$, then
$\dim\,\mathcal {I}_{0}=0$, so $\mathcal{I}_0\cap\mathcal {F}\langle
\widetilde{\bu}_f\rangle\{u_0\} \neq\{0\}$. Thus, if
$\mathcal{I}_0\cap\mathcal {F}\langle
\widetilde{\bu}_f\rangle\{u_0\} =\{0\}$, then $d>0$. It remains to
show that if $d>0$, then $\mathcal{I}_0\cap\mathcal {F}\langle
\widetilde{\bu}_f\rangle\{u_0\} =\{0\}$. Suppose the contrary, then
there exists a nonzero $\delta$-polynomial $p(\widetilde{\bu}_f,
u_{0}) \in \mathcal {I}_{0}\cap \mathcal {F}\{ \widetilde{\bu}_f,
u_{0}\}$. So $p(\widetilde{\bu}_f, -f_0(\xi))=0$. Then,
$\phi=-f_0(\xi)$ is $\delta$-algebraic over $\mathcal {F}\langle
\widetilde{\bu}_f\rangle $. Denote the coefficient of any monomial
$y_i$ in $f$ to be $u_{i0}$. So for any fixed $i$ when $u_{i0}$ is
specialized to $-1$  and all the other $u \in \widetilde{\bu}_f$
specialized to zero, by Theorem~\ref{th-specil}, we conclude that
$\overline{\phi}=\xi_{i}(i=1, \ldots, n)$ is $\delta$-algebraic over
$\mathcal {F}$, which contradicts to the fact that $\mathcal{I}$ has
a positive dimension. So $ \mathcal {I}_{0}\cap \mathcal {F}\langle
\widetilde{\bu}_f\rangle \{ u_{0}\}=\{0\} $.\qed

Now, we will prove the first key result of this paper, which  shows
that by adding a generic $\delta$-polynomial to a prime
$\delta$-ideal, the new $\delta$-ideal is still prime and its
dimension decreases by one. This is generally not valid if the
$\delta$-polynomial is not generic as shown in Examples
\ref{ex-dim1} and \ref{ex-dim2}.

\begin{theorem}\label{th-inter}
Let $\mathcal {I}$ be a prime $\delta$-ideal in $\mathcal {F}\{\Y\}$
with dimension $d$ and $f$ a generic $\delta$-polynomial. If $d>0$,
then $\mathcal {I}_1=[\mathcal {I}, f]$ is a prime $\delta$-ideal in
$\mathcal {F}\langle\bu_f\rangle\{\Y\}$ with dimension $d-1$. And if
$d=0$, then $\mathcal {I}_1$ is the unit ideal in $\mathcal
{F}\langle\bu_f\rangle\{\Y\}$.
\end{theorem}
\proof Firstly, we consider the case $d>0$. Let $(\xi_{1},\ldots,
\xi_{n})$ be a generic point of $\mathcal {I}$ over $\ff$ that is
free from the pure $\delta$-extension field $\ff\langle
\bu_f\rangle$ over $\ff$ and $f=u_{0}+f_0$ where
$f_0=\sum\limits_{\deg(m)\geq 1} u_{m}m$.
By Lemma \ref{lm-pu}, $\mathcal{I}_0=[\mathcal {I},
f]\subset\mathcal {F}\langle \widetilde{\bu}_f\rangle\{\Y,u_0\}$ is
a prime $\delta$-ideal, where $\widetilde{\bu}_f$ is defined in
\bref{eq-buf}.
By Lemma \ref{lm-pu},
 $\mathcal {I}_{0}\cap \mathcal {F}\langle \widetilde{\bu}_f\rangle
\{ u_{0}\}=\{0\}$.  So $\mathcal {I}_1=[\mathcal {I}_{0}]$ in
$\mathcal {F}\langle\bu_f\rangle\{\Y\}$ is not the unit ideal, for
if not, $1\in\CI_1$. Then there exist $\delta$-polynomials $g_{i}\in
\CI\subset\ff\{\Y\},$ $H_{ij},G_k\in\mathcal {F}\langle
\bu_f\rangle\{\Y\}$ such that
$1=\sum_{i,j}H_{ij}g_i^{(j)}+\sum_kG_kf^{(k)}$. By clearing
denominators in the above expression, we obtain a nonzero
$\delta$-polynomial in $\mathcal {I}_{0}\cap \mathcal
{F}\{\bu_{f}\}$, a contradiction.

Now we claim that $\mathcal{I}_1$ is prime and $\mathcal {I}_1\cap
\mathcal {F}\langle \widetilde{\bu}_f\rangle\{\Y,u_{0}\}=\mathcal
{I}_{0}$.
 Suppose $g,h\in\mathcal {F}\langle\bu_f\rangle\{\Y\}$ and $gh\in \CI_1$.
 By collecting denominators, there exist $D_1,D_2\in\mathcal {F}\{ \bu_f\}$ such that $D_1g,D_2h\in\mathcal {F}\{ \bu_f,\Y\}$.
 And $(D_1g)\cdot(D_2h)=D_1D_2(gh)\in\CI_1$ still holds.
 Similarly to the procedure above, we can find a nonzero $D\in \mathcal {F}\{ \bu_f\}$ such that $D\cdot(D_1g)\cdot(D_2h)\in\CI_0$.
 Since $\CI_0$ is prime and $\mathcal {I}_{0}\cap \mathcal {F}\langle \widetilde{\bu}_f\rangle\{ u_{0}\}=\{0\}$,
 $D_1g\in\CI_0$ or $D_2h\in\CI_0$. It follows that $g\in\CI_1$ or $h\in\CI_1$.
 Since  $\CI_1$ is not the unit ideal, $\CI_1$ is prime.
 And for any $g\in\mathcal {I}_1\cap \mathcal {F}\langle \widetilde{\bu}_f\rangle
 \{\Y,u_{0}\}$, there exists some $D\in \mathcal {F}\{ \bu_f\}$ such that $Dg\in \CI_0$, so $g\in\CI_0$.
 Thus, $\mathcal {I}_1\cap \mathcal {F}\langle \widetilde{\bu}_f\rangle
 \{\Y,u_{0}\}=\mathcal {I}_{0}$.

Suppose $\xi_{1},    \ldots,    \xi_{d}$ are $\delta$-independent
over $\mathcal {F}$. Then, $\{y_{1},\ldots, y_{d}\}$ is a parametric
set of $\mathcal {I}$. Thus each $y_{d+i}(i=1,\ldots,n-d)$ is
$\delta$-dependent with $y_{1},\ldots,y_{d}$ modulo $\mathcal
{I}_1$, since $\mathcal {I}\subset \mathcal {I}_1$. By Lemma
\ref{lm-pu}, $\dim \mathcal {I}_{0}=d$. Then $u_{0}, y_{1}, \ldots,
y_{d}$ are $\delta$-dependent modulo $\mathcal {I}_{0}$,  so
$\{y_{1},\ldots, y_{d}\}$ is $\delta$-dependent modulo $\mathcal
{I}_1$. Thus $\dim \mathcal {I}_1\leq d-1$. Now we claim $y_{1},
\ldots, y_{d-1}$ are $\delta$-independent modulo $\mathcal {I}_1$,
which proves $\dim \mathcal {I}_1=d-1$. Suppose the contrary:
$y_{1}, \ldots, y_{d-1}$ are $\delta$-dependent modulo $\mathcal
{I}_{1}$.
Thus there exists a nonzero $\delta$-polynomial $p(y_{1}, \ldots,
y_{d-1}) \in \mathcal {I}_{1}$. Take $p \in \mathcal
{F}\{\widetilde{\bu}_f, y_{1}, \ldots, y_{d-1}, u_{0}\}$. Then
\[ p(\widetilde{\bu}_f,   \xi_{1}, \ldots, \xi_{d-1}, -f_0(\xi) )=0. \] That
is, $\xi_{1},  \ldots, \xi_{d-1}, -f_0(\xi) $ are $\delta$-dependent
over $\mathcal {F} \langle \widetilde{\bu}_f \rangle $. Now we
specialize $u_{d0}$ to $-1$, and the other $u\in \widetilde{\bu}_f$
to zero, where $u_{d0}$ refers to the coefficient of the monomial
$y_d$ in $f$. Then $-f_0(\xi) $ is specialized to $\xi_{d}$. By
Theorem~\ref{th-specil}, $\xi_{1},  \ldots,  \xi_{d}$ are
$\delta$-dependent over $\mathcal {F}$, which is a contradiction. So
in this case $\dim \mathcal {I}_1=d-1$.

Now, it remains to show the case $d=0$. Since $d=0$, by Lemma
\ref{lm-pu}, $\mathcal {I}_{0}\cap \mathcal {F}\langle
\widetilde{\bu}_f \rangle \{ u_{0}\}\neq\{0\}$. So $\mathcal
{I}_{0}\cap \mathcal {F}\langle\bu_{f}\rangle\neq\{0\} $, and
consequently $\mathcal {I}_1=[\mathcal {I}_{0}]$ in $\mathcal
{F}\langle\bu_f\rangle\{\Y\}$ is the unit ideal. \qed

A special case of Theorem \ref{th-inter} is  particularly interesting
and its algebraic counterpart is often listed as a theorem in
algebraic geometry textbooks \cite[p.54, p.110]{hodge}.

\begin{theorem} \label{th-inter1}
Let $\mathcal {I}$ be a prime $\delta$-ideal in $\mathcal {F}\{\Y\}$
with dimension $d>0$. Let $\{u_{0}, u_{1}, \ldots,
u_{n}\}\subset\ee$ be a set of $\delta$-$\ff$-indeterminates. Then
$\mathcal {I}_1=[\mathcal {I}, u_{0}+u_{1}y_{1}+\cdots+u_{n}y_{n}]$
is a prime $\delta$-ideal in $\mathcal {F}\langle u_{0}, u_{1},
\ldots, u_{n}\rangle\{\Y\}$ with dimension $d-1$.
\end{theorem}

Theorem \ref{th-inter} is also valid for a wider class of
$\delta$-polynomials. A $\delta$-polynomial $f$  in $\Y$ is said to
be {\em quasi-generic}, if 1) the coefficients of $f$ as a
$\delta$-polynomial in $y_{1},\ldots,y_{n}$ are
$\delta$-indeterminates and 2) in addition to the degree zero term,
for each $1\le i\le n$, $f$ also contains at least one
$\delta$-monomial in $\mathcal {F}\{y_{i}\}\setminus \mathcal{F}$.
%
For instance, $f = u_0 + u_1 y_1 + u_2y_1y_2$ is not quasi-generic,
because $f$ contains no monomials in $\mathcal {F}\{y_{2}\}\setminus
\mathcal{F}$.

The proof for Theorem \ref{th-inter} can be easily adapted to prove the
following result.

\begin{cor}\label{cor-i}
Let $\mathcal {I}$ be a prime $\delta$-ideal in $\mathcal {F}\{\Y\}$
with dimension $d$ and $f$   a quasi-generic $\delta$-polynomial
with $\bu_f$ as the set of coefficients. If $d>0$, then $\mathcal
{I}_1=[\mathcal {I}, f]$ is a prime $\delta$-ideal in $\mathcal
{F}\langle\bu_f\rangle\{\Y\}$ with dimension $d-1$.  And if $d=0$,
then $\mathcal {I}_1$ is the unit ideal in $\mathcal {F}\langle
\bu_f\rangle\{\Y\}$.
\end{cor}

By {\em independent generic (resp. quasi-generic)
$\delta$-polynomials}, we mean that the coefficients of all of them
are $\delta$-independent over $\ff.$ As a direct consequence, we can
show that the dimension conjecture is valid for quasi-generic
$\delta$-polynomials.

\begin{theorem}[Generic Dimension Theorem]\label{th-gdt}
Let $f_{1}, \ldots, f_{r}$ be independent quasi-generic
$\delta$-polynomials  in $\mathcal {F}\langle\bu\rangle\{\Y\}$ with
$r\leq n$ and $\bu$ the set of coefficients of all $f_i$. Then
$[f_{1}, \ldots, f_{r}]\subset\ff\langle\bu\rangle\{\Y\}$ is a prime
$\delta$-ideal with dimension  $n-r$. And if $r>n$, $[f_{1}, \ldots,
f_{r}]$ is the unit ideal.
\end{theorem}
\proof We prove the theorem by induction. Let $\mathcal {I}=[0]$.
When $r=1$,  by Corollary \ref{cor-i},    $[f_{1}]$ is prime with
dimension $n-1$. Assuming this holds for $r-1$,   now consider the
case $r\leq n$. By the hypothesis,    $[f_{1}, \ldots, f_{r-1}]$ is
a prime $\delta$-ideal with dimension $n-r+1$. Note that the
coefficients of $f_r$ are  $\delta$-indeterminates over
$\ff\langle\bu_{f_1},\ldots,\bu_{f_{r-1}}\rangle$.
Using Corollary \ref{cor-i} again,   $[f_{1},    \ldots, f_{r}]$ is
a prime $\delta$-ideal with dimension $n-r$. When $r>n$, since
$[f_{1}, \ldots, f_{n}]$ is of dimension zero, by Corollary
\ref{cor-i}, $[f_{1}, \ldots, f_{r}]$ is the unit ideal.\qed

\subsection{Order of a system of generic differential polynomials } \label{SECT:pm}
In this section, we consider the order of the intersection of a
$\delta$-variety by a generic $\delta$-hypersurface. Before proving
the main result, we give a series of lemmas and theorems.

\begin{lemma}\label{le-n-1} Let $\mathcal {I}$ be a prime $\delta$-ideal in $\mathcal
{F}\{\Y\}$ with dimension $n-1$. Suppose
$\{f\}$ is a characteristic set of $\mathcal {I}$ w.r.t. some
ranking $\mathscr{R}$ and $f$ is irreducible. Then for any other
ranking $\mathscr{\overline{R}}$, $\{f\}$ is also a characteristic set of
$\mathcal {I}$.
\end{lemma}
\proof Denote $f$ to be $\overline{f}$ under the ranking
$\mathscr{\overline{R}}$. By Theorem \ref{th-sat},  $\mathcal{I} =
\sat(f)$ and $\mathcal{\overline{I}} = \sat(\overline{f})$ are prime
$\delta$-ideals with $f$ and $\overline{f}$ as characteristic sets
respectively. We need to show that  $\mathcal{I}
=\mathcal{\overline{I}}$.
Let $S$ be the separant of $f$. Then for $g\in\sat(f)$, we have $S^m
g = hf+h_{1}f'+\ldots+h_{s}f^{(s)}$ for $m,s \in \mathbb{N}$. Then,
$S^m g\in\sat(\overline{f})$. Since $\sat(\overline{f})$ is prime,
we need only to show that $S$ is not in $\sat(\overline{f})$.
Suppose the contrary, $S\in\sat(\overline{f})$. Since $S$ is
partially reduced w.r.t. $\overline{f}$, we have $S = h
\overline{f}$ for a $\delta$-polynomial $h$, which is impossible
since $S = \frac{\partial{f}}{\partial{u_f}}$. So $\mathcal{I}
\subseteq \mathcal{\overline{I}}$. Similarly, we can prove that
$\mathcal{I} \supseteq \mathcal{\overline{I}}$, thus $\mathcal{I} =
\mathcal{\overline{I}}$.\qed

If $\mathcal {S}$ is any set of $\delta$-polynomials in $\ee\{\Y\}$,
then its set of zeros in $\ee^n$ is called the $\delta$-variety of
$\mathcal {S}$, still denoted by $\mathbb{V}(\mathcal {S})$. The
following lemma generalizes a result in \cite[p.5]{cohn} to the case
of positive dimensions.
\begin{lemma} \label{th-order+1}
Let the $\delta$-variety of a system $\mathcal {S}$ of
$\delta$-polynomials in $\mathcal {F}\{\Y\}$ have a component $V$ of
dimension $d$ and order $h$. Let $\mathcal{\overline{S}}$ be
obtained from $\mathcal {S}$ by replacing $y_{1}^{(k)}$ by
$y_{1}^{(k+1)}(k=0, 1,\ldots)$ in all of the $\delta$-polynomials of
$\mathcal {S}$. Then the $\delta$-variety of $\overline{\mathcal
{S}}$ has a component $\overline{V}$ of dimension $d$ and order
$h_1$ such that $h\leq h_1\leq h+1$. Moreover, if there exists a
parametric set $\U$ not containing $y_{1}$ such that the relative
order of $\mathbb{I}(V)$ w.r.t. $\U$ is $h$, then the order of
$\overline{V}$ is $h+1$; otherwise,  the order of $\overline{V}$ is
$h$. In particular, if $d=0$, then $\ord(\overline{V})=h+1.$
\end{lemma}
\proof  Let $(\xi_{1},  \ldots,    \xi_{n})$ be a generic point of
$V$ and $\mathcal {I}=\I(V)\in \mathcal {F}\{\Y\}$. It is clear that
$[z'-\xi_{1}]$ is a prime $\delta$-ideal in $\mathcal {F}\langle
\xi_{1}, \ldots,    \xi_{n}\rangle \{ z\}$. Let $\eta$ be a generic
point of $[z'-\xi_{1}]$. Then $(\eta, \xi_{2}, \ldots, \xi_{n})$ is
a point of $\mathcal {\overline{S}}$. Suppose this point lies in a
component $\overline{V}$ of $\overline{\mathcal {S}}$, which has a
generic point $(\eta_{1}, \ldots,  \eta_{n})$. Then
$(\eta_{1},\eta_{2}, \ldots, \eta_{n})$ is specialized to $(\eta,
\xi_{2}, \ldots, \xi_{n})$ and correspondingly $(\eta_{1}',\eta_{2},
\ldots, \eta_{n})$ is specialized to $(\xi_{1}, \xi_{2}, \ldots,
\xi_{n})$. Since $(\xi_{1}, \ldots, \xi_{n})$ is a generic point of
$V$ and $(\eta_{1}',\eta_{2}, \ldots, \eta_{n})$ is a zero of
$\mathcal {S}$, the latter specialization is generic, that is,
$(\eta_{1}',\eta_{2}, \ldots, \eta_{n})$ is a generic point of $V$.
We claim that any parametric set $\U$ of $\mathcal {I}$ is a
parametric set of $\I(\overline{V})$, and $\ord_{\U}\mathcal
{I}\leq\ord_{\U}\I(\overline{V})\leq \ord_{\U}\mathcal {I}+1$, which
follows that $\dim(\overline{V})=d$ and by
Theorem~\ref{th-maxorder}, $h\leq\ord(\overline{V})\leq h+1$. Let
$\U$ be any parametric set of $\mathcal {I}$. We consider the
following two cases.

Case 1: $y_{1}\notin \U$. Suppose $\U$ is the set of $y_{2}, \ldots,
y_{d+1}$. By Corollary \ref{cor-order1}， we have
\[
 \ord_{y_{2},    \ldots,    y_{d+1}}\mathcal {I}=\trdeg \mathcal
{F}\langle \xi_{1},    \ldots,    \xi_{d},    \xi_{d+1},    \ldots,
\xi_{n}\rangle/\mathcal {F}\langle \xi_{2},    \ldots,
\xi_{d+1}\rangle.
\]
 Since $\xi_{2},    \ldots,    \xi_{d+1}$ are $\delta$-independent over $\mathcal {F}$,    $\eta_{2}, \ldots,  \eta_{d+1}$
must be $\delta$-independent over $\mathcal {F}$, i.e.
$\I(\overline{V})\cap \mathcal {F}\{\U\}=\{0\}$.
\begin{eqnarray}
 &\quad& \trdeg\, \mathcal {F}\langle \eta_{1},    \eta_{2},\ldots, \eta_{n}\rangle/\mathcal
{F}\langle \eta_{2},    \ldots,    \eta_{d+1}\rangle \nonumber \\
&\geq & \trdeg\, \mathcal {F}\langle \eta,    \xi_{2},\ldots,
  \xi_{n}\rangle/\mathcal {F}\langle \xi_{2},
\ldots,  \xi_{d+1}\rangle \nonumber \\
& &\quad (\hbox{for}\,(\eta_{1},\eta_{2},\ldots,
\eta_{n})\,\hbox{can be specialized to}\,
(\eta,\xi_{2},\ldots,\xi_{n}))\nonumber
\\ &=& \trdeg\, \mathcal {F}\langle \xi_{1},\ldots, \xi_{n}\rangle/\mathcal {F}\langle \xi_{2},
\ldots, \xi_{d+1}\rangle  +\trdeg \,\mathcal {F}\langle \xi_{1},
\ldots, \xi_{n}\rangle \langle \eta\rangle/\mathcal {F}\langle
\xi_{1}, \ldots, \xi_{n}\rangle \nonumber \\ &=&\ord_{y_{2},\ldots,
y_{d+1}}\mathcal {I}+1 \nonumber
\end{eqnarray} and
\begin{eqnarray}
& & \trdeg\, \mathcal {F}\langle \eta_{1},\ldots,\eta_{n}
\rangle/\mathcal
{F}\langle \eta_{2},    \ldots, \eta_{d+1}\rangle \nonumber \\
&\leq & 1+\trdeg\, \mathcal {F}\langle \eta_{1}', \eta_{2}\ldots,
\eta_{n}\rangle/
\mathcal {F}\langle \eta_{2}, \ldots, \eta_{d+1}\rangle \nonumber \\
&=& 1+\ord_{y_{2}, \ldots, y_{d+1}}\mathcal {I} \nonumber
\end{eqnarray}
So $\trdeg\, \mathcal {F}\langle \eta_{1},\ldots,\eta_{n}
\rangle/\mathcal {F}\langle \eta_{2},    \ldots,
\eta_{d+1}\rangle=1+\ord_{y_{2}, \ldots, y_{d+1}}\mathcal
{I}<\infty$. Thus $\overline{V}$ is of dimension $d$ and
$\{y_{2},\ldots,y_{d+1}\}$ is a parametric set of
$\I(\overline{V})$. Moreover,  the relative order of
$\I(\overline{V})$ w.r.t. $y_{2}, \ldots, y_{d+1}$ is $\ord_{y_{2},
\ldots, y_{d+1}}\mathcal {I}+1$.

Case 2: $y_{1}\in \U$. Suppose $\U=\{y_{1}, \ldots, y_{d}\}$. Then
by Corollary \ref{cor-order1}, $\ord_{\U}\mathcal {I}=\ord_{y_{1},
\ldots, y_{d}}\mathcal {I}=\trdeg\, \mathcal {F}\langle \xi_{1},
\ldots, \xi_{n}\rangle/\mathcal {F}\langle \xi_{1}, \ldots,$
$\xi_{d}\rangle$.  We have seen that $(\eta_{1}', \eta_{2}, \ldots,
\eta_{n})$ is a generic point of $V$. Since $\trdeg\,\mathcal
{F}\langle \xi_{1},\ldots,\xi_{n}\rangle\langle \eta\rangle
/\mathcal {F}\langle \xi_{1},\ldots,\xi_{n}\rangle=1$ and
$(\eta_{1},\eta_{2},$ $ \ldots, \eta_{n})$ can be specialized to
$(\eta, \xi_{2}, \ldots, \xi_{n})$, $\eta_{1}$ is algebraically
independent over $\mathcal {F}\langle \eta_{1}',\eta_{2},$
$\ldots,\eta_{n}\rangle$. So
\begin{eqnarray}
 &\quad& \trdeg\, \mathcal {F}\langle
\eta_{1},    \ldots,    \eta_{n}\rangle/\mathcal {F}\langle
\eta_{1},  \ldots,    \eta_{d} \rangle \nonumber \\ &= & \trdeg\,
\mathcal {F}(\eta_{1})\langle \eta_{1}',    \eta_{2}\ldots,
\eta_{n}\rangle/\mathcal {F}(\eta_{1})\langle \eta_{1}', \eta_{2},
\ldots,    \eta_{d}\rangle
 \nonumber \\ &=& \trdeg\, \mathcal
{F}\langle \eta_{1}',    \eta_{2},    \ldots,
\eta_{n}\rangle/\mathcal {F}\langle \eta_{1}', \eta_{2}, \ldots,
\eta_{d}\rangle \,  \nonumber \\
& & \quad(\hbox{for} \, \trdeg\,\mathcal {F}\langle
\eta_{1}', \ldots, \eta_{n}\rangle(\eta_{1})/\mathcal
{F}\langle \eta_{1}', \ldots, \eta_{n}\rangle=1) \nonumber \\
&=& \ord_{y_{1}, \ldots, y_{d}}\mathcal {I} \nonumber
\end{eqnarray}
Since $(\eta_{1},\ldots,\eta_{d})$ can be specialized to
$(\eta,\xi_{2},\ldots,\xi_{d})$ over $\mathcal {F}$,
$d\geq\dtrdeg\,\mathcal {F}\langle \eta_{1},$
$\ldots,\eta_{d}\rangle/\mathcal {F} \geq\dtrdeg\,\mathcal
{F}\langle \eta,\xi_{2},\ldots,\xi_{d}\rangle/\mathcal {F} \geq
\dtrdeg\,\mathcal {F}\langle
\xi_{1},\xi_{2},\ldots,\xi_{d}\rangle/\mathcal {F}=d$. Since
$\trdeg\, \mathcal {F}\langle \eta_{1},    \ldots,
\eta_{n}\rangle/\mathcal {F}\langle \eta_{1},  \ldots,    \eta_{d}
\rangle <\infty$, we have $\dtrdeg\,\mathcal {F}\langle
\eta_{1},\ldots,\eta_{d},\ldots,\eta_{n}\rangle/\mathcal {F}=d$.
 Thus in this case, $\dim(\overline{V})=d$ and $\U=\{y_1,\ldots,y_d\}$ is a parametric set of $\I(\overline{V})$ with
$\ord_{y_{1}, \ldots, y_{d}}\I(\overline{V})=\ord_{y_{1}, \ldots,
y_{d}}\mathcal {I}$.

Consider the two cases together, we can see $\dim(\overline{V})=d$.
And by Theorem~\ref{th-maxorder}, $h\leq\ord(\overline{V})\leq h+1$.
Moreover, if there exists a parametric set $\U$ not containing
$y_{1}$ such that the relative order of $\mathbb{I}(V)$ w.r.t. $\U$
is $h$, then the order of $\overline{V}$ is $h+1$; otherwise,  the
order of $\overline{V}$ is $h$. In particular, if $d=0$, then
$y_1\notin \U=\emptyset$. From case 1,
$\ord(\overline{V})=\ord(V)+1=h+1.$
\qed

Let $\mathcal {G}\subset\ee$ be a $\delta$-extension field of
$\mathcal {F}$. By a {\em $\delta$-$\ff$-isomorphism} of $\mathcal
{G}$, we mean a $\delta$-isomorphic mapping of $\mathcal {G}$ onto a
$\delta$-field $\mathcal {G}'\subset\ee$ such that (a) $\mathcal
{G}'$ is an extension of $\mathcal {F}$, (b) the $\delta$-isomorphic
mapping leaves each element of $\mathcal {F}$ invariant.
By means of well-ordering methods,    it is easy to show that a
$\delta$-$\mathcal {F}$-isomorphism of $\mathcal {G}$  can be
extended to a $\delta$-$\ff$-automorphism of $\ee$. We will use the
following result about $\delta$-isomorphism.

\begin{theorem} \cite{kol2} \label{th-primitive} Let $\mathcal
{G}\subset\ee$ be a $\delta$-extension field of $\mathcal {F}$ and
$\gamma \in \mathcal {G}$. A necessary and sufficient condition that
$\gamma$ be a primitive element of $\mathcal {G}$, i.e. $\mathcal
{G}=\mathcal {F}\langle \gamma\rangle$,    is that no $\mathcal
{F}$-isomorphism of $\mathcal {G}$ other than the identity leaves
$\gamma$ invariant.
\end{theorem}

The following theorem as well as Theorem~\ref{th-inter} prove
Theorem \ref{th-main1}.
\begin{theorem}\label{th-ord}
Let $\mathcal {I}$ be a prime $\delta$-$\ff$-ideal with dimension
$d>0$ and order $h$, and $f$ a generic $\delta$-polynomial of order
$s$. Then $\mathcal {I}_1=[\mathcal {I},f]$ is a prime
$\delta$-ideal in $\mathcal {F}\langle\bu_f\rangle\{\Y\}$ with
dimension $d-1$ and order $h+s$.
\end{theorem}
\proof By Theorem~\ref{th-inter},    $\mathcal {I}_1$ is prime with
dimension $d-1$. Now we prove the order of $\mathcal {I}_1$ is $h+s$.

Let $\mathscr{A}$ be a characteristic set of $\mathcal {I}$ w.r.t.
an orderly ranking $\mathscr{R}$ with $y_{1}, \ldots,    y_{d}$ as a
parametric set. By Theorem \ref{def-order-sadik},
$\ord(\mathscr{A})=h$.
Suppose $\xi=(\xi_{1}, \ldots, \xi_{n})$ is a generic point  of
$\mathcal {I}$ that is free from the pure extension field
$\ff\langle\bu_f\rangle$ over $\ff$. Let $f=u_{0}+f_0$ where
$f_0=\sum\limits_{\deg(m)\geq 1} u_{m}m$.
Let $\mathcal{I}_{0}=[\mathcal {I},f]$ in $\mathcal {F}\langle
\widetilde{\bu}_f\rangle\{y_{1}, \ldots, y_{n}, u_{0}\}$, where
$\widetilde{\bu}_f$ is defined in \bref{eq-buf}.
By Lemma \ref{lm-pu}, $\mathcal {I}_{0}$ is a prime $\delta$-ideal
of dimension $d$ with a generic zero $(\xi_{1}, \ldots, \xi_{n},
-f_0(\xi))$, and $u_{0}$ is $\delta$-independent modulo $\mathcal
{I}_{0}$. $\mathcal {I}_{0}$ and $\mathcal {I}_1$ have such
relations: Any characteristic set of $\mathcal {I}_{0}$ with $u_{0}$
in the parametric set is a characteristic set of $\mathcal {I}_1$,
and conversely, any characteristic set of $\mathcal {I}_1$, by
clearing denominators in $\ff\langle
\widetilde{\bu}_f\rangle\{u_0\}$ when necessary, is a characteristic
set of $\mathcal {I}_{0}$ with $u_{0}$ in the parametric set.
By Theorem \ref{th-maxorder},  we have
$\ord(\mathcal {I}_1)\leq \ord(\mathcal {I}_{0})$.

We claim that $\ord(\mathcal {I}_{0})\le h+s$.
As a consequence,  $\ord(\mathcal {I}_1)\leq h+s$.
To prove this claim, let $\mathcal {I}_{0}^{(i)}=[\mathcal {I},
u_{0}^{(i)}+f_0]\,(i=0,\ldots,s)$ in  $\mathcal {F}\langle
\widetilde{\bu}_f \rangle\{y_{1}, \ldots,    y_{n},   u_{0}\}. $
Similarly to the proof of Lemma~\ref{lm-pu}, $\mathcal
{I}_{0}^{(i)}$ is a prime $\delta$-ideal of dimension $d$. Let
$\bar{f}$ be the pseudo remainder of $u_{0}^{(s)}+f_0$ w.r.t.
$\mathcal {A}$ under the ranking $\mathscr{R}$. Clearly,
$\ord(\bar{f},u_{0})=s.$ It is obvious that for some orderly
ranking, $\{\mathcal{A}, \bar{f} \}$ is a characteristic set of
$\mathcal {I}_{0}^{(s)}$ with $y_1,\ldots,y_d$ as a parametric set.
So $\ord(\mathcal {I}_{0}^{(s)})=h+s$.
%
%
Using Lemma~\ref{th-order+1} $s$ times, we have $\ord(\mathcal {I}_{0})\le\ord(\mathcal{I}_{0}^{(1)})\le\cdots\le \ord(\mathcal
{I}_{0}^{(s)})=h+s$.

Now, it suffices to show $\ord(\mathcal {I}_1)\geq h+s$.
Let $w=u_{0}+\sum_{i=1}^d \sum_{j=0}^s u_{ij} y_{i}^{(j)}$ be a new
$\delta$-indeterminate.  Let $\bu_g$ be the set of coefficients of
$g=w+\sum_{i=d+1}^n \sum_{j=0}^s u_{ij} y_{i}^{(j)}+f_1$ regarded as
a $\delta$-polynomial in $w$ and $\Y$, where $f_1$ is the nonlinear
part of $f$ in $\Y$. We denote $\mathcal{F}_1=\mathcal {F}\langle
\bu_g\rangle$. Then $\mathcal {I}_{2}=[\mathcal {I},g] $ in
$\mathcal {F}_1\{y_{1}, \ldots, y_{n},w\}$ is a prime $\delta$-ideal
with a generic point $(\xi_{1}, \ldots, \xi_{n}, \gamma)$ where
$\gamma=-\sum_{i=d+1}^n \sum_{j=0}^s u_{ij}
\xi_{i}^{(j)}-f_1(\xi_{1},\ldots,\xi_{n})$. We claim that $\gamma$
is a primitive element of $\mathcal {F}_1\langle \xi_{1}, \ldots,
\xi_{d}\rangle \langle \xi_{d+1}, \ldots, \xi_{n}\rangle$ over
$\mathcal {F}_1\langle \xi_{1}, \ldots, \xi_{d}\rangle$. By
Theorem~\ref{th-primitive}, it suffices to show that no $\mathcal
{F}_1\langle \xi_{1}, \ldots, \xi_{d}\rangle$-isomorphism of
$\mathcal {F}_1\langle \xi_{1}, \ldots, \xi_{d}\rangle \langle
\xi_{d+1}, \ldots, \xi_{n}\rangle$ other than the identity leaves
$\gamma$ invariant. Let $\varphi$ be any $\delta$-$\mathcal
{F}_1\langle \xi_{1}, \ldots, \xi_{d}\rangle$-isomorphism of
$\mathcal {F}_1\langle \xi_{1}, \ldots, \xi_{d}\rangle \langle
\xi_{d+1}, \ldots, $ $\xi_{n}\rangle$ which leaves $\gamma$
invariant, and $\varphi(\xi_{d+i})=\eta_{d+i}\,(i=1,\ldots,n-d)$.
Since each $\xi_{d+i}\,(i=1,\ldots,n-d)$ is $\delta$-algebraic over
$\mathcal {F}\langle \xi_{1}, \ldots, \xi_{d} \rangle$ and $\varphi$
is an isomorphism leaving each element of $\mathcal {F}\langle
\xi_{1}, \ldots, \xi_{d} \rangle$ invariant, we can see that each
$\eta_{d+i}\,(i=1,\ldots,n-d)$ is also $\delta$-algebraic over
$\mathcal {F}\langle \xi_{1}, \ldots, \xi_{d} \rangle$. So,
$\eta_{d+i}\,(i=1,\ldots,n-d)$ are also free from
$\mathcal{F}\langle\bu_f\rangle$.
From $\varphi(\gamma)=\gamma$, we have $-\sum_{i=d+1}^n \sum_{j=0}^s
u_{ij}
\eta_{i}^{(j)}-f_1(\xi_{1},\ldots,\xi_d,\eta_{d+1},\ldots,\eta_{n})=-\sum_{i=d+1}^n
\sum_{j=0}^s u_{ij}\xi_{i}^{(j)}
-f_1(\xi_{1},\ldots,\xi_d,\xi_{d+1},\ldots,\xi_{n})$ which can be
rewritten as:
 $$
\sum_{i=d+1}^n \sum_{j=0}^s u_{ij} (\xi_{i}^{(j)}-\eta_{i}^{(j)})+
f_1(\xi_1,\ldots,\xi_n)-f_1(\xi_1,\ldots,\xi_d,\eta_{d+1},\ldots,\eta_n)=0.$$
Since $\bu_f$ are $\delta$-indeterminates over
$\ff\langle\xi_1,\ldots,\xi_n,\eta_{d+1},\ldots,\eta_n\rangle$, we
have $\xi_i-\eta_i=0(i=d+1,\ldots,n)$. So $\varphi$ must be the
identity map and the claim follows.

Since $\mathcal {F}_1\langle \xi_{1}, \ldots, \xi_{d}, \xi_{d+1},
\ldots, \xi_{n}\rangle=\mathcal {F}_1\langle \xi_{1}, \ldots,
\xi_{d}\rangle \langle \gamma\rangle$, $\gamma$ is
$\delta$-algebraic over $\mathcal {F}_1\langle \xi_{1},$ $\ldots,$
$\xi_{d}\rangle$ and each $\xi_{d+i} \in \mathcal {F}_1\langle
\xi_{1}, \ldots, \xi_{d}\rangle \langle \gamma\rangle
(i=1,\ldots,n-d)$.
Let $R(\xi_{1},\ldots, \xi_{d}, w)$ be an irreducible
$\delta$-polynomial in $\mathcal {F}_1\langle \xi_{1},$ $\ldots,$
$\xi_{d}\rangle\{w\}$ annulling $\gamma$ of the lowest order. By
clearing denominators when necessary, suppose $R(y_{1},\ldots,
y_{d}, w)$ is an irreducible $\delta$-polynomial in $\mathcal
{F}_1\{y_{1},\ldots,y_{d},w\}$. Clearly, $R(y_{1},\ldots, y_{d},
w)\in \CI_2$. And there exist $A_{i}\in \mathcal {I}_{2}$ with the
form $A_i = P_i(y_{1},\ldots, y_{d},w)y_{d+i}+ Q_i(y_{1},\ldots,
y_{d},w)$ $(i=1, \ldots, n-d)$, which are reduced w.r.t. $R$.  Since
$\mathcal {I}_{2}\cap \mathcal {F}_1\{y_{1},\ldots,y_{d},w\}$ is a
$d$-dimensional prime $\delta$-ideal, by Lemma~\ref{le-n-1}, $\{R\}$
is its characteristic set w.r.t. any ranking. So for the elimination
ranking $y_{1}\prec \ldots \prec y_{d}\prec w\prec y_{d+1}\prec
\ldots\prec y_{n}$,    a characteristic set of $ \mathcal {I}_{2}$
is $\{R(y_{1}, \ldots, y_{d}, w), A_{1}, \ldots, A_{n-d}\}$. Since
$\mathcal {F}_1\langle \xi_{1}, \ldots,
\xi_{n},\gamma\rangle=\mathcal {F}_1\langle \xi_{1},  \ldots,
\xi_{n}\rangle$, by Corollary~\ref{cor-order1}, $\ord_{y_{1},
\ldots, y_{d}}\mathcal{I}_{2}= \ord_{y_{1}, \ldots, y_{d}}(\mathcal
{I})=\ord(\mathcal{A})=h$. Thus, $\ord(R, w)=h. $

Let $\U = \{u_{ij}:\, i=1, \ldots, d; j=0,\ldots,s\}$. In $\mathcal
{F}_1\langle \U \rangle\{w, y_{1}, \ldots,    y_{n}\}$,
$\mathcal{I}_{2}$ is also prime with $R(y_{1}, \ldots, y_{d},w),
A_{1}, \ldots, A_{n-d}$ as a characteristic set w.r.t. the
elimination ranking $y_{1} \prec \ldots \prec y_{d}\prec w\prec
y_{d+1}\prec \ldots\prec y_{n}$. Let

\[ \begin{array}{ccc} \phi:\mathcal {F}_1\langle \U
\rangle\{y_{1}, \ldots, y_{n},w\} &\longrightarrow &\mathcal
{F}_1\langle \U \rangle\{y_{1},\ldots, y_{n}, u_{0}\}
\\ w & \quad & u_{0}+\sum_{i=1}^d \sum_{j=0}^s u_{ij}
y_{i}^{(j)}\\ y_{i}& \quad & y_{i} \end{array} \]

\noindent be a $\delta$-homomorphism over $\mathcal {F}_1\langle \U
\rangle$. Clearly,    this is a $\delta$-isomorphism which maps $
\mathcal {I}_{2}$ to $\mathcal {I}_{0}$. It is obvious that
$\mathcal {I}_{0}$ has $\phi(R),    \phi(A_{1}), \ldots,
\phi(A_{n-d})$ as a characteristic set w.r.t. the elimination
ranking $y_{1}\prec \ldots \prec y_{d} \prec u_{0}\prec y_{d+1}\prec
\ldots\prec y_{n}$ with $\rank(\phi(A_{i}))=y_{d+i}(i=1, \ldots,
n-d)$. We claim that $\ord(\phi(R),    y_{1})\geq h+s$. If $\ord(R
,y_{1})\geq h+s$, rewrite $R$ in the form $R=\sum_{\psi_{\nu}(w)\neq
1} p_{\nu}(y_{1}, \ldots, y_{d})\psi_{\nu}(w)+p(y_{1}, \ldots,
y_{d})$ where $\psi_{\nu}(w)$ are monomials in $w$ and its
derivatives.
Then
\begin{eqnarray}
\phi(R)&=& \sum_{\psi_{\nu}\neq 1} p_{\nu}(y_{1},    \ldots,
y_{d})\psi_{\nu}(u_{0}+\sum_{i=1}^d
\sum_{j=0}^s u_{ij}y_{i}^{(j)})+p(y_{1},    \ldots,    y_{d}) \nonumber \\
&=&\sum_{\psi_{\nu}\neq
1}p_{\nu}(y_{1},    \ldots,    y_{d})\psi_{\nu}(u_{0})+p(y_{1},    \ldots,    y_{d}) \nonumber \\
& &+\hbox{terms involving}\, u_{ij}(i=1,\ldots,
d;j=0,\ldots,s) \, \hbox{ and their derivatives.} \nonumber
\end{eqnarray}
\noindent Clearly, in this case we have $\ord(\phi(R), y_{1})\geq
\max\{\ord(p_{\nu}, y_{1}),\ord(p, y_{1})\}=\ord(R,$ $ y_{1})\geq h+s$.
If $\ord(R ,y_{1})<h+s$, rewrite $R$ as a polynomial in $w^{(h)}$,
that is, $R=I_{l}(w^{(h)})^l+I_{l-1}(w^{(h)})^{l-1}+\cdots+I_{0}$.
Then $\phi(R)=\phi(I_{l})[(u_{0}+\sum_{i=1}^d \sum_{j=0}^s
u_{ij}y_{i}^{(j)})^{(h)}]^l+\phi(I_{l-1})[(u_{0}+\sum_{i=1}^d
\sum_{j=0}^s u_{ij}y_{i}^{(j)})^{(h)}]^{l-1}+\cdots+\phi(I_{0})$. Since
$\ord(\phi(I_k),y_1)<h+s$ ($k=0,\ldots,l$), we have exactly
$\ord(\phi(R), y_{1})=h+s$.  Thus, consider the two cases together,
$\ord(\phi(R), y_{1})\geq h+s$.

Since $\mathcal {I}_{0} \cap \mathcal {F}_1\langle\U\rangle\{y_{1},
\ldots, y_{d},u_0\}$ is a $d$-dimensional prime $\delta$-ideal, by
Lemma \ref{le-n-1}, $\{\phi(R)\}$ is its characteristic set w.r.t.
any ranking,  in particular, for the elimination ranking $u_{0}\prec
y_{2}\prec \ldots \prec y_{d} \prec y_{1}$. So w.r.t. the
elimination ranking $u_{0}\prec y_{2}\prec \ldots \prec y_{d} \prec
y_{1}\prec y_{d+1} \prec \ldots \prec y_{n}$,    $\{\phi(R),
\phi(A_{1}), \ldots, \phi(A_{n-d})\}$ is a characteristic set of
$\mathcal{I}_{0}$,  thus a characteristic set of $\mathcal {I}_1$.
By Theorem \ref{th-maxorder}, $\ord(\mathcal{I}_1)\geq\ord_{y_{2},
\ldots, y_{d}}\mathcal{I}_1\geq h+s$.

Thus, the order of $\mathcal{I}_1$ is $h+s$. \qed

As a consequence, Theorem \ref{th-inter1} can be strengthened as follows.

\begin{theorem} \label{th-inter2}
Let $\mathcal {I}$ be a prime $\delta$-ideal in $\mathcal {F}\{\Y\}$
with dimension $d>0$ and order $h$. Let $\{u_{0}, u_{1}, \ldots,
u_{n}\}\subset \ee $ be a set of $\delta$-$\ff$-indeterminates. Then
$\mathcal {I}_1=[\mathcal {I}, u_{0}+u_{1}y_{1}+\cdots+u_{n}y_{n}]$
is a prime $\delta$-ideal in $\mathcal {F}\langle u_{0}, u_{1},
\ldots, u_{n}\rangle\{\Y\}$ with dimension $d-1$ and order $h$.
\end{theorem}

As another consequence, the dimension theorem for generic
$\delta$-polynomials can be strengthened as follows.

\begin{theorem}
Let $f_{1},  \ldots, f_{r} (r\le n)$ be independent generic
$\delta$-polynomials with  each $ f_{i}$ of order $s_{i}$. Then
$[f_{1}, \ldots, f_{r}]$ is a prime $\delta$-ideal with dimension
$n-r$ and order $\sum_{i=1}^r s_{i}$ over
$\ff\langle\bu_{f_1},\ldots,\bu_{f_r}\rangle$.
\end{theorem}
\proof We will prove the theorem  by induction on $r$. Let
$\CI=[0]\subset\ff\{\Y\}$. Clearly, $\CI$ is a prime $\delta$-ideal
of dimension $n$ and order $0$.
For $r=1$, by Theorem~\ref{th-ord}, $[f_1]=[\CI,f_1]$ is a prime
$\delta$-ideal of dimension $n-1$ and order $s_1$. So the assertion
holds for $r=1$.
Now suppose the assertion holds for $r-1$, we now prove it for $r$.
By the hypothesis, $\CI_{r-1}=[f_1,\ldots,f_{r-1}]$ is a prime
$\delta$-ideal of dimension $n-r+1$ and order $\sum_{i=1}^{r-1}
s_{i}$ over $\ff\langle\bu_{f_1},\ldots,\bu_{f_{r-1}}\rangle$.
Since $f_1,\ldots,f_r$ are independent generic $\delta$-polynomials,
using Theorem~\ref{th-ord} again, $\CI_r=[f_1,\ldots,f_r]$ is a
prime $\delta$-ideal of dimension $n-r$ and order $\sum_{i=1}^r
s_{i}$ over $\ff\langle\bu_{f_1},\ldots,\bu_{f_r}\rangle$.  Thus,
the theorem is proved. \qed

\begin{remark}
When $f$ is a quasi-generic $\delta$-polynomial, Theorem
\ref{th-ord} may not be true. A counter example is as follows. Let
$\mathcal {I}=[y_{2},\ldots,y_{n}] \in \mathcal{F}\{\Y\}$ and
$f=u_{0}+u_{1}y_{1}+u_{2}y_{2}''+\cdots+u_{n}y_{n}''$. Clearly, $f$
is a quasi-generic $\delta$-polynomial and $[\mathcal {I}, f]$ is a
prime $\delta$-ideal of dimension 0. But $\ord([\mathcal {I},
f])=\ord(\mathcal{I})=0 \neq \ord(\mathcal{I})+\ord(f)=2$.
\end{remark}

\section{Chow form for an irreducible differential variety}
\label{sec-chow}

In this section, we define the differential Chow form  and establish its properties by proving
Theorem \ref{th-main2}.

\subsection{Definition of differential Chow form}
Throughout this section, we assume that  $V\subset\ee^n$ is an
irreducible $\delta$-$\ff$-variety with dimension $d$ and
$\CI=\I(V)\subset\ff\{\Y\}$. Let
 \begin{eqnarray}\label{eq-bu}
 \widehat{\bu}&=&\{u_{ij} \,(i=0,  \ldots, d;j=0, \ldots,   n)\}
 \end{eqnarray}be
 $(d+1)(n+1)$ $\delta$-$\ff$-indeterminates in $\ee$.

By saying that a point $\eta\in\ee^n$ is free from the pure
$\delta$-extension field $\ff\langle\widehat{\bu}\rangle$ over
$\ff$, we mean $\widehat{\bu}$ are $\delta$-indeterminates over
$\ff\langle\eta\rangle$.
Let
 $\xi=(\xi_{1}, \ldots, \xi_{n})$
be a generic point of $V$ which is free from
$\ff\langle\widehat{\bu}\rangle$ and $\zeta_{0}, \zeta_{1}, \ldots,$
$ \zeta_{d}$ elements of $\mathcal {F}\langle\bu, \xi_{1},$ $
\ldots, \xi_{n} \rangle$:
\begin{eqnarray} \label{equ1}
\zeta_{\sigma}=- \sum_{\rho=1}^n u_{ \sigma \rho}\xi_{\rho}\,
(\sigma=0,    \ldots,    d)
\end{eqnarray}
where
 \begin{eqnarray*}
 {\bu}&=&\{u_{ij} \,(i=0,  \ldots, d;j=1, \ldots,    n)\}.
 \end{eqnarray*}
 We will show that the $\delta$-transcendence degree of
$\zeta_{0}, \ldots,    \zeta_{d}$ over $\mathcal
{F}\langle\bu\rangle$ is $d$.
\begin{lemma} \label{lem-d}
 $\dtrdeg\,\mathcal {F}\langle\bu\rangle \langle
\zeta_{0},    \ldots,    \zeta_{d} \rangle/\mathcal {F}\langle
\bu\rangle=d$. Furthermore, if $d>0$,  $\zeta_{1},    \ldots,
\zeta_{d}$ are $\delta$-independent over $\mathcal {F}\langle
\bu\rangle$.
\end{lemma}
\proof  By Lemma \ref{le-fieldextension}, $\dtrdeg\, \mathcal
{F}\langle\bu\rangle \langle \xi_{1},    \ldots,
\xi_{n}\rangle/\mathcal {F}\langle\bu\rangle=d$. Since the $d+1$
elements $\zeta_{0}, \ldots, \zeta_{d}$ belong to $\mathcal
{F}\langle\bu, \xi_{1},    \ldots, \xi_{n}\rangle$, they are
$\delta$-dependent over $\mathcal {F}\langle\bu\rangle$. Then, we
have $\dtrdeg \,\mathcal {F}\langle\bu\rangle \langle \zeta_{0},
\ldots,    \zeta_{d} \rangle/\mathcal {F}\langle\bu\rangle\leq d$.
Thus, if $d=0$, we have  $\dtrdeg \,\mathcal {F}\langle\bu\rangle \langle
\zeta_{0}, \ldots,    \zeta_{d} \rangle/$ $\mathcal
{F}\langle\bu\rangle=0$.

Now, suppose $d>0$.
We claim that $\zeta_{1},    \ldots,    \zeta_{d}$ are
$\delta$-independent over $\mathcal {F}\langle\bu\rangle$,
 thus it follows that $\dtrdeg\,\mathcal {F}\langle\bu\rangle \langle \zeta_{0},    \ldots,    \zeta_{d}
\rangle/\mathcal {F}\langle\bu\rangle=d$.  Suppose the contrary.
Since $\zeta_i \in \mathcal {F}\{\bu, \xi_{1},    \ldots,
\xi_{n}\}$,   when we  specialize $u_{ij}$ to $ -\delta_{k_{i}j}$
($j=1,\ldots,n$, $k_i \in \{1,\ldots,n\}$), $\zeta_{i}$ will be
specialized to $ \xi_{k_{i}}$. Then from Theorem~\ref{th-specil}, we
conclude that $\xi_{k_1}, \ldots, \xi_{k_d}$ are $\delta$-dependent
over $\mathcal {F}$. Since we can choose $k_1,\ldots,k_d$ so that
$\xi_{k_1}, \ldots, \xi_{k_d}$ are $\delta$-independent over
$\mathcal {F}$, it amounts to a contradiction. Thus the claim is
proved.  \qed

Let $\I_{\zeta}$ be the prime $\delta$-ideal in the
$\delta$-polynomial ring $\mathcal{R} = \mathcal {F}\langle\bu
\rangle\{z_{0},$ $ \ldots, z_{d}\}$ having $\zeta=(\zeta_{0},
\ldots, \zeta_{d})$ as a generic point. By Lemma~\ref{lem-d}, the
dimension of $\I_{\zeta}$ is $d$. By Theorem \ref{relativehubert},
the characteristic set of
 $\I_{\zeta}$ w.r.t. any ranking  consists of an irreducible $\delta$-polynomial
 $f(z_{0},    \ldots,  z_{d})$ in $\mathcal {R}$ and
\begin{equation}\label{eq-def} \I_{\zeta} = \sat(f).\end{equation}
Since the coefficients of $f(z_{0},    \ldots,  z_{d})$ are elements
in $\mathcal {F}\langle\bu \rangle$, without loss of generality, we
assume that $f(\bu; z_{0}, \ldots,z_{d})$ is an irreducible
$\delta$-polynomial in $\mathcal {F}\{\bu; z_{0}, \ldots, z_{d}\}$.
We shall subsequently replace $\{z_{0},    \ldots,    z_{d}\}$ by
$\{u_{00}, \ldots,    u_{d0}\} \subset \widehat{\bu}$,   and obtain
\begin{equation}\label{eq-cf}
F(\bu_{0},\bu_{1},    \ldots, \bu_{d})=f(\bu;u_{00},    \ldots,
u_{d0})\in\ff\{\widehat{\bu}\},
\end{equation}
where $\widehat{\bu}$ is from \bref{eq-bu} and $\bu_{i}=(u_{i0},
\ldots, u_{in})$ for $i=0,\ldots,d$.

\begin{definition}[Differential Chow form]
The $\delta$-polynomial defined in \bref{eq-cf} is called the {\em
differential Chow form} of $V$ or the prime $\delta$-ideal
$\CI=\I(V)$.
\end{definition}

A {\em generic $\delta$-hyperplane} is the zero set of
$u_{0}+u_{1}y_{1}+\cdots+u_{n}y_{n}=0$ where $u_{i}$ are
$\delta$-$\ff$-indeterminates. The following result shows that the
differential Chow form can be obtained by intersecting $\CI$ with
$d+1$ generic $\delta$-hyperplanes.

\begin{lemma} \label{lm-cf1}
Using the notations introduced above, let $\mathcal{I}=
\mathbb{I}(V)$ and
\begin{equation}\label{eq-genhp}
 \P_{i}=u_{i0}+u_{i1}y_{1}+\cdots+u_{in}y_{n}\,(i=0,\ldots,d)
 \end{equation}
where $u_{ij}$ are from \bref{eq-bu}.
Then $[\mathcal {I},\P_{0},\P_1,\ldots,\P_{d}]$ is a prime
$\delta$-ideal in $\mathcal {F}\langle\bu\rangle\{u_{00},u_{10},$
$\ldots,u_{d0},y_{1},\ldots,y_{n}\}$ and  $[\mathcal
{I},\P_{0},\P_1,\ldots,\P_{d}]\cap \mathcal
{F}\langle\bu\rangle\{u_{00},u_{10},\ldots,u_{d0}\}=\sat(f)$.
\end{lemma}
\proof Similarly to the proof of Lemma~\ref{lm-pu}, it is easy to
show that $[\mathcal {I},\P_{0},\P_1,\ldots,\P_{d}]$ is a prime
$\delta$-ideal with a generic zero $(\zeta,\xi)$. Denote $[\mathcal
{I},\P_{0},\P_1,\ldots,\P_{d}]$ by $\I_{\zeta,\xi}$. Then,
$\I_{\zeta,\xi}\cap \mathcal
{F}\langle\bu\rangle\{u_{00},u_{10},\ldots,u_{d0}\}$ is a prime
$\delta$-ideal with a generic zero $\zeta$, which implies
$\I_{\zeta,\xi}\cap \mathcal
{F}\langle\bu\rangle\{u_{00},u_{10},\ldots,u_{d0}\}=\I_{\zeta}=\sat(f)$.\qed

In the following context, we will denote $[\mathcal
{I},\P_{0},\P_1,\ldots,\P_{d}]$ by $\I_{\zeta,\xi}$.
\begin{remark}\label{rm-comp}
From Lemma \ref{lm-cf1}, we have two observations.
Firstly, the Chow form for a $\delta$-variety is independent of  the
generic point used in \bref{equ1}. The $\delta$-ideals $\I_{\zeta}$
and $\I_{\zeta,\xi}$ are also independent of the choice of $\xi$.
Secondly, we can compute the Chow form of $V$ with the
$\delta$-elimination algorithms
\cite{boulier2010,ardm1,ritt,sit,alexey} if we know a set of
finitely many generating $\delta$-polynomials for $V$.
Furthermore, given a characteristic set $\mathcal {A}$ of
$\mathbb{I}(V)$, we can also compute its Chow form. Indeed, from
Lemma~\ref{lm-cf1}, it suffices to compute a characteristic set of
$\I_{\zeta,\xi}$ w.r.t. a ranking $\U<<\Y$ (elimination ranking
between elements of $\U=\{u_{00},\ldots,u_{d0}\}$ and $\Y$. It is
clear that $\I_{\zeta,\xi}$ has a characteristic set $\{\mathcal
{A},\P_{0},\ldots,\P_{d} \}$ w.r.t. a ranking $\Y<<\U$. Then using
the algorithms given by Boulier et al \cite{boulier2010} and
Golubitsky et al \cite{alexey} for transforming a characteristic
decomposition of a radical $\delta$-ideal from one ranking to
another, we can obtain the Chow form.
\end{remark}


\begin{example}\label{example1}
Consider the case $n=1$. By Theorem \ref{relativehubert}, a prime
$\delta$-ideal in $\ff\{y_1\}$ is of the form $\sat(p)$ where
$p\in\ff\{y_1\}$ is irreducible. The zero set of $\sat(p)$ is an
irreducible $\delta$-$\ff$-variety in the {\em affine line}:
$\mathbb{A}^1(\ee)$. Let $\P_0 = u_{00} + u_{01}y_1$ and $\bu_0 =
(u_{00}, u_{01})$. Then the Chow form of $\sat(p)$ is
$F(\bu_0)=u_{00}^d p(-\frac{u_{01}}{u_{00}})$ where $d$ is a natural
number such that $F(\bu_0)$ is an irreducible $\delta$-polynomial in
$\ff\{\bu_0\}$.
For instance, let $p=(y_1')^2-4y_1$. Then the Chow form of $\sat(p)$
is $F(\bu_{0})=u_{01}^2
(u_{00}')^2-2u_{01}u_{01}'u_{00}u_{00}'+(u_{01}')^2
u_{00}^2+4u_{01}^3u_{00}$.
\end{example}

\begin{example}\label{ex-cf1}
If $V$ is an irreducible $\delta$-variety of dimension $n-1$ and its
corresponding prime $\delta$-ideal is $\mathcal {I}=\sat(p)
\subset\mathcal {F}\{\Y\}$. Then its Chow form is
$F(\bu_{0},\ldots,\bu_{n-1})=D^m p\big(
\frac{D_{1}}{D},\ldots,\frac{D_{n}}{D}\big)$, where
 \[ D=\left|\begin{array}{llcl}
u_{01}&u_{02}& \ldots &u_{0n}\\
u_{11}&u_{12}&\ldots&u_{1n}\\
\multicolumn{4}{c} \dotfill \\
u_{n-1,1}&u_{n-1,2}&\ldots&u_{n-1,n}
\end{array}  \right| \]
and $D_{i}(i=1,\ldots,n)$ is the determinant of the matrix formed by
replacing the $i$-th column of $D$ by the column vector
$(-u_{00},-u_{10},\ldots,-u_{n-1,0})^T$, and $m$ is the minimal
integer such that $D^m p\big(
\frac{D_{1}}{D},\ldots,\frac{D_{n}}{D}\big) \in \mathcal
{F}\{\bu;u_{00},\ldots,u_{d0}\}$.
\end{example}

\begin{example} \label{example3}
Let $V$ be the irreducible $\delta$-variety corresponding to
$[y_{1}'+1,y_{2}'] \in \Q(t)\{y_{1},y_{2}\}$. It is of dimension
zero and the Chow form of $V$ is $F(\bu_{0})=u_{02}
u_{01}'u_{00}''+u_{01}''u_{01}u_{02}-2u_{02}(u_{01}')^2-u_{02}u_{01}''u_{00}'-u_{02}'u_{01}u_{00}''-u_{02}''u_{01}^2-u_{01}'u_{02}''u_{00}+2u_{02}'u_{01}u_{01}'
+u_{02}''u_{01}u_{00}'+u_{01}''u_{02}'u_{00}$, where
$\bu_0=(u_{00},u_{01},u_{02})$.
\end{example}

\subsection{The order of differential Chow form}
In this section, we will show that the order of the differential
Chow form is the same as that of the corresponding $\delta$-variety.
%
%
Before this, we give the following lemmas.
\begin{lemma} \label{le-divisible}
Let $\zeta_{0},    \zeta_{1},    \ldots,    \zeta_{d}$ be defined in
\bref{equ1},    and $f(\bu;    u_{00},    \ldots,    u_{d0})$
 the Chow form of $V$. Then for any $p(\bu;    u_{00},
\ldots, u_{d0}) \in \mathcal {F}\langle\bu \rangle \{u_{00}, \ldots,
u_{d0}\}$ with $\ord(f)=\ord(p)$,    such that $p(\bu; \zeta_{0},
\ldots, \zeta_{d})=0$,    we have $p(\bu;    u_{00}, \ldots,
u_{d0})=f(\bu;  u_{00},    \ldots,    u_{d0})$ $h(\bu; u_{00},
\ldots,    u_{d0})$,    where $h(\bu;u_{00},\ldots, u_{d0})$ is in
$\mathcal {F}\langle\bu \rangle\{u_{00}, \ldots,    u_{d0}\}$.
\end{lemma}
\proof  Since $\{f\}$ is the characteristic set of
$\mathbb{I}_{\zeta}$ w.r.t. an orderly ranking, and $p \in
\mathbb{I}_{\zeta}$ with $\ord(f)=\ord(p)$,   so $\init_{f}^m p=fh$
for some $m\in\mathbb{N}$. Since $f$ is irreducible, we can see that
$f$ divides $p$. \qed

The Chow form $f(\bu; u_{00}, \ldots, u_{d0})$ has certain symmetric
properties as shown by the following results.
\begin{lemma} \label{le-sym}
Let $F(\bu_{0},\bu_{1},\ldots,\bu_{d})$ be the Chow form of an
irreducible $\delta$-$\ff$-variety $V$ and $F^*(\bu_{0},$
$\bu_{1},\ldots,$ $\bu_{d})$ obtained from $F$ by interchanging
$\bu_\rho$ and $\bu_\tau$. Then $F^*$ and $F$ differ at most by a
sign. Furthermore, $\ord(F,u_{ij})$ $(i=0,\ldots,d; j=0,1,\ldots,n)$
are the same for all $u_{ij}$ occurring in $F$. In particular,
$u_{i0}\,(i=0,\ldots,d)$ appear effectively in $F$. And a necessary
and sufficient condition for some $u_{ij}\,(j>0)$ not occurring
effectively in $F$ is $y_{j} \in \mathbb{I}(V).$
\end{lemma}
\proof Consider the $\delta$-automorphism $\phi$ of $\mathcal
{F}\langle\xi_1,\ldots,\xi_n\rangle\langle\bu\rangle$ over $\mathcal
{F}\langle \xi_1,\ldots,$ $\xi_n\rangle$:
$\phi(u_{ij})=u_{ij}^{*}=\left\{ \begin{array}{ll}
u_{ij},&i\neq\rho,\tau\\u_{\tau j},&i=\rho\\u_{\rho j},&i=\tau
\end{array} \right.$. Of course, $\phi(\zeta_{i})=\zeta_{i}^{*}=\left\{ \begin{array}{ll}
\zeta_{i},&i\neq\rho,\tau\\\zeta_{\tau},&i=\rho\\\zeta_{\rho},&i=\tau
\end{array} \right.$. Then
$\phi(f(\bu;\zeta_{0},\ldots,\zeta_\rho,\ldots,\zeta_\tau,\ldots,\zeta_{d}))=f(\bu^{*};\zeta_{0},\ldots,\zeta_\tau,$
$\ldots,$ $\zeta_\rho$, $\ldots,\zeta_{d})=0$. Instead of $f(\bu;
z_{0},$ $ \ldots, z_{d})$, we obtain another $\delta$-polynomial
 $p(\bu; z_{0}, \ldots, z_{d})=f(\bu^{*};z_{0},$ $ \ldots, z_{\tau},
\ldots, z_{\rho},$ $\ldots, z_{d})\in\mathbb{I}_\zeta$. Since the
two $\delta$-polynomials $f$ and $p$ have the same
 order and degree, and as algebraic polynomials they have the same content,  by  Lemma \ref{le-divisible},
     $f(\bu^{*}; z_{0},    \ldots,$   $ z_{\tau},$ $ \ldots,z_{\rho},    \ldots,    z_{d})$
 can only differ by a sign with
 $f(\bu;    z_{0},    \ldots,    z_{\rho},    \ldots,$   $ z_{\tau},    \ldots,    z_{d})$.
 So we conclude that
 $F(\bu_{0},    \bu_{1},    \ldots,    \bu_{d})$ produces at most a change of sign if
 we interchange $\bu_{\rho}$ with $\bu_{\tau}$. In particular, each $u_{i0}$ appears effectively in $F$ and $\ord(F, u_{i0})$ are the same for all
 $i=0,1,\ldots,d$.
Suppose $\ord(F,u_{i0})=s$. For $j \neq 0$, we consider
$\ord(F,u_{ij})$. If $\ord(F,u_{ij})=l>s$, then we differentiate
$f(\bu;\zeta_{0},\ldots,\zeta_{d})=0$ on both sides w.r.t.
$u_{ij}^{(l)}$. Thus $\frac{\partial f}{\partial
u_{ij}^{(l)}}(\bu;\zeta_{0},\ldots,\zeta_{d})=0$, which amounts to a
contradiction by Lemma~\ref{le-divisible}.  If $\ord(F,u_{ij})=l<s$,
then we differentiate $f(\bu;\zeta_{0},\ldots,\zeta_{d})=0$ on both
sides w.r.t. $u_{ij}^{(s)}$. Thus $\frac{\partial f}{\partial
u_{i0}^{(s)}}(\bu;\zeta_0,\ldots,\zeta_d)\cdot(-\xi_{j})=0$. Since
$\frac{\partial f}{\partial
u_{i0}^{(s)}}(\bu;\zeta_0,\ldots,\zeta_d)\neq0$, we have
$\xi_{j}=0$. And $y_{j} \in \mathbb{I}(V) \Longleftrightarrow
\xi_{j}=0 \Longleftrightarrow \zeta_{i}$ is free of $u_{ij}$
$\Longleftrightarrow \frac{\partial f}{\partial u_{ij}^{(k)}}=0$ for
all $k \in \mathbb{Z}^{+}$ $\Longleftrightarrow$  $u_{ij}$ does not
appear in $F$.
From the above, if $u_{ij}$ occurs effectively in $F$,
$\ord(F,u_{ij})=s$, which completes the theorem.
\qed

The {\em order of the Chow form} is defined to be $\ord(f)=\ord(F)=\ord(f,
u_{i0})$ for any $i\in\{0,\ldots,d\}$. By Lemma \ref{le-sym},
$\ord(f)$ is equal to $\ord(F,u_{ij})$ for those $u_{ij}$ occurring
in $F$.

The following lemma gives a property for the $\delta$-ideal
$\I_{\zeta,\xi}$ defined in Lemma \ref{lm-cf1}.
\begin{lemma} \label{lm-cf2}
Let $F(\bu_{0}, \bu_{1}, \ldots, \bu_{d})=f(\bu;$ $u_{00},
u_{10},\ldots, u_{d0})$ be the Chow form of a prime
$\delta$-$\ff$-ideal $\mathcal {I}$ and  $s=\ord(f)$.
Then
$$\A=\{f, S_fy_1-\frac{\partial f}{\partial u_{01}^{(s)}},\, \ldots,
        S_fy_n-\frac{\partial f}{\partial u_{0n}^{(s)}}\}
$$
is a characteristic set of the prime $\delta$-ideal $\I_{\zeta,\xi}
= [\mathcal {I},\P_{0},\P_1,\ldots,\P_{d}]$ in $\mathcal
{F}\langle\bu\rangle\{u_{00},u_{10},$ $\ldots,u_{d0},$ $\Y\}$ w.r.t.
the elimination ranking $u_{d0}\prec \ldots \prec u_{00}\prec
y_{1}\prec\ldots\prec y_{n}$, where $S_f = \frac{\partial
f}{\partial u_{00}^{(s)}}$.
\end{lemma}
\proof From Lemma~\ref{lm-cf1}, $\I_{\zeta,\xi} $ is a prime
$\delta$-ideal of dimension $d$ with a generic point
$(\zeta_{0},\ldots,\zeta_{d},\xi_{1},\ldots,\xi_{n})$. From
Lemma~\ref{lem-d}, $u_{10},\ldots,u_{d0}$ is a parametric set of
$\I_{\zeta,\xi} $. If we differentiate
$f(\bu;\zeta_{0},\ldots,\zeta_{d})=0$ w.r.t. $u_{0\rho}^{(s)}$
($\rho=1,\cdots,n$),  then we have $\overline{\frac{\partial
f}{\partial u_{0\rho}^{(s)}}}-\xi_{\rho} \overline{S}_f=0$, where
$\overline{\frac{\partial f}{\partial u_{0\rho}^{(s)}}}$ and
$\overline{S}_f$ are obtained by replacing $(u_{00},\ldots,u_{d0})$
with $(\zeta_{0},\ldots,\zeta_{d})$ in $ \frac{\partial f}{\partial
u_{0\rho}^{(s)}}$ and $S_f$ respectively. So $T_{\rho}=S_f
y_{\rho}-\frac{\partial f}{\partial u_{0\rho}^{(s)}} \in
\I_{\zeta,\xi} $. Since $f$ is irreducible, we have  $S_f\not\in
\I_{\zeta,\xi} $.
Also note that $T_i$ is linear in $y_i$. $\mathcal {A}$ must be a
characteristic set\footnote{Here $\A$ is a $\delta$-chain. See
Remark \ref{re-1}.} of $\I_{\zeta,\xi} $ w.r.t. the elimination
ranking $u_{d0}\prec \ldots \prec u_{00}\prec y_{1}\prec\ldots\prec
y_{n}$.  \qed

\vskip5pt Now we give the first main property of the differential
Chow form.
\begin{theorem} \label{th-choworder}
Let $\mathcal {I}$ be a prime $\delta$-$\ff$-ideal with dimension
$d$, and $f(\bu;$ $u_{00}, u_{10}, \ldots,$ $u_{d0})$ its Chow form.
Then $\ord(f)=\ord(\mathcal {I})$.
\end{theorem}
\proof  Use the notations  $\xi_{i},\zeta_{i}$, and $\P_{i}$
introduced in \bref{equ1} and \bref{eq-genhp}.
 Let $\mathcal {I}_{d}=[\mathcal
{I},\P_{1},\ldots,\P_{d}] \subset \mathcal {F}\langle
\bu_{1},\ldots,\bu_{d}\rangle \{\Y\}$. By Theorem \ref{th-inter2},
$\mathcal {I}_{d}$ is a prime $\delta$-ideal with dimension 0 and
the same order as $\mathcal {I}$.

Let $ \mathcal{I}_1=[\mathcal {I},\P_{0},\ldots,\P_{d}]=[\mathcal
{I}_{d},\P_0] \subset \mathcal {F}\langle
\bu_1,\ldots,\bu_d;u_{01},\ldots,u_{0n}\rangle
\{u_{00},y_{1},\ldots,y_{n}\}$.
From Lemma~\ref{lm-cf2}, $\mathcal {A}=\{f,S_fy_{1}-\frac{\partial
f}{\partial u_{01}^{(s)}},\, \ldots,\, S_fy_{n}-\frac{\partial
f}{\partial u_{0n}^{(s)}}\}$ is a characteristic set of
$\I_{\zeta,\xi}$. By Lemma~\ref{lem-d}, $u_{10},\ldots,u_{d0}$ is a
parametric set of $\I_{\zeta,\xi}$. So $\mathcal {A}$ is also a
characteristic set of $\mathcal{I}_1$ w.r.t. the elimination ranking
$u_{00}\prec y_{1} \prec \cdots \prec y_{n}$. Since $\dim(\mathcal
{I}_1)=0$, from Corollary~\ref{cor-0 order}, we have $\ord(\mathcal
{I}_1)=\ord(\mathcal {A})=\ord(f)$.

On the other hand, if $(\eta_{1},\ldots,\eta_{n})$ is a generic
point of $\mathcal {I}_{d}$, then $(\eta_{1},\ldots,\eta_{n},\zeta)$
is a generic point of $\mathcal{I}_1$ with $\zeta=-\sum_{j=1}^n
u_{0j}\eta_{j}$ and $\dim(\mathcal {I}_1)=0$. Clearly,
$u_{0k}(k=1,\ldots,n)$ are $\delta$-independent over $\mathcal
{F}\langle\bu_{1},\ldots,\bu_{d},\eta_{1},\ldots,\eta_{n}\rangle$.
Denote $\mathcal {F}\langle\bu_{1},\ldots,\bu_{d},u_{01},$
$\ldots,u_{0n}\rangle$ by $\ff_1$. So for a sufficiently large $t$,
\begin{eqnarray}
\quad &&\omega_{\mathcal {I}_1}(t) = \ord(\mathcal{I}_1)\nonumber \\
&=&\trdeg\,\ff_1(\eta_{i}^{(j)},\zeta^{(j)}:\,i=1,\ldots,n;j\leq
t)\big/\mathcal {F}_1\nonumber \\
&=&\trdeg\,\mathcal {F}_1(\eta_{i}^{(j)}:\,i=1,\ldots,n;j\leq
t)\big/\mathcal {F}_1\nonumber \\
&=&\trdeg\,\mathcal {F}\langle
\bu_{1},\ldots,\bu_{d}\rangle(\eta_{i}^{(j)}:\,i=1,\ldots,n;j\leq
t)\big/\mathcal {F}\langle\bu_{1},\ldots,\bu_{d}\rangle \nonumber
\\&=&\omega_{\mathcal {I}_{d}}(t)=\ord(\mathcal {I}_{d}) \nonumber
\end{eqnarray}
Thus, $\ord(\mathcal{I}_1)=\ord(\mathcal {I}_{d})=\ord(\mathcal
{I})$, and consequently, $\ord(\mathcal {I})=\ord(f)$.   \qed

As a consequence, we can give an equivalent definition for the order
of a prime $\delta$-ideal using Chow forms.
\begin{definition}
Let $\mathcal {I}$ be a prime $\delta$-ideal in $\mathcal {F}\{\Y\}$
with dimension $d$ and $F(\bu_{0},\bu_{1},$ $\ldots,$ $\bu_{d})$ its
Chow form. The order of $\mathcal {I}$ is defined to be the order of
its Chow form.
 \end{definition}

The following result shows that we can recover the generic point $(\xi_1,\ldots,\xi_n)$
of $V$ from its Chow form.
\begin{theorem} \label{th-generic}
Let $f(\bu;u_{00},\ldots,u_{d0})$ be defined as above and $h$ the
order of $V$. Then we have
 \[\xi_{\rho}=\overline{\frac{\partial f}{\partial u_{0\rho}^{(h)}}}\bigg/{\overline{S}_f},\rho=1,\ldots,n\]
where $\overline{\frac{\partial
f}{\partial u_{0\rho}^{(h)}}}$ and $\overline{S}_f$ are obtained by replacing $(u_{00},\ldots,u_{d0})$
by $(\zeta_{0},\ldots,\zeta_{d})$ in $ \frac{\partial f}{\partial
u_{0\rho}^{(h)}}$ and $ \frac{\partial f}{\partial u_{00}^{(h)}}$
respectively.
\end{theorem}
\proof It follows directly from
Lemma~\ref{lm-cf2} and Theorem \ref{th-choworder}. \qed

\begin{remark} \label{rm-resol}
The $\delta$-resolvent \cite{hubert} and \cite[p.34]{ritt} is
closely related with the differential Chow form, although they are
different. From Lemma \ref{lm-cf2} and Theorem~\ref{th-choworder},
we can see that both of them can be used to reduce a $\delta$-ideal
$\CI$ to a $\delta$-hypersurface which has the same order as $\CI$
in certain sense. But they are quite different. Firstly, the
resolvent depends on a parameter set $U$ of $\CI$. To be more
precise, let $Z=\Y\setminus U$. Then, the resolvent is essentially
constructed in $\ff\langle U\rangle\{Z\}$. Secondly, to define the
resolvent, we need only one new indeterminate $w$ and add one linear
$\delta$-polynomial $w - c_1v_1-\cdots-c_p v_p$ to $\CI$, where $Z =
\{v_1,\ldots,v_p\}$. Furthermore, $c_i$ in the above expression are
from $\ff$. Therefore, the resolvent will never be the Chow form.
%
%
Of course, if $\CI$ is of dimension zero and we take $c_i$ in $w -
c_1y_1-\cdots-c_n y_n$ as $\delta$-inderterminates, then the
resolvent is the Chow form of $\CI$. On the other hand, the
resolvent of $\CI$ can be obtained from its Chow form by
specializing some $u_{ij}$ to certain specific values and using
techniques in Theorem \ref{th-specil}.
\end{remark}

\subsection{Differential Chow form is differentially homogenous}
\label{sec-chow2} Following Kolchin \cite{kol4}, we introduce the
concept of $\delta$-homogenous $\delta$-polynomials.
\begin{definition} \label{d-homogenous}
A $\delta$-polynomial $p \in \mathcal
{F}\{y_{0},y_{1},\ldots,y_{n}\}$ is called $\delta$-homogenous of
degree $m$ if for a new $\delta$-$\ff$-indeterminate $\lambda$, we
have $p(\lambda y_{0},\lambda y_{1}\ldots,\lambda
y_{n})=\lambda^{m}p(y_{0},y_{1},$ $\ldots,y_{n}) $.
\end{definition}

The differential analog of Euler's theorem related to
homogenous polynomials is valid.
\begin{theorem}\cite[p.71]{kol} \label{th-dhomo}
A $\delta$-polynomial $f \in \mathcal
{F}\{y_{0},y_{1},\ldots,y_{n}\}$ is $\delta$-homogenous of degree
$m$ if and only if \newline \[\sum_{j=0}^{n} \sum_{k \in \mathbb{N}}
{k+r \choose r} y_{j}^{(k)} \frac{\partial
f(y_{0},\ldots,y_{n})}{\partial y_{j}^{(k+r)} } = \left\{
\begin{array}{ccc} mf & & r = 0 \\ 0 & & r \neq 0 \\ \end{array} \right.\]
\end{theorem}
For the Chow form, we have the following result.
\begin{theorem} \label{chowhomogenous}
Let $F(\bu_{0},\bu_{1},\ldots,\bu_{d})=f(\bu;u_{00},\ldots,u_{d0})$
be the Chow form of an irreducible $\delta$-variety $V$ of dimension
$d$ and order $h$. Then

1)  \[ \sum_{j=0}^{n}u_{\tau j}\frac{\partial f}{\partial u_{\sigma
j}} +\sum_{j=0}^{n}u_{\tau j}'\frac{\partial f}{\partial u_{\sigma
j}'}+\cdots+\sum_{j=0}^{n}u_{\tau j} ^{(h)}\frac{\partial
f}{\partial u_{\sigma j}^{(h)}} = \left\{ \begin{array}{lcl} 0&\quad& \sigma \neq \tau \\
rf &\quad& \sigma = \tau  \end{array} \right.  \] where $r$ is a
nonnegative integer.

2)  $F(\bu_{0},\bu_{1},\ldots,\bu_{d})$ is a $\delta$-homogenous
$\delta$-polynomial of degree $r$ in each set $\bu_{i}$ and $F$ is
of total degree $(d+1)r$.
\end{theorem}
\proof  Let $(\xi_{1},\ldots,\xi_{n})$ be a generic point of $V$ and
$\zeta_{i}=-\sum_{j=1}^{n}u_{ij}\xi_{j}\,(i=0,\ldots,d)$ defined in \bref{equ1}.
From \bref{eq-def}, $f(\bu;u_{00},\ldots,u_{d0})$ is the
characteristic set of the prime $\delta$-ideal $\mathbb{I}_\zeta$.
Since $f(\bu;\zeta_{0},\zeta_{1},\ldots,\zeta_{d})=0$, we have
\[ \begin{array}{*{5}{c@{\:+\:}}c@{\;=\;}cr}
\overline{\frac{\partial f}{\partial u_{\sigma j}}}&\frac{\partial
f}{\partial \zeta_{\sigma }}(-\xi_{j})&\frac{\partial f}{\partial
\zeta_{\sigma
}'}(-\xi_{j}')&\frac{\partial f}{\partial \zeta_{\sigma }''}(-{2 \choose 0}\xi_{j}'')&\ldots& \frac{\partial f}{\partial \zeta_{\sigma }^{(h)}}[-{h \choose 0}\xi_{j}^{(h)}]&0&(0*)\\
\overline{\frac{\partial f}{\partial u_{\sigma
j}'}}&0&\frac{\partial f}{\partial \zeta_{\sigma
}'}(-\xi_{j})&\frac{\partial f}{\partial \zeta_{\sigma }''}(-{2 \choose 1}\xi_{j}')&\ldots& \frac{\partial f}{\partial \zeta_{\sigma }^{(h)}}[-{h \choose 1}\xi_{j}^{(h-1)}]&0&(1*)\\
\overline{\frac{\partial f}{\partial u_{\sigma
j}''}}&0&0&\frac{\partial f}{\partial\zeta_{\sigma }''}(-{2 \choose
2}\xi_{j})&\ldots&\frac{\partial f}{\partial \zeta_{\sigma
}^{(h)}}[-{h \choose 2}\xi_{j}^{(h-2)}]&0&(2*)\\
\multicolumn {7}{c}{\dotfill} \\
\overline{\frac{\partial f}{\partial u_{\sigma j}^{(h)}}}&0&0&0&0&\frac{\partial f}{\partial \zeta_{\sigma }^{(h)}}[-{h \choose h }\xi_{j}^{(0)}]&0&(h*)\\
\end{array}
\]
In the above equations,    $\overline{\frac{\partial f}{\partial
u_{\sigma j}^{(l)}}}$ and $\frac{\partial f}{\partial \zeta_{\sigma
}^{(l)}}$ $(l=0,\ldots,h; j=1,\ldots,n)$ are  obtained by
substituting $\zeta_{i}$ to $u_{i0}\,(i=0, 1,    \ldots,    d)$ in
each $\frac{\partial f}{\partial u_{\sigma j}^{(l)}}$ and
$\frac{\partial f}{\partial u_{\sigma 0}^{(l)}}$ respectively.

Now, let us consider the $\delta$-polynomial $\sum_{j=0}^n
\sum_{k\geq 0} {k+i \choose k}u_{\sigma j}^{(k)}\frac{\partial
f}{\partial u_{\sigma j}^{(k+i)}}$.

 In the case $i=0$, firstly, let $(0*)\times u_{\tau j}
+(1*)\times u_{\tau j}'+\cdots+(h*)\times u_{\tau j}^{(h)}$ and
add them together for $j$ from 1 to $n$. We obtain
\[ \sum_{j=1}^{n} u_{\tau j}\overline{\frac{\partial f}{\partial u_{\sigma
j}}}+\sum_{j=1}^{n} u_{\tau j}' \overline{\frac{\partial f}{\partial
u_{\sigma j}'}}+\cdots+\sum_{j=1}^{n} u_{\tau j}^{(h)}
\overline{\frac{\partial f}{\partial u_{\sigma
j}^{(h)}}}+\zeta_{\tau} \frac{\partial f}{\partial \zeta_{\sigma
}}+\zeta_{\tau}' \frac{\partial f}{\partial \zeta_{\sigma
}'}+\cdots+\zeta_{\tau}^{(h)} \frac{\partial f}{\partial
\zeta_{\sigma }^{(h)}}=0.\]
 So  the  $\delta$-polynomial $\sum\limits^{n}_{j=0} u_{\tau j}\frac{\partial f}{\partial
u_{\sigma j}}+\sum\limits^{n}_{j=0} u_{\tau j}' \frac{\partial
f}{\partial u_{\sigma j}'}+\sum\limits^{n}_{j=0} u_{\tau j}''
\frac{\partial f}{\partial u_{\sigma
j}''}+\cdots+\sum\limits^{n}_{j=0} u_{\tau j}^{(h)}$ $
\frac{\partial f}{\partial u_{\sigma j}^{(h)}} $ vanishes at
$(u_{00},\ldots,u_{d0})=(\zeta_{0},\ldots,\zeta_{d})$. Thus in the
case $\tau=\sigma$, it can be divisible by $f$, i.e.
$\sum\limits_{j=0}^n \sum\limits_{k=0}^h u_{\sigma j}^{(k)}$ $
\frac{\partial f}{\partial u_{\sigma j}^{(k)}}$ $=rf$. By Euler's
theorem, $f$ is an algebraic homogenous $\delta$-polynomial of
degree $r$ in each set of indeterminates
$\bu_{i}=(u_{i0},\ldots,u_{in})$ and their derivatives. But in the
case $\tau \neq \sigma$, since this $\delta$-polynomial is of order
not greater than $f$ and can not be divisible by $f$, by
Lemma~\ref{le-divisible}, it must be identically zero. Thus, we have
proved 1) of the Theorem.

In  the case $i \neq 0$,
{\tiny
\begin{eqnarray} 0 &=& (i*)\times {i \choose
i} u_{\sigma j} +(i+1*)\times {i+1 \choose i} u_{\sigma
j}'+\cdots+(h*)\times {h \choose i} u_{\sigma j}^{(h-i)} \nonumber \\
 &=& {i \choose i}u_{\sigma j}\overline{\frac{\partial f}{\partial u_{\sigma
j}^{(i)}}}+{i+1 \choose i}u_{\sigma j}'\overline{\frac{\partial
f}{\partial u_{\sigma j}^{(i+1)}}}+\cdots+{h \choose i}u_{\sigma
j}^{(h-i)}\overline{\frac{\partial
f}{\partial u_{\sigma j}^{(h)}}} \nonumber \\
 &\quad&+\frac{\partial f}{\partial\zeta
_{\sigma}^{(i)}}\bigg(-{i \choose i}u_{\sigma j} \xi_{j}\bigg)
 +\frac{\partial f}{\partial \zeta _{\sigma}^{(i+1)}}\bigg(-{i+1 \choose
i}u_{\sigma j} \xi_{j}'-{i+1 \choose i}{i+1 \choose i+1}u_{\sigma
j}' \xi_{j}\bigg)   \nonumber \\     &\quad&        + \cdots \cdots\nonumber \\
&\quad &+\frac{\partial f}{\partial \zeta _{\sigma}^{(h)}}\bigg(-{i
\choose i}{h \choose i}u_{\sigma j} \xi_{j}^{(h-i)}-{i+1 \choose
i}{h \choose i+1}u_{\sigma j}' \xi_{j}^{(h-i-1)}-\cdots-{h \choose
i}{h \choose h}u_{\sigma j}^{(h-i)} \xi_{j}\bigg) \nonumber \\
&=& {i \choose i}u_{\sigma j}\overline{\frac{\partial f}{\partial
u_{\sigma j}^{(i)}}}+{i+1 \choose i}u_{\sigma
j}'\overline{\frac{\partial f}{
\partial u_{\sigma j}^{(i+1)}}}+\cdots+{h \choose i}u_{\sigma
j}^{(h-i)}\overline{\frac{\partial
f}{\partial u_{\sigma j}^{(h)}}}\nonumber \\
&\quad& +{i \choose i}\frac{\partial f}{\partial \zeta
_{\sigma}^{(i)}}(-u_{\sigma j} \xi_{j})+{i+1 \choose
i}\frac{\partial f}{\partial \zeta _{\sigma}^{(i+1)}}(-u_{\sigma j}
\xi_{j})'+\cdots+{h \choose i}\frac{\partial f}{\partial\zeta _{
\sigma}^{(h)}}(-u_{\sigma j} \xi_{j})^{(h-i)} \nonumber
\end{eqnarray}}

Therefore, $\sum_{j=1}^n {i \choose i}u_{\sigma j}\overline{\frac{\partial
f}{\partial u_{\sigma j}^{(i)}}}+\sum_{j=1}^n {i+1 \choose
i}u_{\sigma j}'\overline{\frac{\partial f}{\partial u_{\sigma
j}^{(i+1)}}}+\cdots+\sum_{j=1}^n {h \choose i}u_{\sigma
j}^{(h-i)}\overline{\frac{\partial f}{\partial u_{\sigma
j}^{(h)}}}+{i \choose i}\zeta_{\sigma}\frac{\partial f}{\partial
\zeta_{\sigma}^{(i)}}$ $+{i+1 \choose i}$
$\zeta_{\sigma}'\frac{\partial f}{\partial
\zeta_{\sigma}^{(i+1)}}+\cdots+{h \choose
i}\zeta_{\sigma}^{(h-i)}\frac{\partial f}{\partial
\zeta_{\sigma}^{(h)}}=0$.

Thus,  the $\delta$-polynomial $\sum_{j=0}^n {i \choose i}u_{\sigma
j}\frac{\partial f}{\partial u_{\sigma j}^{(i)}}+\sum_{j=0}^n {i+1
\choose i}u_{\sigma j}'\frac{\partial f}{\partial u_{\sigma
j}^{(i+1)}}+\cdots+\sum_{j=0}^n {h \choose i}$ $ u_{\sigma
j}^{(h-i)}\frac{\partial f}{\partial u_{\sigma j}^{(h)}}$ is
identically zero, for it vanishes at
$(u_{00},\ldots,u_{d0})=(\zeta_{0},\ldots,\zeta_{d})$ and can not be
divisible by $f$.

From the above, we conclude that \[ \sum_{j=0}^n \sum_{k\geq 0} {k+i
\choose i}u_{\sigma j}^{(k)}\frac{\partial f}{\partial u_{\sigma
j}^{(k+i)}}=\left\{
\begin{array}{ccc} 0&\quad& i\neq 0 \\ rf&\quad& i=0
\end{array} \right.\]
From Theorem~\ref{th-dhomo} and the symmetry property of
$F(\bu_{0},\ldots,\bu_{d})$, the theorem is obtained. \qed

Lemma~\ref{le-sym}, Theorem \ref{th-choworder}, and
Theorem~\ref{chowhomogenous} together prove the first property of
Theorem \ref{th-main2}.

\begin{remark}\label{rm-phom}
Using the terminology from \cite{kol51}, the differential Chow form
$F(\bu_{0},\bu_{1},$ $\ldots,$ $\bu_{d})$ is
$(d+1)$-$\delta$-homogenous in $(\bu_{0},\bu_{1},\ldots,\bu_{d})$.
\end{remark}

\begin{definition}\label{def-ddeg}
Let $\mathcal {I}$ be a prime $\delta$-ideal in $\mathcal {F}\{\Y\}$
of dimension $d$ and $F(\bu_{0},\bu_{1},$ $\ldots,$ $\bu_{d})$  its
Chow form. The {\em $\delta$-degree} of $\mathcal {I}$ is defined to
be the $\delta$-homogenous degree $r$ of its Chow form in each
$\bu_{i}$ $(i=0,\ldots,d)$.
\end{definition}

The following result shows that the $\delta$-degree of a
$\delta$-variety $V$ is an invariant of $V$ under invertible linear
transformations.

\begin{lemma}\label{lm-cfl}
Let $A=(a_{ij})$ be an $n\times n$ invertible matrix with
$a_{ij}\in\mathcal {F}$ and $F(\bu_0,\bu_1,\ldots,\bu_d)$ the Chow
form of an irreducible $\delta$-variety $V$ with dimension $d$. Then
the Chow form of the image $\delta$-variety of $V$ under the linear
transformation $\Y=A\X$ is
$F^{A}(\bv_{0},\ldots,\bv_{d})=F(\bv_{0}B,\ldots,\bv_{d}B)$, where
$B=\left(\begin{array}{cc} 1 & 0 \ldots 0 \\{0 \atop {\vdots \atop
0}}& A
\end{array}\right)$ and $\bu_{i}$ and $\bv_{i}$ are
regarded as row vectors.
\end{lemma}
\proof Let $\xi=(\xi_{1},\ldots,\xi_{n})$ be a generic point of $V$.
Under the linear transformation $\Y=A\X$, $\xi$ is mapped to
$\eta=(\eta_{1},\ldots,\eta_{n})$ with $\eta_i=\sum_{j=1}^n
a_{ij}\xi_j$. Under this transformation $V$ is mapped to a
$\delta$-variety $V^A$ whose generic point is $\eta$.
Note that
$F^{A}(\bv_{0},\ldots,\bv_{d})=f^{A}(v_{ij};v_{00},\ldots,v_{d0})=f(\sum_{k=1}^n
v_{ik}a_{kj};v_{00},\ldots,v_{d0})$
 and $f^{A}(v_{ij};-\sum_{k=1}^n v_{0k}\eta_{k},\ldots,-\sum_{k=1}^n v_{dk}\eta_{k})=f(\sum_{k=1}^n v_{ik}a_{kj};$ $-\sum_{k=1}^n v_{0k}\eta_{k},\ldots,$ $-\sum_{k=1}^n v_{dk}\eta_{k})
 $ $=f(\sum_{k=1}^n v_{ik}a_{kj};-\sum_{j=1}^n(\sum_{k=1}^n
v_{0k}a_{kj})\xi_{j},\ldots,-\sum_{j=1}^n(\sum_{k=1}^n v_{dk}$
\newline $a_{kj})\xi_{j})=0.$ Since $V^A$ is of the same dimension
and order as $V$  and  $F^{A}$ is irreducible, from the definition
of the Chow form, the claim is proved.\qed

\begin{definition} \label{denomination}
Let  $p$ be a  $\delta$-polynomial in $\mathcal {F}\{\Y\}$. Then the
smallest number $r$ such that $y_{0}^r
p(y_{1}/y_{0},\ldots,y_{n}/y_{0})\in \mathcal
{F}\{y_{0},y_{1},\ldots,y_{n}\} $ is called the denomination of $p$,
which is denoted by $\den(p)$.
\end{definition}

\begin{example} In the case $d=n-1$ and $n>1$, if $\{p(y_{1},\ldots,y_{n})\}$ is a
characteristic set of $\mathcal {I}$ w.r.t. any ranking, then by
Example \ref{ex-cf1} the $\delta$-degree of $\mathcal {I}$ cannot
exceed the denomination of $p$. So  the denomination of $p$ gives an
upper bound of $\delta$-degree of $\CI$. But, we do not know whether
they are the same.
\end{example}

\begin{example}\label{ex-n12}
 In the case  $n=1$, we have $d=0$. If $\{p(y)\}$ is a
characteristic set of $\mathcal {I}$ w.r.t. any ranking, then the
$\delta$-degree of $\mathcal {I}$ is exactly equal to the
denomination  of $p$. Now we can give a bound for the
$\delta$-degree of $\mathcal {I}$ from the original equation of $p$
without computing its denomination.

For a $\delta$-monomial
$\theta(y)=y^{l_{0}}(y')^{l_{1}}(y'')^{l_{2}}\ldots
(y^{(s)})^{l_{s}}$, define its {\em weighted degree} to be
$l_{0}+2l_{1}+\cdots+(s+1)l_{s}$, denoted by
$\wdeg(\theta(y))=l_{0}+2l_{1}+\cdots+(s+1)l_{s}$. For a
$\delta$-polynomial $p \in \mathcal {F}\{y\}$, we can define its
weighted degree to be the maximum of the  weighted degrees of all
the $\delta$-monomials effectively appearing in $p$. Clearly, the
denomination of $p$ cannot exceed its  weighted degree. And we have
examples for which $\den(p)<\wdeg(p)$. Let $p=2y'^2-yy''$. Then,
$\den(p)=3$ and $\wdeg(p)=4$. The Chow form of $\sat(p)$ is
$F(\bu)=u_{0}u_{1}u_{0}''-u_{0}^2 u_{1}''-2u_{0}'^2
u_{1}+2u_{0}u_{0}'u_{1}'$, where $\bu=(u_{0},u_{1})$. So the
$\delta$-degree of $\mathcal{I}=\sat(p)$ is 3 which is less than
$\wdeg(p)$.

Now we contrast the above $p$ with a $\delta$-polynomial
$q=y'^2-yy''$ that is different from $p$ by only a single
coefficient. Then $\den(q)=4$ and $\wdeg(q)=4$. The Chow form of
$\sat(q)$ is $G(\bu_0)=u_{0}u_{1}^2u_{0}''-u_{0}^2
u_1u_{1}''-u_{0}'^2 u_{1}^2+u_{0}^2u_{1}'^2$, so the $\delta$-degree
of $\sat(q)$ is 4, which is equal to the weighted degree of $q$.
Thus, the weighted degree is a sharp bound.
%
\end{example}

%
%
%

\subsection{Factorization of differential Chow form }
\label{sec-fac} In the algebraic case, the Chow form can be factored
into the product of linear polynomials with the generic points of
the variety as coefficients. In this section, we will show that
there is a differential analog to this result.

Let
$$\tilde{\bu}=\cup_{i=0}^d \bu_i\setminus \{u_{00}\}$$ and
$\mathcal {F}_{0}=\mathcal
{F}\langle\tilde{\bu}\rangle(u_{00}^{(0)},\ldots,u_{00}^{(h-1)})$.
Consider the Chow form $f(\bu;u_{00},u_{10},\ldots,$ $u_{d0})$ as an
irreducible algebraic polynomial $p(u_{00}^{(h)})$ in $\mathcal
{F}_{0}[u_{00}^{(h)}]$. Let
$g=\deg(p,u_{00}^{(h)})=\deg(f,u_{00}^{(h)})$.
In a suitable algebraic extension field of $\mathcal {F}_{0}$,
$p(u_{00}^{(h)})=0$ has $g$ roots $\gamma_{1},\ldots,\gamma_{g}$.
Thus
\begin{equation}\label{eq-fac00}
f(\bu;u_{00},u_{10},\ldots,u_{d0})=A(\bu_{0},\bu_{1},\ldots,\bu_{d})\prod^g_{\tau=1}(u_{00}^{(h)}-\gamma_{\tau})
\end{equation}
where  $A(\bu_{0},\bu_{1},\ldots,\bu_{d})$ is in $\mathcal
{F}\{\bu_0,\ldots,\bu_d\}$.
For each $\tau$ such that $1\le \tau\le g$, let
\begin{equation}\label{eq-ftau}
 \mathcal{F}_\tau=\mathcal {F}_0(\gamma_{\tau})=\mathcal
{F}\langle\tilde{\bu}\rangle(u_{00}^{(0)},\ldots,u_{00}^{(h-1)},\gamma_{\tau})
\end{equation}
be an algebraic extension of $\mathcal {F}_0$ defined by
$p(u_{00}^{(h)})=0$.
We will define a derivation $\delta_{\tau}$ on $\mathcal{F}_\tau$ so
that $(\mathcal{F}_\tau,\delta_\tau)$ becomes a differential field.
This can be done in a very natural way. For $e\in \mathcal
{F}\langle\tilde{\bu}\rangle$, define $\delta_{\tau} e = \delta
e=e'$. Define $\delta_{\tau}^{i} u_{00} = u_{00}^{(i)}$ for
$i=0,\ldots,h-1$ and $$\delta_{\tau}^{h} u_{00} = \gamma_{\tau}.$$
Since $f$, regarded as an algebraic polynomial $p$ in
$u_{00}^{(h)}$, is a minimal polynomial of $\gamma_{\tau}$,
$\sep_{f}=\frac{\partial f}{\partial u_{00}^{(h)}}$ does not vanish
at $u_{00}^{(h)}=\gamma_{\tau}$. Now, we define the derivatives of
$\delta_{\tau}^{i} u_{00}$ for $i> h$ by induction. Firstly, since
$p(\gamma_{\tau})=0$, $\delta_\tau
(p(\gamma_{\tau}))=\sep_{f}\big|_{u_{00}^{(h)}=\gamma_\tau}\delta_\tau
(\gamma_{\tau})+T\big|_{u_{00}^{(h)}=\gamma_\tau}=0$, where
$T=f'-\sep_fu_{00}^{(h)}$. We define $\delta_{\tau}^{h+1} u_{00}$ to
be $\delta_\tau
(\gamma_{\tau})=-\frac{T}{\sep_{f}}\Big|_{u_{00}^{(h)}=\gamma_{\tau}}$.
Supposing the derivatives of $\delta_{\tau}^{j} u_{00}$ with order
less than $j<i$ have been defined, we now  define $\delta_{\tau}^{i}
u_{00}$. Since $f^{(i)}=\sep_{f}u_{00}^{(h+i)}+T_{i}$ is linear in
$u_{00}^{(h+i)}$, we define $\delta_{\tau}^{i} u_{00}$ to be
$-\frac{T_{i}}{\sep_{f}}\Big|_{u_{00}^{(h+j)}=\delta_\tau^{h+j}u_{00},
j <i}$.
%

In this way, $(\ff_\tau,\delta_\tau)$ is a differential field which
can be considered as a finitely differential extension field of
$(\ff\langle\tilde{\bu}\rangle,\delta)$. Recall that
$\ff\langle\tilde{\bu}\rangle$ is a finitely $\delta$-extension
field of $\ff$ contained in $\ee$. By the definition of universal
$\delta$-extension field, there exists a $\delta$-extension field
$\ff^*\subset\ee$ of $\ff\langle\tilde{\bu}\rangle$ and a
$\delta$-isomorphism $\varphi_\tau$ over
$\ff\langle\tilde{\bu}\rangle$ from $(\ff_\tau,\delta_\tau)$ to
$(\ff^*,\delta)$. Summing up the above results, we have

\begin{lemma}\label{lm-ftau}
$(\ff_\tau, \delta_\tau)$ defined above is a finitely differential
extension field of  $\ff\langle\tilde{\bu}\rangle$, which is
$\delta$-$\ff\langle\tilde{\bu}\rangle$-isomorphic to a subfield of
$\ee$.
\end{lemma}

Let $p$ be a $\delta$-polynomial in
$\mathcal{F}\{\bu_0,\bu_1,\ldots,\bu_d\} =
\ff\{\tilde{\bu},u_{00}\}$.
For convenience,  by the symbol
$p\Big|_{u_{00}^{(h)}=\gamma_{\tau}}$ or saying replacing
$u_{00}^{(h)}$ by $\gamma_{\tau}$, we mean substituting
$u_{00}^{(h+i)}$ by $\delta_\tau^{i}\gamma_{\tau}\,(i\geq0)$ in $p$.
Similarly, by saying $p$ vanishes at $u_{00}^{(h)}=\gamma_{\tau}$,
we mean $p\Big|_{u_{00}^{(h)}=\gamma_{\tau}}=0$. It is easy to prove
the following lemma.
%
%

\begin{lemma}\label{lm-ftaup}
Let $p$ be a $\delta$-polynomial in $\ff\{\tilde{\bu},u_{00}\}$ and
$f$ the Chow form for a prime $\delta$-$\ff$-ideal $\CI$ of
dimension $d$. Then $p\in\sat(f)$ if and only if $p$ vanishes at
$u_{00}^{(h)}=\gamma_{\tau}$.
\end{lemma}

When a $\delta$-polynomial $h\in\ff\langle\tilde{\bu}\rangle\{\Y\}$
vanishes at a point $\eta\in\ff_\tau^n$, it is easy to see that $h$
vanishes at $\varphi_\tau(\eta)\in\ee^n$. For convenience, by saying
$\eta$ is in a $\delta$-variety $V$ over
$\ff\langle\tilde{\bu}\rangle$, we mean $\varphi_\tau(\eta)\in V$.
%

\begin{remark}\label{rm-ftau}
In order to make $\ff_\tau$ a differential field, we need to
introduce a differential operator $\delta_\tau$ which is related to
$\gamma_\tau$ and there does not exist a single differential
operator to make all $\ff_\tau(\tau=1,\ldots,g)$ differential
fields.
This natural phenomenon related with nonlinear differential
equations seems not used before.
For instance, let $p=y'^2-4y$. Then $\CI=\sat(p)$ is a prime
$\delta$-ideal in $\Q(t)\{y\}$ and let $\overline{\ff}$ be the
differential rational function field of $\CI$.
By factoring $p = (y' -2\sqrt{y})(y' +2\sqrt{y})$, we can define two
more differential fields: $\ff_1=\Q(t)(\sqrt{y})$ with a
differential operator $\delta_1 y = 2\sqrt{y}$ and
$\ff_2=\Q(t)(\sqrt{y})$ with a differential operator $\delta_2 y =
-2\sqrt{y}$. Note that $\ff_1$ and $\ff_2$ are not compatible,
although each of them is isomorphic to a subfield of $\ee$. Finally,
both $\ff_1$ and $\ff_2$ are isomorphic to  $\overline{\ff}$.
\end{remark}

With these preparations, we now give the following theorem.
\begin{theorem}\label{th-fac1}
Let  $F(\bu_{0},\bu_{1},\ldots,\bu_{d})=f(\bu;u_{00},\ldots,u_{d0})$
 be the Chow form of an irreducible $\delta$-$\ff$-variety of
dimension $d$ and order $h$. Then, there exist $\xi_{\tau
1},\ldots,\xi_{\tau n}$ in a differential extension field
$(\ff_\tau,\delta_\tau)$ ($\tau=1,\ldots,g$)  of $(\mathcal
{F}\langle\tilde{\bu}\rangle,\delta)$  such that
\begin{eqnarray}\label{eq-fac10}
F(\bu_{0},\bu_{1},\ldots,\bu_{d})=A(\bu_{0},\bu_{1},\ldots,\bu_{d})\prod^g_{\tau=1}(u_{00}+\sum_{\rho=1}^n
u_{0\rho}\xi_{\tau \rho})^{(h)}
\end{eqnarray}
where  $A(\bu_{0},\bu_{1},\ldots,\bu_{d})$ is in $\mathcal
{F}\{\bu_0,\ldots,\bu_d\}$, $\tilde{\bu}=\cup_{i=0}^d\bu_i\backslash
u_{00}$ and $g=\deg(f,u_{00}^{(h)})$.
Note that equation \bref{eq-fac10} is formal and should be
understood in the following precise meaning:
$(u_{00}+\sum_{\rho=1}^n u_{0\rho}\xi_{\tau \rho})^{(h)}
\stackrel{\triangle}{=}
\delta^{h}u_{00}+\delta_\tau^h(\sum_{\rho=1}^n u_{0\rho}\xi_{\tau
\rho}).$
\end{theorem}
\proof We will follow the notations introduced in the proof of Lemma
\ref{lm-ftau}. Since $f$ is irreducible, we have $f_{\tau
0}=\frac{\partial f}{\partial u_{00}^{(h)}}\Big|
_{u_{00}^{(h)}=\gamma_{\tau}} \neq 0$.
Let $\xi_{\tau\rho}=f_{\tau \rho}\big/f_{\tau0}(\rho=1,\ldots,n)$,
where $f_{\tau \rho}=\frac{\partial f}{\partial u_{0
\rho}^{(h)}}\Big| _{u_{00}^{(h)}=\gamma_{\tau}}$. Note that $f_{\tau
\rho}$ and  $\xi_{\tau\rho}$ are in $\ff_\tau$.
We will prove \[\gamma_{\tau}=-\delta_\tau^h(u_{01}\xi_{\tau
1}+u_{02}\xi_{\tau2}+\cdots+u_{0n}\xi_{\tau n}).\]
Differentiating the equality
$f(\bu;\zeta_{0},\zeta_{1},\ldots,\zeta_{d})=0$ w.r.t.
$u_{0\rho}^{(h)}$, we have
\[\overline{\frac{\partial f}{\partial u_{0\rho}^{(h)}}}+
 \overline{\frac{\partial f}{\partial u_{00}^{(h)}}}(-\xi_{\rho})=0,\]
where $\overline{\frac{\partial f}{\partial u_{0\rho}^{(h)}}}$ are
obtained by substituting $\zeta_{i}$ to $u_{i0}\,(i=0, 1, \ldots,
d)$ in $\frac{\partial f}{\partial u_{ 0\rho}^{(h)}}$.
Multiplying $u_{0\rho}$ to the above equation and for $\rho$ from 1
to $n$, adding them together, we have
\[\sum_{\rho=1}^n u_{0\rho}\overline{\frac{\partial f}{\partial u_{0\rho}^{(h)}}}+
 \overline{\frac{\partial f}{\partial u_{00}^{(h)}}}(-\sum_{\rho=1}^n u_{0\rho}\xi_{\rho})
 =\sum_{\rho=1}^n u_{0\rho}\overline{\frac{\partial f}{\partial u_{0\rho}^{(h)}}}
  +\zeta_0 \overline{\frac{\partial f}{\partial u_{00}^{(h)}}}=0.\]
Thus, $q=\sum_{\rho=1}^n u_{0\rho}\frac{\partial f}{\partial
u_{0\rho}^{(h)}}+u_{00}\frac{\partial f}{\partial u_{00}^{(h)}} \in
\,\sat(f)$. Since $q$ is of order not greater than $f$, it must be
divisible by $f$. Since $q$ and $f$ have the same degree, there
exists an $a\in\mathcal {F}$ such that $q=af$. Setting
$u_{00}^{(h)}=\gamma_{\tau}$ in both sides of $q=af$, we have
$\sum_{\rho=1}^n u_{0\rho}f_{\tau \rho}+u_{00}f_{\tau 0}=0$. Hence,
as an algebraic equation, we have
\begin{equation}\label{eq-ff2}
u_{00}+\sum_{\rho=1}^n u_{0\rho}\xi_{\tau \rho}=0
\end{equation}
under the constraint $u_{00}^{(h)}=\gamma_\tau$. Equivalently, the
above equation is valid in $(\ff_\tau,\delta_\tau)$.
As a consequence, $\gamma_{\tau}=-\delta_\tau^{h}(\sum_{\rho=1}^n
u_{0\rho}\xi_{\tau \rho})$. Substituting them into
equation~\bref{eq-fac00}, the theorem is proved.\qed

\begin{remark} The factors in equation (\ref{eq-fac10}) are the
$h$-th derivatives of the factors in the corresponding factorization
of the algebraic Chow form \cite[p.37]{hodge}.
\end{remark}

For an element $\eta=(\eta_1,\ldots,\eta_n)$, denote its truncation
up to order $k$ as
$\eta^{[k]}=(\eta_1,\ldots,\eta_n,\ldots,\eta_1^{(k)},$
$\ldots,\eta_n^{(k)})$.

In the proof of Theorem \ref{th-fac1}, some equations are valid in
the algebraic case only. To avoid confusion, we introduce the
following notations:
\begin{equation}\label{eq-apoly}
 \begin{array}{l}
^a\P^{(0)}_0 =\,^a\P_0:=u_{00}+u_{01}y_1+\cdots+u_{0n}y_n \\
^a\P^{(1)}_0 =\,
^a\P'_0:=u_{00}'+u_{01}'y_1+u_{01}y'_1+\cdots+u_{0n}'y_n+u_{0n}y'_n
\\ \cdots \\
 ^a\P^{(s)}_0:=u_{00}^{(s)} + \sum_{j=1}^n \sum_{k=0}^s {s \choose k}u_{0j}^{(k)}y_j^{(s-k)}\end{array}
\end{equation}
which are considered to be algebraic polynomials in $\mathcal
{F}(\bu_{0}^{[s]},\ldots,\bu_{n}^{[s]})[\Y^{[s]}]$, and
$u_{ij}^{(k)}, y_i^{(j)}$ are treated as algebraic indeterminates.
A point $\eta=(\eta_1,\ldots,\eta_n)$ is said to be lying on
$^a\P^{(k)}_0$  if regarded as an algebraic point, $\eta^{[k]}$
%
is a zero of $^a\P^{(k)}_0$. As a consequence of \bref{eq-ff2} in
the proof of Theorem \ref {th-fac1}, we have

\begin{cor}\label{co-f1}
$(\xi_{\tau_1},\ldots,\xi_{\tau n},\ldots,\delta_\tau^{(h-1)}\xi_{\tau
1},\ldots,\delta_\tau^{(h-1)}\xi_{\tau n})\,(\tau=1,\ldots,g)$ are
common zeros of \,$^a\P_0=0,\,^a\P'_0$ $=0,\ldots,$
$^a\P^{(h-1)}_0=0$
 where $\delta_{\tau}$ is defined in
Lemma \ref{lm-ftau}.
\end{cor}

\begin{example}\label{ex-ff1} Continue from Example~\ref{example1}. In this example,
$F(\bu)=u_{01}^2 (u_{00}')^2-2u_{01}u_{01}'u_{00}u_{00}'$ $+(u_{01}')^2
u_{00}^2+4u_{01}^3u_{00}$, so $g=2$. And
$F(\bu)=u_{01}^2 (u_{00}'-u_{01}'\frac{u_{00}}{u_{01}}+2\sqrt{-1}\sqrt{u_{00}u_{01}})
               (u_{00}'-u_{01}'\frac{u_{00}}{u_{01}}-2\sqrt{-1}\sqrt{u_{00}u_{01}})$.
So we can obtain $\gamma_\tau=u_{01}'\frac{u_{00}}{u_{01}}\mp
2\sqrt{-1}\sqrt{u_{00}u_{01}}\,(\tau=1,2)$. Following the  proof of
Theorem~\ref{th-fac1}, we obtain that $\xi_{11}=-u_{00}/u_{01}$ with
$u_{00},u_{01}$ satisfying the relation $\delta_1
{u_{00}=\gamma_1=\frac{u_{00}}{u_{01}}
u_{01}'-2\sqrt{-1}\sqrt{u_{00}u_{01}}}$, and
$\xi_{21}=-u_{00}/u_{01}$ with $u_{00},u_{01}$ satisfying
$\delta_2{u_{00}=\gamma_2=\frac{u_{00}}{u_{01}}
u_{01}'+2\sqrt{-1}\sqrt{u_{00}u_{01}}}$.
In other words, $\xi_{11}$ and $\xi_{21}$ are in $\ff_1$ and $\ff_2$
respectively.
Thus, $F(\bu)$ has the factorization $F(\bu)=u_{01}^2 (u_{00}'+\xi_{11}u_{01}'+2\sqrt{-1}\sqrt{u_{00}u_{01}})
               (u_{00}'+\xi_{21}u_{01}'-2\sqrt{-1}\sqrt{u_{00}u_{01}})
               =u_{01}^2(u_{00}+\xi_{11}u_{01})'(u_{00}+\xi_{21}u_{01})'$.
Note that both $\xi_{i1}\,(i=1,2)$  satisfy $^a\P_{0} = u_{00}+
u_{01}\xi_{i1}=0$, but
$^a\P_{0}^{(1)} = u_{00}' + u_{01}'\xi_{i1} + u_{01}\delta_i
\xi_{i1}\ne0$.
\end{example}

\begin{lemma} \label{lm-free}
In equation \bref{eq-fac10}, $A(\bu_{0},\bu_{1},\ldots,\bu_{d})$ is free of
$u_{0i}^{(h)}\,(i=1,\ldots,n)$.
\end{lemma}
\proof Since $f$ is homogenous in the indeterminates $u_{0i}$ and
its derivatives up to the order $h$, we have
\[\sum_{\rho=1}^n u_{0\rho}^{(h)} \frac{\partial f}{\partial
u_{0\rho}^{(h)}}+u_{00}^{(h)} \frac{\partial f}{\partial
u_{00}^{(h)}}+\sum_{k=0}^{h-1}\sum_{\rho=0}^n u_{0\rho}^{(k)}
\frac{\partial f}{\partial u_{0\rho}^{(k)}}=rf,r \in \mathbb{N}.\]

In this equation, let $u_{00}^{(h)}=\gamma_{\tau}$, we obtain
\[\sum_{\rho=1}^n u_{0\rho}^{(h)} f_{\tau\rho}+\gamma_{\tau} f_{\tau 0}+\sum_{k=0}^{h-1}\sum_{\rho=0}^n u_{0\rho}^{(k)}
\widehat{\frac{\partial f}{\partial u_{0\rho}^{(k)}}}=0.\]
 where
$\widehat{\frac{\partial f}{\partial u_{0\rho}^{(k)}}}$ means
replacing $u_{00}^{(h)}$ by $\gamma_{\tau}$ in $\frac{\partial
f}{\partial u_{0\rho}^{(k)}}$.
Consequently, \[\gamma_{\tau}=-\sum_{\rho=1}^n
u_{0\rho}^{(h)}\xi_{\tau\rho}-\sum_{k=0}^{h-1}\sum_{\rho=0}^n
u_{0\rho}^{(k)} \widehat{\frac{\partial f}{\partial
u_{0\rho}^{(k)}}} \bigg/f_{\tau0}.
\]
Hence,
\begin{equation}\label{eq-fac11}
f(\bu;u_{00},u_{10},\ldots,u_{d0})=A\prod^g_{\tau=1}\Big(u_{00}^{(h)}+
 \sum_{\rho=1}^n u_{0\rho}^{(h)}\xi_{\tau\rho}+\sum_{k=0}^{h-1}\sum_{\rho=0}^n
u_{0\rho}^{(k)}  \widehat{\frac{\partial f}{\partial
u_{0\rho}^{(k)}}}\bigg/f_{\tau0}\Big).
\end{equation}
%

We claim that $\xi_{\tau \rho}$
and $\sum_{k=0}^{h-1}\sum_{\rho=0}^n u_{0\rho}^{(k)} $
$\widehat{\frac{\partial f}{\partial
u_{0\rho}^{(k)}}}\big/f_{\tau0}$ are algebraically independent  of
$u_{0i}^{(h)}\,(i=1,\ldots,n)$. Firstly, since $\xi_{\rho}$ is
algebraically independent of $u_{0i}^{(h)}(i=1,\ldots,n)$, and by Theorem~\ref{th-generic},
\[\xi_{\rho}=\frac{\partial f}{\partial
u_{0\rho}^{(h)}}\bigg/\frac{\partial f}{\partial
u_{00}^{(h)}}\bigg|_{(u_{00},\ldots,u_{d0})=(\zeta_{0},\ldots,\zeta_{d})}=
\overline{\frac{\partial f}{\partial u_{0\rho}^{(h)}}}
\bigg/\overline{\frac{\partial f}{\partial u_{00}^{(h)}}},
\]
we have
\[\frac{\partial \xi_{\rho}}{\partial
u_{0i}^{(h)}}=\frac{\frac{\partial }{\partial
u_{0i}^{(h)}}\bigg(\overline{\frac{\partial f}{\partial
u_{0\rho}^{(h)}}}\bigg)\overline{\frac{\partial f}{\partial
u_{00}^{(h)}}} -\frac{\partial }{\partial
u_{0i}^{(h)}}\bigg(\overline{\frac{\partial f}{\partial
u_{00}^{(h)}}}\bigg)\overline{\frac{\partial f}{\partial
u_{0\rho}^{(h)}}}}{(\overline{\frac{\partial f}{\partial
u_{00}^{(h)}}})^2}=0,\] where $\overline{\frac{\partial f}{\partial
u_{0\rho}^{(h)}}}$ is obtained by replacing $(u_{00},\ldots,u_{d0})$
by $(\zeta_{0},\ldots,\zeta_{d})$ in $ \frac{\partial f}{\partial
u_{0\rho}^{(h)}}$.
\newline
Or equivalently
\[\frac{\partial }{\partial
u_{0i}^{(h)}}\bigg(\frac{\partial f}{\partial
u_{0\rho}^{(h)}}\bigg)\frac{\partial f}{\partial
u_{00}^{(h)}}-\frac{\partial }{\partial
u_{0i}^{(h)}}\bigg(\frac{\partial f}{\partial
u_{00}^{(h)}}\bigg)\frac{\partial f}{\partial u_{0\rho}^{(h)}}\in
\sat(f).\] Set $u_{00}^{(h)}=\gamma_{\tau}$, we have
\[ \frac{\partial }{\partial
u_{0i}^{(h)}}\bigg(f_{\tau \rho}\bigg)f_{\tau0}-\frac{\partial
}{\partial u_{0i}^{(h)}}\bigg(f_{\tau0}\bigg)f_{\tau \rho}=0.\]
Thus,
\begin{equation}\label{eq-fac12}
\frac{\partial \xi_{\tau\rho}}{\partial
u_{0i}^{(h)}}=\frac{\partial (f_{\tau\rho}/f_{\tau0})}{\partial
u_{0i}^{(h)}} =0.
\end{equation}

Secondly,  set $(u_{00},\ldots,u_{d0})=(\zeta_0,\ldots,\zeta_{d})$
in the equation
\[\sum_{\rho=1}^n u_{0\rho}^{(h)} \frac{\partial f}{\partial
u_{0\rho}^{(h)}}+u_{00}^{(h)} \frac{\partial f}{\partial
u_{00}^{(h)}}+\sum_{k=0}^{h-1}\sum_{\rho=0}^n u_{0\rho}^{(k)}
\frac{\partial f}{\partial u_{0\rho}^{(k)}}=rf,r \in \mathbb{N}.\]
We have
\[\sum_{\rho=1}^n
u_{0\rho}^{(h)} \overline{\frac{\partial f}{\partial
u_{0\rho}^{(h)}}}+\zeta_{0}^{(h)} \overline{\frac{\partial
f}{\partial u_{00}^{(h)}}}+\overline{\sum_{k=0}^{h-1}\sum_{\rho=0}^n
u_{0\rho}^{(k)} \frac{\partial f}{\partial u_{0\rho}^{(k)}}}=0.\] By
Theorem \ref{th-generic},
\[\sum_{\rho=1}^n
u_{0\rho}^{(h)} \xi_{\rho}+\zeta_{0}^{(h)}
+\overline{\sum_{k=0}^{h-1}\sum_{\rho=0}^n u_{0\rho}^{(k)}
\frac{\partial f}{\partial
u_{0\rho}^{(k)}}}\bigg/\overline{\frac{\partial f}{\partial
u_{00}^{(h)}}}=0.\] Then,
\[\frac{\partial}{\partial u_{0i}^{(h)}}\bigg(\overline{\sum_{k=0}^{h-1}\sum_{\rho=0}^n u_{0\rho}^{(k)}
\frac{\partial f}{\partial
u_{0\rho}^{(k)}}}\bigg/\overline{\frac{\partial f}{\partial
u_{00}^{(h)}}} \bigg)=-\xi_{i}-(-\xi_{i})=0 \]\[ =
\frac{\frac{\partial }{\partial
u_{0i}^{(h)}}\bigg(\overline{\sum_{k=0}^{h-1}\sum_{\rho=0}^n
u_{0\rho}^{(k)} \frac{\partial f}{\partial
u_{0\rho}^{(k)}}}\bigg)\overline{\frac{\partial f}{\partial
u_{00}^{(h)}}} -\frac{\partial }{\partial
u_{0i}^{(h)}}\bigg(\overline{\frac{\partial f}{\partial
u_{00}^{(h)}}} \bigg)\bigg(\overline{\sum_{k=0}^{h-1}\sum_{\rho=0}^n
u_{0\rho}^{(k)} \frac{\partial f}{\partial
u_{0\rho}^{(k)}}}\bigg)}{\left(\overline{\frac{\partial f}{\partial
u_{00}^{(h)}}} \right)^2}.
\]
Thus, we have $\frac{\partial }{\partial
u_{0i}^{(h)}}\bigg(\sum_{k=0}^{h-1}\sum_{\rho=0}^n u_{0\rho}^{(k)}
\frac{\partial f}{\partial u_{0\rho}^{(k)}}\bigg)\frac{\partial
f}{\partial u_{00}^{(h)}}-\frac{\partial }{\partial
u_{0i}^{(h)}}\bigg(\frac{\partial f}{\partial
u_{00}^{(h)}}\bigg)\bigg(\sum_{k=0}^{h-1}\sum_{\rho=0}^n
u_{0\rho}^{(k)}$ $ \frac{\partial f}{\partial u_{0\rho}^{(k)}}\bigg)
\in \sat(f)$.
From Lemma \ref{lm-ftaup}, by setting $u_{00}^{(h)}=\gamma_{\tau}$
in the above relation, we obtain
\begin{equation}\label{eq-fac13}
\frac{\partial
\big(\sum_{k=0}^{h-1}\sum_{\rho=0}^n u_{0\rho}^{(k)}
\widehat{\frac{\partial f}{\partial
u_{0\rho}^{(k)}}}\Big/f_{\tau0}\big)}{\partial u_{0i}^{(h)}}=0.
\end{equation}

From \bref{eq-fac12} and \bref{eq-fac13}, $\xi_{\tau \rho}$ and
$\sum_{k=0}^{h-1}\sum_{\rho=0}^n u_{0\rho}^{(k)}
\widehat{\frac{\partial f}{\partial u_{0\rho}^{(k)}}}\Big/f_{\tau0}$
are algebraically independent  of $u_{0i}^{(h)}\,(i=1,\ldots,$ $n)$.
Then the symmetric functions of $\xi_{\tau1},\ldots,\xi_{\tau
n},\sum_{k=0}^{h-1}$ $\sum_{\rho=0}^n u_{0\rho}^{(k)}
\widehat{\frac{\partial f}{\partial u_{0\rho}^{(k)}}}\Big/f_{\tau0}$
are rational functions in the set of indeterminates
$\{u_{ik},\ldots,u_{ik}^{(h)},$ $u_{0k},\ldots,u_{0k}^{(h-1)}:\,
i=1,\ldots,d;$ $\,k=0,\ldots,n\}$ only. Therefore,
$\prod^g_{\tau=1}(u_{00}^{(h)}-\gamma_{\tau})=\frac{\phi}{\psi}$
where $\psi$ is free of $u_{0i}^{(h)}\,(i=1,\ldots,n)$ and
$\gcd(\phi,\psi)=1$. Thus $A\phi=f\psi$. Since $f$ is irreducible,
we conclude that $A=\psi$ is free of $u_{0i}^{(h)}\,(i=1,\ldots,n)$.
\qed

Note that the factorization \bref{eq-fac10} is formal in the sense
that different factors are from different differential fields
$\ff_\tau$. The following result shows how to obtain a factorization
in the same extension field.
\begin{theorem}\label{lm-uni}
The  quantities $\xi_{\tau 1},\ldots,\xi_{\tau n}$ in
\bref{eq-fac10} are unique and \bref{eq-fac11} is a factorization of
$F$ as an algebraic polynomial in $u_{00}^{(h)},\ldots,u_{0n}^{(h)}$ in an extension
field of $\mathcal
{F}(u_{ik},\ldots,u_{ik}^{(h)},$ $u_{0k},\ldots,u_{0k}^{(h-1)}:i=1,\ldots,d;\,k=0,\ldots,n)$.

\end{theorem}
\proof From Lemma \ref{lm-free}, equations \bref{eq-fac12} and
\bref{eq-fac13}, we can see that $A(\bu_0,\ldots,\bu_d)$, $\xi_{\tau
j}$, and $\sum_{k=0}^{h-1}\sum_{\rho=0}^n u_{0\rho}^{(k)}
\widehat{\frac{\partial f}{\partial
u_{0\rho}^{(k)}}}\Big/f_{\tau0}$ are free of
$u_{0i}^{(h)}(i=1,\ldots,n)$. Then, \bref{eq-fac11} is a
factorization of the Chow form $F(\bu_0,\ldots,\bu_n)$ in the
polynomial ring $\mathcal
{F}(u_{ik},\ldots,u_{ik}^{(h)},$ $u_{0k},$ $\ldots,u_{0k}^{(h-1)}:i=1,\ldots,d;\,k=0,\ldots,n)[u_{00}^{(h)},\ldots,
u_{0n}^{(h)}]$. Thus, the factorization \bref{eq-fac11} must be
unique and hence $\xi_{\tau i}$.\qed

\subsection{Leading differential degree of an irreducible differential variety}
In this subsection, we will define the leading differential degree
for a prime $\delta$-ideal and give its geometric meaning.

\begin{definition}\label{def-ldeg}
Let  $F(\bu_{0},\bu_{1},\ldots,\bu_{d})=f(\bu;u_{00},\ldots,u_{d0})$
be the Chow form of a prime $\delta$-$\ff$-ideal $\CI$ of dimension
$d$ and order $h$.
By Lemma \ref{lm-cfl}, the number $g=\deg(f,u_{00}^{(h)})$ is an
invariant of $\CI$ under invertible linear transformations, which is
called the {\em leading differential degree} of $\CI$.
\end{definition}

From \bref{eq-fac10}, there exist $g$ points
$(\xi_{\tau1},\ldots,\xi_{\tau n})\,(\tau=1,\ldots,g)$, which have
interesting properties.

\begin{theorem} \label{th-zerochowform}
The points $(\xi_{\tau1},\ldots,\xi_{\tau n})\,(\tau=1,\ldots,g)$ in
\bref{eq-fac10} are  generic points of the $\delta$-$\ff$-variety
$V$. If $d>0$, they also satisfy the equations \[u_{\sigma
0}+\sum_{\rho=1}^n u_{\sigma \rho}y_{
\rho}=0\,(\sigma=1,\ldots,d).\]
\end{theorem}
\proof  Suppose $\phi(y_{1},\ldots,y_{n}) \in \mathcal {F}\{\Y\}$ is
any $\delta$-polynomial vanishing on $V$. Then
$\phi(\xi_{1},\ldots,\xi_{n})=0$. From Theorem~\ref{th-generic},
$\xi_{\rho}=\overline{\frac{\partial f}{\partial
 u_{0\rho}^{(h)}}}\big/\overline{\frac{\partial{f}}{\partial{u}_{00}^{(h)}}}$,
%
%
we have \[\phi\Big(\overline{\frac{\partial f}{\partial
u_{01}^{(h)}}}\bigg/\overline{\frac{\partial{f}}{\partial
u_{00}^{(h)}}}, \ldots,\overline{\frac{\partial f}{\partial
u_{0n}^{(h)}}}\bigg/\overline{\frac{\partial{f}}{\partial
u_{00}^{(h)}}}\Big)=0,\]
where $\overline{\frac{\partial f}{\partial u_{0\rho}^{(h)}}}$ are
obtained by substituting $\zeta_{i}$ to $u_{i0}\,(i=0, 1, \ldots,
d)$ in $\frac{\partial f}{\partial u_{ 0\rho}^{(h)}}$.

Hence, $ \phi(\frac{\partial f}{\partial
u_{01}^{(h)}}\Big/\frac{\partial f}{\partial
u_{00}^{(h)}},\ldots,\frac{\partial f}{\partial u_{0n}^{(h)}}\Big/
\frac{\partial f}{\partial u_{00}^{(h)}})$ vanishes for
$(u_{00},\ldots,u_{d0})=(\zeta_{0},\ldots,\zeta_{d})$. Then there
exists an $ m \in \mathbb{N}$, such that  $(\frac{\partial
f}{\partial u_{00}^{(h)}})^m \phi(\frac{\partial f}{\partial
u_{01}^{(h)}}\Big/\frac{\partial f}{\partial
u_{00}^{(h)}},\ldots,\frac{\partial f}{\partial
u_{0n}^{(h)}}\Big/\frac{\partial f}{\partial u_{00}^{(h)}}) \in
\sat(f)$ in $\ff\{\tilde{\bu},u_{00}\}$.
%
%
By Lemma \ref{lm-ftaup}, we have
$(f_{\tau0})^m \phi(\xi_{\tau1},\ldots,\xi_{\tau n})=0$, thus
$\phi(\xi_{\tau1},\ldots,\xi_{\tau n})=0$, which means that
$(\xi_{\tau1},\ldots,\xi_{\tau n})$ $ \in V$.

Conversely, for any $p \in \mathcal {F}\{\Y\}$ such that
$p(\xi_{\tau1},\ldots,\xi_{\tau n})=0$, there exists an $l \in
\mathbb{N}$ such that $\widetilde{p}=(\frac{\partial f}{\partial
u_{00}^{(h)}})^l p(\frac{\partial f}{\partial
u_{01}^{(h)}}\Big/\frac{\partial f}{\partial
u_{00}^{(h)}},\ldots,\frac{\partial f}{\partial
u_{0n}^{(h)}}\Big/\frac{\partial f}{\partial u_{00}^{(h)}})$ is in
$\ff\{\tilde{\bu},u_{00}\}$, which vanishes at
$u_{00}^{(h)}=\gamma_{\tau}$. By Lemma \ref{lm-ftaup},
$\widetilde{p} \in \sat(f)$.
Now treating $\widetilde{p}$ as a $\delta$-polynomial in
$\ff\langle\bu\rangle\{u_{00},\ldots,u_{d0}\}$, we have
$\widetilde{p}(\xi_{1},\ldots,$ $\xi_{n})=0$ and hence
 $p(\xi_{1},\ldots,$ $\xi_{n})=0$.
So $(\xi_{\tau1},\ldots,\xi_{\tau n})$ is a generic point of $V$.

Since $\overline{\frac{\partial f}{\partial
 u_{0\rho}^{(h)}}}+\overline{\frac{\partial{f}}{\partial{u}_{00}^{(h)}}}(-\xi_{\rho})=0$,
we have $\sum_{\rho=1}^n u_{\sigma \rho}\overline{\frac{\partial
f}{\partial
u_{0\rho}^{(h)}}}+\zeta_{\sigma}\overline{\frac{\partial{f}}{\partial{u}_{00}^{(h)}}}=0$.
Thus, $\sum_{\rho=0}^n u_{\sigma \rho}\frac{\partial f}{\partial
u_{0\rho}^{(h)}}$ vanishes at
$(u_{00},\ldots,u_{d0})=(\zeta_{0},\ldots,\zeta_{d})$. In the case
$\sigma \neq 0$, $\sum_{\rho=0}^n u_{\sigma \rho}$ $\frac{\partial
f}{\partial u_{0\rho}^{(h)}}=0$. Consequently,
$u_{\sigma0}+\sum_{\rho=1}^n u_{\sigma \rho}\xi_{\tau
\rho}=0(\sigma=1,\ldots,d)$. \qed


\begin{example}\label{ex-ff11} Continue from Example~\ref{ex-ff1}.
We have $\xi_{11}=-u_{00}/u_{01}$ under the condition
${u_{00}'=\gamma_1=\frac{u_{00}}{u_{01}}
u_{01}'-2\sqrt{-1}\sqrt{u_{00}u_{01}}}$. Then $\xi_{11}' =
2\sqrt{-\frac{u_{00}}{u_{01}}}$ and  $\xi_{11}$ is a zero of the
original $\delta$-ideal $\sat(y_1'^2-4y_1)$.
%
%
\end{example}

\vskip5pt Now, we will prove a result which gives the geometrical
meaning of the leading $\delta$-degree.

%
Suppose $F(\bu_0,\ldots,\bu_d)$ is the Chow form of $V$ which is of
dimension $d$, order $h$, and leading differential degree $g$.
Recall that by saying a point $\eta=(\eta_1,\ldots,\eta_n)$ lying on
$^a\P^{(k)}_0$ we mean that
$\eta^{[k]}=(\eta_1,\ldots,\eta_n,\ldots,\eta_1^{(k)},$
$\ldots,\eta_n^{(k)})$ is an algebraic zero of $^a\P^{(k)}_0$.
Theorem \ref{th-zerochowform} and Corollary \ref{co-f1} show that
$(\xi_{\tau1},\ldots,\xi_{\tau n})(\tau=1,\ldots,g)$ are
intersection points of $V$ and $\mathbb{P}_i=0(i=1,\ldots,d)$ as
well as $^a\P^{(k)}_0=0\,(k=0,\ldots,h-1)$.
In the next theorem, we will prove the converse of this result, that
is,
 $(\xi_{\tau1},\ldots,\xi_{\tau
n})\,(\tau=1,\ldots,g)$ are the only elements in $V$ which are also
on $\P_i=0\, (i=1,\ldots,d)$ as well as on
$^a\P^{(k)}_0=0\,(k=0,\ldots,h-1)$. Intuitively, we use
$\P_i=0\,(i=1,\ldots,d)$ to decrease the dimension of $V$ to zero
and use $^a\P^{(k)}_0=0\,(k=0,\ldots,h-1)$ to determine the $h$
arbitrary constants in the solutions of the zero dimensional
$\delta$-variety.

\begin{theorem}\label{th-g}
$(\xi_{\tau1},\ldots,\xi_{\tau n})\,(\tau=1,\ldots,g)$ defined in
\bref{eq-fac10} are the only elements of $V$ which also lie on
$\P_i(i=1,\ldots,d)$\footnote{If $d=0$, $\P_i(i=1,\ldots,d)$ is
empty.} as well as on $^a\P_0^{(j)}(j=0,\ldots,h-1)$.
\end{theorem}
\proof  Firstly, by Theorem \ref{th-zerochowform} and Corollary
\ref{co-f1}, $(\xi_{\tau1},\ldots,\xi_{\tau n})\,(\tau=1,\ldots,g)$
are solutions of $\mathbb{I}(V)$ and $\P_1,\ldots,\P_d$ which also
lie on $^a\P_0,^a\P'_0,\ldots,^a\P^{(h-1)}_0$. It suffices to show
that the number of solutions of $\mathbb{I}(V)$ and
$\P_1,\ldots,\P_d$ which also lie on
$^a\P_0,^a\P'_0,\ldots,^a\P^{(h-1)}_0$ does not exceed $g$.

Let
$\mathcal{J}=[\mathbb{I}(V),\P_1,\ldots,\P_d]\subset\mathcal{F}^\star\{\Y\}$,
where $\mathcal{F}^\star=\ff\langle\cup _{i=1}^d\bu_i\rangle$. By
Theorem \ref{th-inter2}, $\mathcal{J}$ is a prime $\delta$-ideal of
dimension zero and order $h$. Let $\mathcal{J}^{<h>} =
\mathcal{J}\cap \mathcal{F}^\star[\Y^{[h]}]$.
Since  $\mathcal{J}$ is of dimension zero and order $h$, its
$\delta$-dimension polynomial is of the form $\omega(t)=h,$ for $t
\geq h$. So $\mathcal{J}^{<h>}$ is an algebraic prime ideal of
dimension $h$.

Let
$\mathcal{J}_0=(\mathcal{J}^{<h>},^a\P_0,^a\P'_0,\ldots,^a\P^{(h-1)}_0
)$ be an algebraic ideal in the polynomial ring
$\ff_0[\Y^{[h]},u_{00},\ldots,u_{00}^{(h-1)}]$, where
$\ff_0=\mathcal {F}^\star( \cup_{j=1}^nu_{0j}^{[h-1]})$ and
$u_{0j}^{[h-1]}$ denote the set
$\{u_{0j},u_{0j}',\ldots,u_{0j}^{(h-1)}\}$. Similarly to the proof
of Lemma~\ref{lm-pu}, we can show that $\mathcal{J}_0$ is an
algebraic prime ideal of dimension $h$. If we can prove that
$u_{00}^{[h-1]}=\{u_{00},\ldots,u_{00}^{(h-1)}\}$ is a parametric
set of $\mathcal{J}_0$, then it is  clear that $\mathcal{J}_1=
(\mathcal{J}_0)$ is an algebraic prime ideal of dimension zero in
$\ff_0(u_{00}^{[h-1]})[\Y^{[h]}]=\ff^*(\bu_{0}^{[h-1]})[\Y^{[h]}]$.
So we need to prove that $u_{00}^{[h-1]}$ is a parametric set of
$\mathcal{J}_0$. Suppose the contrary, then there exists a nonzero
$\delta$-polynomial involving only
$\{u_{00},\ldots,u_{00}^{(h-1)}\}$ as well as the other $u$, which
belongs to $\mathcal{J}_0$. Such a $\delta$-polynomial also belongs
to $[\mathcal{J},\P_0]\in \mathcal {F}^\star\langle
u_{01},\ldots,u_{0n}\rangle\{y_1,\ldots,y_n,$ $u_{00}\}$. From the
proof of Theorem~\ref{th-choworder}, $ \{F,S_Fy_{1}-\frac{\partial
F}{\partial u_{01}^{(s)}},\, \ldots,\, S_Fy_{n}-\frac{\partial
F}{\partial u_{0n}^{(s)}}\}$ is a characteristic set of
$[\mathcal{J},\P_0]$ w.r.t. the elimination ranking $u_{00}\prec
y_{1} \prec \ldots \prec y_{n}$.  So this $\delta$-polynomial can be
reduced to zero by $F(\bu_0,\ldots,\bu_d)$. But $\ord(F,u_{00})=h$,
a contradiction. So we have  proved that $\mathcal {J}_1$ is an
algebraic prime ideal of dimension zero.

Clearly, $\mathcal {J}_2=(\mathcal {J}_1,^a\P_0^{(h)})\in \mathcal
{F}_0\big(u_{00}^{[h-1]},\cup_{j=1}^nu_{0j}^{[h]}\big)[\Y^{[h]},$
$u_{00}^{(h)}]$ is  an algebraic prime ideal of dimension zero.
Then, there exists an irreducible polynomial involving only
$u_{ij}^{(k)}$ and $u_{00}^{(h)}$. Similarly as above, it also
belongs to $[\mathcal{J},\P_0]$, thus it can be divisible by $F$.
Since $F$ is irreducible, it differs from $F$ only by a factor in
$\mathcal {F}$. Thus, $F=f(\bu;u_{00},u_{10},\ldots,u_{d0})\in
\mathcal{J}_2$.

Let $(\xi_1,\ldots,\xi_n)$ ba a generic point of $V$ and
$\zeta_{i}=-\sum_{j=1}^n u_{ij}\xi_j\,(i=0,\ldots,d)$. By Lemma
\ref{lm-cf1}, the $\delta$-ideal
$[\mathbb{I}(V),\P_1,\ldots,\P_d,\P_0]$ in $\mathcal
{F}\langle\bu\rangle\{y_1,\ldots,y_n,u_{00},\ldots,u_{d0}\}$ has a
generic point $(\xi_1,\ldots,\xi_n,\zeta_0,\ldots,\zeta_d)$. Since
$f(\bu;\zeta_0,\ldots,\zeta_d)=0$, differentiate both sides of this
identity w.r.t. $u_{0j}^{(k)}$, we have the following identities
 $$\overline{\frac{\partial f}{\partial u_{0j}^{(k)}}}+
  \sum_{l=k}^h \overline{\frac{\partial f}{\partial u_{00}^{(l)}}}
  \bigg(-{l \choose k}\xi_{j}^{(l-k)}\bigg)=0, \qquad (j=1,\ldots,n;k=0,\ldots,h),$$
where $\overline{\frac{\partial f}{\partial u_{0i}^{(j)}}}$ are
obtained by  substituting $\zeta_i$ to $u_{i0}$ in $\frac{\partial
f}{\partial u_{0i}^{(j)}}$. Let $g_{jk}={h\choose k}\frac{\partial
f}{\partial u_{00}^{(h)}}y_{j}^{(k)}+\sum_{l=1}^k {h-l\choose
k-l}\frac{\partial f}{\partial
u_{00}^{(h-l)}}y_{j}^{(k-l)}-\frac{\partial f}{\partial
u_{0j}^{(h-k)}}\,(j=1,\ldots,n;k=0,\ldots,h)$. Then
$g_{jk}\in[\mathbb{I}(V),\P_1,$ $\ldots,\P_d,\P_0]$ $\subset
[\mathcal {J},\P_0]$, for $g_{jk}$ vanishes at
$(\xi_1,\ldots,\xi_n,$ $\zeta_0,\ldots,\zeta_d)$. Denote the
algebraic ideal $[\mathcal {J},\P_0]\cap
\mathcal{F}^\star(\cup_{j=1}^nu_{0j}^{[h]})[\Y^{[h]},$ $u_{00},$
$\ldots,u_{00}^{(h)}]$ by $[\mathcal {J},\P_0]^{<h>}$. It is clear
that $g_{jk}\in [\mathcal {J},\P_0]^{<h>}$. We will show that
$[\mathcal {J},\P_0]^{<h>}=(\mathcal
{J}^{<h>},^a\P_0,\ldots,^a\P^{(h)}_0)$, which implies that
$g_{jk}\in \mathcal {J}_2$. Let $\eta=(\eta_1,\ldots, \eta_n)$ be a
generic point of $\mathcal{J}$. Then $(\eta_1,\ldots, \eta_n,$
$-\sum_{j=1}^n u_{0j}\eta_j)$ is a generic point of $[\mathcal
{J},\P_0].$ Thus,
 $$(\eta_1,\ldots, \eta_n,\ldots,\eta_1^{(h)}, \ldots, \eta_n^{(h)},-\sum_{j=1}^n
 u_{0j}\eta_j, \ldots,-\sum_{j=1}^n (u_{0j}\eta_j)^{(h)})$$
is a generic point of $[\mathcal {J},\P_0]^{<h>}$. Of course, it is
also a generic point of $(\mathcal {J}^{<h>},^a\P_0,$
$\ldots,^a\P^{(h)}_0)$. So the two ideals are identical. Thus,
$g_{jk}$ belongs to $\mathcal {J}_2$. Note that the coefficient of
$y_{j}^{(k)}$ in $g_{jk}$ is ${h \choose k}S_f={h \choose
k}\frac{\partial f}{\partial u_{00}^{(h)}}$. So $\mathbb{V}(\mathcal
{J}_2)\subseteq
\mathbb{V}(f(\U_0,u_{00}^{(h)}),g_{jk}:j=1,\ldots,n;k=0,\ldots,h)$
and the latter algebraic variety of consists exactly $g$ elements.
Thus, $|\mathbb{V}(\mathcal {J}_1)| =|\mathbb{V}(\mathcal
{J}_2)|\leq g.$ Since every solution of $\mathcal{J}$ which also
lies on $^a\P_0,^a\P'_0,\ldots,^a\P^{(h-1)}_0$, when truncated up to
order $h$,  becomes a solution of $\mathcal {J}_1$, it follows that
the number of solutions of $\mathbb{I}(V)$ and $\P_1,\ldots,\P_d$
which also lie on $^a\P_0,^a\P'_0,\ldots,^a\P^{(h-1)}_0$ does not
exceed $g$. \qed

With Theorems \ref{th-fac1}, \ref{lm-uni}, \ref{th-zerochowform},
and \ref{th-g}, we proved the second and third statements of Theorem
\ref{th-main2}.

From the proof, we can see that a zero dimensional algebraic ideal
is obtained as shown by the following corollary.
\begin{cor}\label{co-zdimi}
Let $\CI$ be a prime $\delta$-ideal in $\ff\{\Y\}$ of dimension $d$,
order $h$, and leading $\delta$-degree $g$. Use the same notations
as Theorem \ref{th-g}. Then
$$\widetilde{\mathcal{I}}=\big([\CI,\P_1,\ldots,\P_d]\cap\ff^\star[\Y^{[h]}],^a\P_0,\ldots,^a\P_0^{(h-1)}\big)\subset
\ff^\star(\bu_0^{[h-1]})[\Y^{[h]}]$$
is an algebraic prime ideal of dimension zero whose solutions are
exactly $(\xi_{\tau1},\ldots,$ $\xi_{\tau n})^{[h]}\,$
$(\tau=1,\ldots,g)$, where $\ff^\star=\ff\langle\cup
_{i=1}^d\bu_i\rangle$.
\end{cor}

\begin{example} \label{ex-t1}
Continue from Example~\ref{ex-ff1}. Note that $d=0$, $h=1$, and
$g=2$.
Let $V$ be the general component of $p=y_1'^2 - 4 y_1=0$. As in
Theorem \ref{th-g}, we introduce the equation $^a\P_0 = u_{00} +
u_{01} y_1$ which intersects $V$ at two points: $\xi_{11} =
-u_{00}/u_{01}$ with $\delta_1 \xi_{11} = 2\sqrt{-u_{00}/u_{01}}$
and $\xi_{21} = -u_{00}/u_{01}$ with $\delta_2 \xi_{21} =
-2\sqrt{-u_{00}/u_{01}}$.
As indicated by Corollary \ref{co-zdimi},
$(\xi_{11},\delta_1\xi_{11})$ and $(\xi_{21},\delta_2\xi_{21})$ are
the only solutions of the algebraic ideal
$\widetilde{\mathcal{I}}=(p,^a\P_0)\subset\Q(u_{00},u_{01})[y_1,y_1']$.
%
%
\end{example}

%
\vskip 10pt

Due to Corollary \ref{co-zdimi}, we can give a differential analog
to the Stickelberger's Theorem in algebraic geometry.

\begin{theorem}[{\bf Stickelberger's Theorem}] \cite[p. 54]{cox2}, \cite{Rouillier}\label{th-astickelberger}
Let $\mathcal{P}\subset \ff[\Y]$ be a zero-dimensional ideal. Denote
$A=\ff[\Y]/\mathcal {P}$. Then $A$ is a finite dimensional vector
space over $\ff.$ For any polynomial $f\in\ff[\Y]$, let $L_f$ be the
$\ff$-linear map:
\[\begin{array}{cccc}L_f:\,& A &\longrightarrow &A\\
                     \quad&          \overline{g}&\quad&  \overline{f g}
  \end{array}
\] where $\overline{g}$ denotes the residue class of $g\in\ff[\Y]$ in $A=\ff[\Y]/\mathcal {P}$.
Then the eigenvalues of $L_f$ are $f(\alpha)$, with multiplicity
$m_\alpha$ where $\alpha\in \mathbb{V}(\mathcal {P})$ and $m_\alpha$
is the multiplicity of $\alpha\in \mathbb{V}(\mathcal {P})$. Thus,
the  determinant of $L_f$ is $\prod_{\alpha\in \mathbb{V}(\mathcal
{P})}f(\alpha)^{m_\alpha}$.
\end{theorem}


\begin{theorem}[{\bf Differential Stickelberger's Theorem}]\label{th-dstickelberger}
Let $\CI$ be a prime $\delta$-ideal in $\ff\{\Y\}$ of dimension
zero, order $h$, and leading $\delta$-degree $g$. Let
$\P_0=u_{00}+u_{01}y_1+\cdots+u_{0n}y_n=0$ be a generic
$\delta$-hyperplane and $^a\P_0,\ldots,^a\P_0^{(h-1)}$ defined in
\bref{eq-apoly}. For any $\delta$-polynomial $p\in \ff\{\Y\}$, let
$s = \max\{h,\ord(p)\}$ and $L_p$ the $\ff(\bu_0^{[h-1]})$-linear
map:
\[\begin{array}{cccc}L_p:\,& \ff(\bu_0^{[h-1]})[\Y^{[s]}]/\widetilde{\mathcal{I}} &\longrightarrow &\ff(\bu_0^{[h-1]})[\Y^{[s]}]/\widetilde{\mathcal{I}}\\
                     \quad&          \overline{g}&\quad&  \overline{gp}
  \end{array}
\]
where
$\widetilde{\mathcal{I}}=(\CI\cap\ff[\Y^{[s]}],^a\P_0,\ldots,^a\P_0^{(h-1)})\subset
\ff(\bu_0^{[h-1]})[\Y^{[s]}]$ and $\overline{g}$ denotes the residue
class of $g\in\ff(\bu_0^{[h-1]})[\Y^{[s]}]$ in
$\ff(\bu_0^{[h-1]})[\Y^{[s]}]/\widetilde{\mathcal{I}}$. Then the
eigenvalues of $L_p$ are $p(\xi_\tau)$ and the  determinant of $L_p$
is $\prod_{\tau=1}^g p(\xi_\tau)$.
\end{theorem}
\proof By Theorem \ref{th-astickelberger}, it suffices to show that
$\widetilde{\mathcal{I}}$ is a prime  $\delta$-ideal in
$\ff(\bu_0^{[h-1]})[\Y^{[s]}]$ of dimension zero, and
$\mathbb{V}(\widetilde{\mathcal{I}})=\{\xi_{\tau}^{[s]}=(\xi_{\tau1},\ldots,\xi_{\tau
n},\ldots,\delta_\tau^{s}\xi_{\tau1},\ldots,\delta_\tau^{s}\xi_{\tau
n}):\tau=1,\ldots,g\}$. If $s=h$, then this is a direct consequence
of Corollary~\ref{co-zdimi}.
If $s>h$, $\xi_{\tau}^{[s]}$ clearly vanishes
$^a\P_0^{(i)}(i=0,\ldots,h-1)$. By Theorem \ref{th-g},
$\xi_{\tau}^{[s]}$ are also zeros of $\CI\cap\ff[\Y^{[s]}]$. By
Theorem~\ref{th-g}, $\xi_{\tau}^{[s]}$ are the only zeros of
$\widetilde{\mathcal{I}}$.\qed

\subsection{Relations between the differential Chow form and the variety}
\label{sec-var} In the algebraic case,  we can obtain the defining
equations of a variety from its Chow form. But in the differential
case, this is not valid. Now we proceed as follows to obtain a
weaker result. Recall that a $\delta$-variety is unmixed if all of
its components have the same dimension.

\begin{lemma} \label{lemma-genericthrough0}
Let $ V$ be an irreducible $\delta$-$\ff$-variety of dimension $d>0$
and $(0,\ldots,0)$ $\notin V$. Then, the intersection of $V$ with
a generic $\delta$-hyperplane passing through $(0,\ldots,0)$
 is either empty  or unmixed of dimension  $d-1$. Moreover, in the case
 $d>1$, it is exactly unmixed of dimension $d-1$.
\end{lemma}
\proof Let $\mathcal {I}=\mathbb{I}(V )$ be the prime
$\delta$-$\ff$-ideal corresponding to $V$. A generic
$\delta$-hyperplane passing through $(0,\ldots,0)$ is
$u_{1}y_{1}+u_{2}y_{2}+\cdots+u_{n}y_{n}$ where the $u_{i}\in\ee$
are $\delta$-$\ff$-indeterminates. Since $(0,0,\ldots,0) \notin V$,
we have
\begin{eqnarray}
&\quad& V\cap \dZero(u_{1}y_{1}+u_{2}y_{2}+\cdots+u_{n}y_{n})
\nonumber \\ &=& \dZero(\mathcal
{I},u_{1}y_{1}+u_{2}y_{2}+\cdots+u_{n}y_{n}) \nonumber \\&=&
\bigcup^n_{i=1}\dZero([\mathcal
{I},u_{1}y_{1}+u_{2}y_{2}+\cdots+u_{n}y_{n}]/y_{i})\nonumber
\\ &=&\bigcup^n_{i=1}\dZero([\mathcal
{I},u_{1}y_{1}+u_{2}y_{2}+\cdots+u_{n}y_{n}]:y_{i}^\infty) \nonumber
\end{eqnarray}
 Suppose a generic point of $V$ is
$(\xi_{1},\ldots,\xi_{n})$. Since $(0,0,\ldots,0) \notin V$, there
exists at least one $i \in \{1,\ldots,n\}$, such that
$\xi_{i}\neq0$. Of course, $\xi_{i}=0$ means $\dZero([\mathcal
{I},u_{1}y_{1}+u_{2}y_{2}+\cdots+u_{n}y_{n}]:y_{i}^\infty)=\emptyset
$. So we need only to consider the case when $\xi_{i}\neq 0$.
Without loss of generality, we suppose $\xi_{1}\neq0$.

Let \[\mathcal {Q}=[\mathcal
{I},u_{1}y_{1}+u_{2}y_{2}+\cdots+u_{n}y_{n}]:y_{1}^\infty \subseteq
\mathcal {F}\langle
u_{1},\ldots,u_{n}\rangle\{y_{1},\ldots,y_{n}\}\] and \[\mathcal
{Q}_{0}=[\mathcal
{I},u_{1}y_{1}+u_{2}y_{2}+\cdots+u_{n}y_{n}]:y_{1}^\infty \subseteq
\mathcal {F}\langle u_{2},\ldots,u_{n}\rangle\{y_{1},\ldots,y_{n}
,u_{1}\}\ \] It is easy to verify that
$(\xi_{1},\ldots,\xi_{n},-\frac{u_{2}\xi_{2}+\cdots+u_{n}\xi_{n}}{\xi_{1}})$
is generic point of $\mathcal {Q}_{0}$
and $\dim(\mathcal {Q}_{0})=d$ by following the proof of
Lemma~\ref{lm-pu}. Now we discuss in three cases.

Case 1) $\mathcal {I}\cap \mathcal {F}\{y_{1}\} \neq \{0\}$, that
is, $\xi_{1}$ is $\delta$-algebraic over $\mathcal {F}$. We have
\[\dim V=d=\dtrdeg \,\mathcal {F}\langle \xi_{1},\ldots,\xi_{n} \rangle \big/ \mathcal {F}=\dtrdeg \,\mathcal {F}\langle \xi_{1}\rangle\langle \xi_{2},\ldots,\xi_{n} \rangle \big/ \mathcal {F}\langle\xi_{1}\rangle.\]
Suppose $\xi_{2},\ldots,\xi_{d+1}$ are $\delta$-independent over
$\mathcal {F}\langle \xi_{1}\rangle$.

Firstly, $\mathcal {Q}_{0}\cap \mathcal {F}\langle
u_{2},\ldots,u_{n}\rangle\{u_{1}\}=\{0\}$. For if not, we have a
nonzero $\delta$-polynomial $h(u_{2},\ldots,u_{n},u_{1})\in \mathcal
{F}\{ u_{2},\ldots,u_{n},u_{1}\}$ such that $h(u_{2},\ldots,u_{n},
-\frac{u_{2}\xi_{2}+\cdots+u_{n}\xi_{n}}{\xi_{1}})=0$. For a fixed
$i$ between 2 and $n$, if we specialize $u_{i}$ to $-1$, $u_{j}\,(j
\neq i)$ to 0, then by Theorem \ref{th-specil}, $\xi_{i}/\xi_{1}$ is
$\delta$-algebraic over $\mathcal {F}$. So each
$\xi_{i}\,(i=1,\ldots,n)$ is $\delta$-algebraic over $\mathcal {F}$,
which contradicts to the fact $d>0$. It follows that $\mathcal {Q}$
is not the unit ideal and $\dim(\mathcal {Q})\geq 0$.

Secondly, since $y_{2},\ldots,y_{d+1}$ is a parametric set of
$\mathcal {I}$, it is also a parametric set for $\mathcal {Q}_{0}$.
So $y_{2},\ldots,y_{d+1},u_{1}$ are $\delta$-dependent modulo
$\mathcal {Q}_{0}$. Since $\mathcal {Q}_{0}\cap \mathcal {F}\langle
u_{2},\ldots,u_{n}\rangle\{u_{1}\}=\{0\}$, we know that
$y_{2},\ldots,y_{d+1}$ are $\delta$-dependent modulo $\mathcal {Q}$.
Using the fact that each remaining $y_{i}$ and
$y_{2},\ldots,y_{d+1}$ are $\delta$-dependent modulo $Q$, we obtain
$\dim (\mathcal {Q})\leq d-1$. If $d=1$, then $\dim(\mathcal
{Q})=0=d-1$ follows. Now for $d>1$, we claim that $\dim (\mathcal
{Q})= d-1$ by proving that $y_{2},\ldots,y_{d}$ are
$\delta$-independent modulo $\mathcal {Q}$. For if not, there exists
$ 0 \neq h(y_{2},\ldots,y_{d},u_{1})\in \mathcal {Q}_{0}$ such that
$h(\xi_{2},\ldots,\xi_{d},-\frac{u_{2}\xi_{2}+\cdots+u_{n}\xi_{n}}{\xi_{1}})=0$.
By Theorem \ref{th-specil}, we can specialize $u_{d+1}$ to $-1$, the
other $u_{i}$ to zero, and conclude that
$\xi_{2},\ldots,\xi_{d},\frac{\xi_{d+1}}{\xi_{1}}$ are
$\delta$-dependent over $\mathcal {F}$. Since
$\xi_{2},\ldots,\xi_{d}$ are $\delta$-independent over $\mathcal
{F}$, $\xi_{d+1}$ is $\delta$-algebraic over $\mathcal {F}\langle
\xi_{1},\ldots,\xi_d\rangle$, which is a contradiction. Thus $\dim
\mathcal {Q}=d-1 $.

 Case 2) $d>1$ and $\xi_{1}$ is $\delta$-transcendental over
$\mathcal {F}$. In this case, we suppose a $\delta$-transcendence
basis is $\xi_{1},\ldots,\xi_{d}$.

Firstly, $\mathcal {Q}_{0}\cap \mathcal {F} \langle
u_{2},\ldots,u_{n}\rangle\{u_{1}\}=\{0\}$. For if not, as proceeded
in the preceding case, we conclude that $\xi_{i}/\xi_{1}$ is
$\delta$-algebraic over $\mathcal {F}$, that is, $\xi_{i},\xi_{1}$
are $\delta$-algebraic over $\mathcal {F}$, which contradicts to the
fact $d>1$. So $\mathcal {Q}$ is a nontrivial prime $\delta$-ideal.

Secondly, $\dim(\mathcal {Q})=d-1$, for on the one hand  from the
fact that $y_{1},y_{2},\ldots,y_{d},u_{1}$ are $\delta$-dependent
modulo $\mathcal {Q}_{0}$, we have $\dim(\mathcal {Q})\leq d-1$, and
on the other hand, from the fact that $y_{2},\ldots,y_{d},u_{1}$ are
$\delta$-independent modulo $\mathcal {Q}_{0}$, it comes
$\dim(\mathcal {Q})\geq d-1$.

Case 3) $d=1$ and $\xi_{1}$ is $\delta$-transcendental over
$\mathcal {F}$.
If $\mathcal {Q}_{0}\cap \mathcal {F}\langle
u_{2},\ldots,u_{n}\rangle$ $\{u_{1}\}\ne\{0\}$,   the intersection is empty.
If $\mathcal {Q}\neq [1]$, similar to case 2, we can easily
prove that the intersection is of dimension zero.

So for each $i \in \{1,\ldots,n\}$ such that $\xi_{i}\neq 0$, we can
show that $\dZero([\mathcal
{I},u_{1}y_{1}+u_{2}y_{2}+\cdots+u_{n}y_{n}]:y_{i}^\infty) $ is
either empty or of dimension $d-1$ similarly as the above steps for
the case $i=1$. And if $d >1$, it is exactly of dimension $d-1$.
Thus the theorem is proved. \qed

%

By saying independent generic $\delta$-hyperplanes, we mean that the
coefficients of these $\delta$-hyperplanes are
$\delta$-indedeterminates in $\ee$. The following result gives an
equivalent condition  for a point to be in a $\delta$-variety.
\begin{theorem} \label{th-x'generic}
Let $V$ be a $\delta$-$\ff$-variety of dimension $d$. Then
$\overline{x} \in V$  if and only if $d+1$ independent generic
$\delta$-hyperplanes $\P_{0},\P_{1},\ldots,\P_{d}$ passing through
$\overline{x}$ meet $V$.
\end{theorem}
\proof The necessity of the condition is obviously true. We now
consider the sufficiency. We adjoin the coordinates of
$\overline{x}$ to $\mathcal {F}$, and denote
$\mathcal{\overline{F}}$ to be the $\delta$-field thus obtained.
Regarded as a $\delta$-variety over $\mathcal {\overline{F}}$, $V$
is the sum of a finite number of irreducible $\delta$-varieties
$\overline{V}_{i}$, which are of dimension $d$  \cite[p.51]{ritt}.
Suppose $\overline{x} \notin V$, and therefore does not lie in any
component of $V$. We now prove that any $d+1$ independent generic
$\delta$-hyperplanes passing through $\overline{x}$ do not meet
$\overline{V}_{i}$. Without loss of generality, suppose
$\overline{x}=(0,0,\ldots,0)$. Then a generic $\delta$-hyperplane
passing through $\overline{x}$ is $s_{1}y_{1}+\cdots+s_{n}y_{n}$
where $s_{i}\in\ee$ are $\delta$-$\ff$-indeterminates. We proceed by
induction on $d$.

If $d=0$, then for $(a_{1},\ldots,a_{n})\in V$, each $a_{i}$ is
$\delta$-algebraic over $\mathcal {F}$. If $V\cap
\dZero(s_{1}y_{1}+\cdots+s_{n}y_{n}) \neq \emptyset $, then there
exists some $(a_{1},\ldots,a_{n}) \in V$ such that
$s_{1}a_{1}+\cdots+s_{n}a_{n}=0$.  Since the $s_{i}$ are
$\delta$-independent over $\mathcal {F}$. Thus
$(a_{1},\ldots,a_{n})=(0,\ldots,0)$, a contradiction to the fact
that $\overline{x} \notin V$. Thus the theorem is proved when $d=0$.

We therefore assume the truth of the theorem for $\delta$-varieties
of dimension less than $d$, and consider a $\delta$-variety $V$ of
dimension $d$. Let $\P_{0},\ldots,\P_{d}$ be $d+1$ independent
generic $\delta$-hyperplanes passing through $\overline{x}$. The
equation $\P_{d}$ can be written as $s_{1}y_{1}+\cdots+s_{n}y_{n}=0$
with $s_{i}$ $\delta$-$\ff$-indeterminates. From
Lemma~\ref{lemma-genericthrough0}, $\mathbb{P}_d=0$ meets $V$ in a
$\delta$-variety $\mathbb{W}$ of dimension less than $d$. By the
hypothesis of the induction, $\P_{0},\ldots,\P_{d-1}$ do not meet
$\mathbb{W}$, it follows that $V$ does not meet
$\P_{0},\ldots,\P_{d}$. Therefore the theorem is proved.\qed

The following result proves the fourth statement of Theorem
\ref{th-main2}.
\begin{theorem} \label{th-sf}
Let $F(\bu_{0},\bu_{1},\ldots,\bu_{d})$ be the Chow form of $V$ and
$S_{F}=\frac{\partial F}{\partial u_{00}^{(h)}}$.
Suppose that $\bu_i\,(i=0,\ldots,d)$ are $\delta$-specialized over
$\ff$ to sets $\bv_i$ of specific elements in  $\mathcal{E}$ and
$\overline{\P}_{i}\,(i=0,\ldots,d)$ are obtained by substituting
$\bu_i$ by $\bv_i$ in $\P_i$.
If\, $\overline{\P}_i=0(i=0,\ldots,d)$ meet $V$, then
$F(\bv_{0},\ldots,\bv_{d})$ $=0$. Furthermore, if
$F(\bv_{0},\ldots,\bv_{d})=0$ and $S_{F}(\bv_{0},\ldots,\bv_{d})\neq
0$, then the $d+1$ $\delta$-hyperplanes $\overline{\P}_{i}=0$
$(i=0,\ldots,d)$ meet $V$.
\end{theorem}
\proof Let $\mathcal{I}=\mathbb{I}(V)\subseteq \mathcal
{F}\{\Y\} $,  $\I_{\zeta,\xi}=[\mathcal
{I},\P_{0},\ldots,\P_{d}]\subseteq \mathcal
{F}\langle\bu\rangle\{y_{1},\ldots,y_{n},$ $u_{00},\ldots,u_{d0}\}$,
and
  $\mathcal{I}_1=[\I_{\zeta,\xi}]\subseteq \mathcal {F}\langle\bu_0,\ldots,\bu_{d}\rangle\{y_{1},\ldots,y_{n}\}$.
By Lemma \ref{lm-cf2},
$\{F,\frac{\partial F}{\partial u_{00}^{(h)}}y_{1}-\frac{\partial
F}{\partial u_{01}^{(h)}},\ldots,$ $\frac{\partial F}{\partial
u_{00}^{(h)}}y_{n}-\frac{\partial F}{\partial u_{0n}^{(h)}}\}$ is a
characteristic set of $\I_{\zeta,\xi}$ w.r.t. the elimination ranking
$u_{d0}\prec \cdots \prec u_{00}\prec y_{1}\prec\cdots\prec y_{n}$.
Since $F$ is irreducible, $\I_{\zeta,\xi}=[F,S_{F}y_{1}-\frac{\partial
F}{\partial u_{01}^{(h)}},\ldots,S_{F}y_{n}-\frac{\partial
F}{\partial u_{0n}^{(h)}}]:S_{F}^\infty$ with $S_{F}=\frac{\partial
F}{\partial u_{00}^{(h)}}$.

When $\bu_i$ are  $\delta$-specialized to $\bv_i$, $\mathcal {I}_1$
becomes a $\delta$-ideal in
$\mathcal{F}\langle\bv_0,\ldots,\bv_d\rangle\{\Y\}$. If
$\overline{\P}_{0},\ldots,\overline{\P}_{d}$ meet $V$, then
$\overline{\mathcal {I}_1}=[\mathcal
{I},\overline{\P}_{0},\ldots,\overline{\P}_{d}]\neq [1]$ which
implies $F(\bv_{0},\ldots,\bv_{d})=0$ since $F\in \I_{\zeta,\xi}$.

If $S_{F}(\bv_{0},\ldots,\bv_{d})\neq 0$ and
$F(\bv_{0},\ldots,\bv_{d})=0$, then let
$\overline{y}_{i}=(\frac{\partial F}{\partial
u_{0i}^{(h)}}(\bv_{0},\ldots,$ $\bv_{d}))\big/(S_{F}(\bv_{0},$
$\ldots,\bv_{d}))(i=1,\ldots,n)$. We claim that
$(\overline{y}_{1},\ldots,\overline{y}_{n})$ lies in $V$ and the
$d+1$ $\delta$-hyperplanes
$\overline{\P}_{0},\ldots,\overline{\P}_{d}$, which implies that
$\overline{\P}_{0},\ldots,\overline{\P}_{d}$ meet $V$.

Firstly, let $p$ be any $\delta$-polynomial in  $\mathcal {I}$. Then
$p \in \I_{\zeta,\xi}$, so there exists an integer $m$ such that
$S_{F}^m p\in [F,S_{F}y_{1}-\frac{\partial F}{\partial
u_{01}^{(h)}},\ldots,S_{F}y_{n}-\frac{\partial F}{\partial
u_{0n}^{(h)}}]$. If we specialize $u_{ij} \rightarrow v_{ij},u_{i0}
\rightarrow v_{i0}$ and let $y_{i}=\overline{y}_{i}$, then we have
$S_{F}^m
(\bv_{0},\ldots,\bv_{d})p(\overline{y}_{1},\ldots,\overline{y}_{n})=0$,
so $p(\overline{y}_{1},\ldots,\overline{y}_{n})=0$. That is,
$(\overline{y}_{1},\ldots,\overline{y}_{n}) \in V$. Secondly, since
$\P_{i} \in \I_{\zeta,\xi}$, similarly as the above, it follows that
$(\overline{y}_{1},\ldots,\overline{y}_{n}) $ lies in
$\overline{\P}_{i}$. So $\overline{\P}_{0},\ldots,\overline{\P}_{d}$
meet $V$.   \qed


\begin{remark} \label{re-sf}
Let $X$ be the set of all $(n-d-1)$-dimensional linear spaces in
$\ff^n$ that meet an irreducible $\delta$-$\ff$-variety $V$ of
dimension $d$. From Theorem \ref{th-sf},  $X\subset\V(F)$ and
$X\setminus \V(S_F) = \V(F)\setminus \V(S_F)$. That is, a ``major"
part of $X$ is known to be $\V(F)\setminus \V(S_F)$. An interesting
problem is to see whether $X$ is a $\delta$-variety for a projective
$\delta$-variety $V$. In \cite{kol51}, Kolchin showed that this
problem has a positive answer in a special case, that is, $V$ is a
projective algebraic variety which is treated as a zero dimensional
differential projective variety in certain way.
\end{remark}

\vskip5pt Similar to the algebraic case  \cite[p.22]{hodge}, we can
show that a generic $\delta$-hyperplane passing through a given
point $x=(x_{1},x_{2},\ldots,x_{n})$ is of the form
$a_{0}+a_{1}y_{1}+\cdots+a_{n}y_{n}=0$ with $a_{i}=\sum_{j=0}^n
s_{ij}x_{j}(i=0,1,\ldots,n)$, where $x_{0}=1$ and $S=(s_{ij})$  is
an $(n+1)\times(n+1)$ skew-symmetric matrix with $s_{ij}(i<j)$
independent $\delta$-$\ff$-indeterminates in $\ee$. That is,
\[\left( \begin{array}{c} a_{0}\\a_{1}\\ \vdots \\ a_{n}\end{array}\right)=S\left( \begin{array}{c} 1\\x_{1}\\ \vdots \\ x_{n}\end{array}\right)\]
For convenience, we denote such a $\delta$-hyperplane by $Sx$ and
say a generic $\delta$-hyperplane passing through a point $x$ is of
the form $Sx$.

Now we write $\bu_{i}=(u_{i0},u_{i1},\ldots,u_{in})^T=S^i Y$ where
$Y=(1,y_{1},\ldots,y_{n})^T$ and the $S^i$ are skew-symmetric
matrices with $s^i_{jk}(j<k)$ independent
$\delta$-$\ff$-indeterminates in $\ee$. Substituting the $\bu_{i}$
in $F(\bu_{0},\bu_{1},\ldots,\bu_{d})$ by these equations, we obtain
a $\delta$-polynomial involving $s^i_{jk}(j<k)$ and the $y_l$.
Regarding this $\delta$-polynomial as a $\delta$-polynomial in
$s^i_{jk}(j<k)$, then we have
$F(\bu_{0},\bu_{1},\ldots,\bu_{d})=F(S^0Y,S^1Y,\ldots,S^dY)=\sum
g_{\phi}(y_{1},\ldots,y_{n}) \phi(s^i_{jk})$ where $\phi(s^i_{jk})$
are different $\delta$-monomials. In this way, we get a finite
number of $\delta$-polynomials $g_{\phi}(y_{1},\ldots,y_{n})$ over
$\mathcal {F}$, which is denoted  by $\mathcal {P}$. Similarly, in
this way, we will get another set $\mathcal {D}$ of
$\delta$-polynomials  from $S_{F}(\bu_{0},\ldots,\bu_{d})$.

\begin{theorem}\label{th-var1}
Let $V$ be an irreducible $\delta$-$\ff$-variety with dimension $d$
and $F(\bu_{0},\bu_{1},$ $\ldots,\bu_{d})$ its Chow form. Then
$V\setminus\dZero(\mathcal {D})=\dZero(\mathcal {P})\setminus
\dZero(\mathcal {D})\ne\emptyset,$ where $\mathcal {P}, \mathcal
{D}$ are the $\delta$-polynomial sets obtained from
$F(\bu_{0},\bu_{1},\ldots,\bu_{d})$ and
$S_{F}(\bu_{0},\bu_{1},\ldots,$ $\bu_{d})$ respectively as above.
\end{theorem}
\proof On the one hand, for any $\overline{x} \in V$,  from
Theorem~\ref{th-x'generic}, any $d+1$ generic $\delta$-hyperplanes
passing through $\overline{x}$ meet $V$. So $S^0 \overline{x},S^1
\overline{x},\ldots,S^d \overline{x}$ meet $V$. By the proof of
Theorem~\ref{th-sf}, $F(S^0 \overline{x},S^1 \overline{x},\ldots,$
$S^d \overline{x})=0$. Since $s^i_{jk}(j<k)$ are
$\delta$-ndeterminates, $\overline{x}\in \dZero(\mathcal {P})$. So
$V\setminus\dZero(\mathcal {D}) \subseteq \dZero(\mathcal
{P})\setminus \dZero(\mathcal {D})$.

On the other hand, for any $\overline{x} \in \dZero(\mathcal
{P})\setminus \dZero(\mathcal {D})$, since any $d+1$ generic
$\delta$-hyperplanes passing through $\overline{x}$ are of the form
$S^0 \overline{x},S^1 \overline{x},\ldots,S^d \overline{x}$ with the
$S^i$ $\delta$-indeterminate skew-symmetric  matrices, we have
$F(S^0 \overline{x},S^1 \overline{x},\ldots,S^d \overline{x})=0$ and
$S_{F}(S^0 \overline{x},S^1 \overline{x},\ldots,S^d \overline{x})$
$\neq0$. From Theorem~\ref{th-sf}, $S^0 \overline{x},S^1
\overline{x},\ldots,S^d \overline{x}$ meet $V$. Thus from
Theorem~\ref{th-x'generic}, $\overline{x} \in V$. Thus
$V\setminus\dZero(\mathcal {D}) = \dZero(\mathcal {P})\setminus
\dZero(\mathcal {D})$.

Now, we show that $V\setminus\dZero(\mathcal {D})\neq \emptyset$.
Suppose the contrary, i.e. $V\subset \dZero(\mathcal {D})$, in
particular, its generic point $(\xi_{1},\ldots,\xi_{n}) \in
\dZero(\mathcal {D})$. Thus, $S_{F}(S^0 \xi,S^1 \xi,\ldots,S^d
\xi)=0$, where $\xi=(1,\xi_{1},\ldots, \xi_{n})$. Recall that
$s^i_{jk}(j<k;\,i=0,1,\ldots,d)$ are independent
$\delta$-indeterminates over $\mathcal
{F}\langle\xi_{1},\ldots,\xi_{n}\rangle$. Now we consider a
$\delta$-endomorphism $\phi$ of $\mathcal
{F}\langle\xi_{1},\ldots,\xi_{n}\rangle\{s^i_{jk}$
$(j<k;\,i=0,1,\ldots,d)\}$ over $\mathcal
{F}\langle\xi_{1},\ldots,\xi_{n}\rangle$ satisfying
$\phi(s^i_{0k})=-s^i_{0k}$ and
$\phi(s^i_{jk})=0\,(j<k;\,j=1,\ldots,n)$. It is clear that
$\phi(S_{F}(S^0 \xi,\ldots,S^d \xi))=S_{F}(s^i_{0k};-\sum_{k=1}^n
s^0_{0k}\xi_{k},\ldots,-\sum_{k=1}^n s^{d}_{0k}\xi_{k})=0$. Denoting
$s^i_{0k}$ by $u_{ik}$, we have $S_{F}(\bu;$
$\zeta_{0},\ldots,\zeta_{d})=0$, thus $S_{F} \in \sat(F)$, which is
a contradiction. So ${V}\setminus\dZero(\mathcal {D})\neq
\emptyset$.\qed

Since $V$ is an irreducible $\delta$-variety, $V\cap\dZero(\mathcal
{D})$ is a subset of $V$ with lower dimension than that of $V$ or
with the same dimension but of lower order. Thus,
$V\setminus\dZero(\mathcal {D})$ is an open set  of $V$ in the
Kolchin topology.

\begin{example}
 Continue from Example~\ref{example1}.
 In this example, $F(\bu_0)=u_{1}^2
(u_{0}')^2-2u_{1}u_{1}'u_{0}u_{0}'+(u_{1}')^2 u_{0}^2+4u_{1}^3u_{0}$
 and $S_{F}(\bu_0)=2u_{1}^2u_{0}'-2u_{1}u_{1}'u_{0}$. Following the
 steps as above, we obtain $\mathcal{P}=\{(y_1')^2-4y_1\}$ and $\mathcal{D}=\{y_1'\}$.
That is, we obtain the defining equation $(y_1')^2-4y_1=0$ for the
$\delta$-variety under the condition $y_1'\ne0$.
\end{example}

\section{Differential Chow variety}
\label{sec-dcv}

In Theorem~\ref{th-main2}, we have listed four properties for the
differential Chow form. In this section, we are going to prove that
these properties are also the sufficient conditions for  a
$\delta$-polynomial $F(\bu_{0},\ldots,\bu_d)$ to be the Chow form
for a $\delta$-variety.
Based on these sufficient conditions, we can define the
$\delta$-Chow quasi-variety for certain classes of
$\delta$-varieties in the sense that  a point in the $\delta$-Chow
quasi-variety represents a $\delta$-variety in the class. In other
words, we give a parametrization of all $\delta$-varieties in the
class. Obviously, this is an extension of the concept of Chow
variety in algebraic case \cite{gelfand,hodge}.

\subsection{Sufficient conditions for a polynomial to be a differential Chow form }
The following result gives sufficient conditions for a
$\delta$-polynomial to be the Chow form of an irreducible
$\delta$-variety. From Theorem~\ref{th-main2}, they are also
necessary conditions.

\begin{theorem} \label{th-suff}
Let $F(\bu_{0},\bu_{1},\ldots,\bu_{d})$ be an irreducible
$\delta$-polynomial in $\mathcal{F}\{\bu_0,\bu_1,$ $\ldots,\bu_d\}$
where $\bu_{i}=(u_{i0},u_{i1},\ldots,u_{in})\,(i=0,\ldots,d)$. If
$F$ satisfies
the following conditions, then it is the Chow form for an
irreducible $\delta$-variety of dimension $d$ and order $h$.

1.
 $F(\bu_{0},\bu_{1},\ldots,\bu_{d})$ is $\delta$-homogenous of the same degree in each set of $\delta$-variables $\bu_{i}$
 and $\ord(F,u_{ij})=h$ for all $u_{ij}$ occurring in $F$.

2.
$F(\bu_{0},\bu_{1},\ldots,\bu_{d})$  can be factored uniquely into
the following form
\begin{eqnarray}F(\bu_{0},\bu_{1},\ldots,\bu_{d})&=&A(\bu_{0},\bu_{1},\ldots,\bu_{d})
 \prod^g_{\tau=1}(u_{00}^{(h)}+\sum_{\rho=1}^n u_{0\rho}^{(h)}\xi_{\tau \rho}+t_{\tau}) \nonumber \\&=&A(\bu_{0},\bu_{1},\ldots,\bu_{d})
 \prod^g_{\tau=1}(u_{00}+\sum_{\rho=1}^n u_{0\rho}\xi_{\tau \rho})^{(h)} \nonumber
\end{eqnarray}
where $g=\deg(F,u_{00}^{(h)})$ and $\xi_{\tau \rho}$ are in a
differential extension field $\ff_\tau$ of
$\ff$.
The first $``="$ is obtained by factoring
$F(\bu_{0},\bu_{1},\ldots,\bu_{d})$ as an algebraic polynomial in
the variables $u_{00}^{(h)},u_{01}^{(h)},\ldots,u_{0n}^{(h)}$, while
the second one is a differential expression 
by defining the derivatives of $\xi_{\tau\rho}$ to be
$$\xi_{\tau\rho}^{(m)}=(\delta\xi_{\tau\rho}^{(m-1)})|_{u_{00}^{(h)}=-\sum_{\rho=1}^n u_{0\rho}^{(h)}\xi_{\tau
\rho}-t_{\tau}}\,(m \geq1)$$ recursively.

3. $\Xi_\tau=(\xi_{\tau1},\ldots,\xi_{\tau n})\,(\tau=1,\ldots,g)$
 are on the $\delta$-hyperplanes $\P_{\sigma}=0\,(\sigma=1,\ldots,d)$
 as well as on the algebraic hyperplanes
 $^a\P_{0}^{(k)}=0\,(k=0,\ldots,h-1)$.

4. For each $\tau$, if
$v_{i0}+v_{i1}\xi_{\tau1}+\cdots+v_{in}\xi_{\tau n}=0$
$(i=0,\ldots,d)$, then $F(\bv_{0},\ldots,\bv_{d})=0$, where
$\bv_{i}=(v_{i0},v_{i1},\ldots,v_{in})$ and $v_{ij}\in\mathcal{E}$.
Equivalently, if $S^0,\ldots,S^d$ are $(n+1)\times (n+1)$
skew-symmetric matrices, each having independent
$\delta$-indeterminates above its principle diagonal, then
$F(S^0\xi_\tau,\ldots,S^d\xi_{\tau})=0,$ where
$\xi_\tau=(1,\xi_{\tau1},\ldots,\xi_{\tau n}).$
\end{theorem}

Before proving the theorem, we need several lemmas.
\begin{lemma} \cite[p.11, Theorem 1]{hodge1}
\label{le-algisomorphism extension} Let $\mathcal {R}$ and $\mathcal
{S}^*$ be two  rings and $\mathcal {R}$  isomorphic to a
 subring $\mathcal {S}$ of $\mathcal {S}^*$. Then there
exists an extension ring $\mathcal {R}^*$ of $\mathcal {R}$ such
that this isomorphism can be extended to an isomorphism between
$\mathcal {R}^*$ and $\mathcal {S}^*$.

\end{lemma}
\begin{lemma}\label{le-genericpointextension}
Let $V$ be an irreducible $\delta$-$\mathcal {F}$-variety of
dimension $d>0$
 and $\P=u_0+u_1y_1+\cdots+u_ny_n$ a generic
$\delta$-hyperplane where $u_i\in\ee$. Then every generic point of
$\mathbb{V}\big(\mathbb{I}(V),\P\big)$ over $\ff_1=\mathcal
{F}\langle u_0,\ldots,u_n\rangle$ is a generic point of $V$ over
$\mathcal {F}.$
\end{lemma}
\proof By Theorem~\ref{th-inter1}, $[\mathbb{I}(V),\P]$ is a prime
$\delta$-ideal of dimension $d-1$ in $\ff_1\{\Y\}$. Let $\eta$ be a
generic point of $\mathbb{V}\big(\mathbb{I}(V),\P\big)$. Then for
any $\delta$-polynomial $p$ in $\mathbb{I}(V)$, we have $p(\eta)=0$.
On the other hand, for any $\delta$-polynomial $p\in \mathcal
{F}\{y_1,\ldots,y_n\}$ such that $p(\eta)=0$, we have $p\in
[\mathbb{I}(V),\P]$. Then $p \equiv  \sum_i h_i \P^{(i)}
\,\mod\,\mathbb{I}(V)$. Substituting $u_0$ by
$-u_1y_1-\cdots-u_ny_n$ in the above equality, we have $p \equiv
0\,\mod \,\mathbb{I}(V)$. Hence $\eta$ is a generic point of
$V$.\qed

In the next result, we will show that if condition 4) from Theorem \ref{th-suff} holds, then the following stronger version
is also valid.

\begin{lemma}\label{lm-sfc4}
Let $F$ satisfy condition 4) of Theorem \ref{th-suff}. Consider $F$  as an algebraic polynomial
$f(u_{\sigma j}^{( k)},u_{0j}^{(l)},$ $u_{00}^{(h)},\ldots,$
$u_{0n}^{(h)})$ in $u_{i j}^{(k)}$ and $(\xi_{\tau
1},\ldots,\xi_{\tau n},\ldots,\xi_{\tau 1}^{(h)},\ldots,\xi_{\tau
n}^{(h)})$ is regarded  as an algebraic point.
If $w_{i0k}+\sum_{j=1}^n \sum_{m=0}^k{k \choose m} w_{i
jm}\delta_\tau^{(k-m)}\xi_{\tau j}=0\,(i=0,\ldots,d;k=0,\ldots,h)$,
then $f(w_{\sigma jk},w_{0jl},w_{00h},$ $\ldots,w_{0nh})=0,$ where
the $w_{ijk}$ are considered as elements in the underlying ordinary
field of $\ee$.
\end{lemma}
\proof Regard $\mathcal {Q}_\tau=[v_{i0}+\sum_{j=1}^n
v_{ij}\xi_{\tau j}:i=0,\ldots,d]$ as a $\delta$-ideal in $\mathcal
{F}\langle \xi_{\tau
j}:j=1,\ldots,n\rangle\{v_{i0},\ldots,v_{in}:i=0,\ldots,d\}$, where
$v_{ij}\in\ee$ are  $\delta$-$\ff$-indeterminates. From condition 4)
of Theorem \ref{th-suff},
$F(\bv_0,\ldots,\bv_d)|_{\mathbb{V}(\mathcal {Q}_\tau)}\equiv 0$. It
is clear that $\mathcal {Q}_\tau$ is a prime $\delta$-ideal and
$\{v_{i0}+\sum_{j=1}^n v_{ij}\xi_{\tau j}:i=0,\ldots,d\}$ is its
characteristic set with $v_{i0}$ as leaders. By the differential
Nullstellensatz, $F(\bv_0,\ldots,\bv_d)\in \mathcal {Q}_\tau$. From
condition 1) of Theorem \ref{th-suff}, $\ord(F,v_{i0})=h$. Then
$F(\bv_0,\ldots,\bv_d)\in(v_{i0}+\sum_{j=1}^n v_{ij}\xi_{\tau
j},\ldots,v_{i0}^{(h)}+\sum_{j=1}^n v_{ij}^{(h)}\xi_{\tau
j}+\sum_{j=1}^n \sum_{m=0}^{h-1}{h \choose m}v_{i j}^{( m)}\xi_{\tau
j}^{(h-m)}:i=0,\ldots,d)$. Regarding the above relation as a pure
algebraic relation,  we can substitute $v_{ij}^{(k)}$ by algebraic
indeterminates $w_{ijk}$ and regard $F$ as an algebraic polynomial.
Then $f(w_{\sigma
jk},w_{0jl},w_{00h},\ldots,w_{0nh})\in(w_{i00}+\sum_{j=1}^n
w_{ij0}\xi_{\tau j},\ldots,w_{i0h}+\sum_{j=1}^n w_{ijh}\xi_{\tau
j}+\sum_{j=1}^n \sum_{m=0}^{h-1}{h \choose m}w_{i j m}\xi_{\tau
j}^{(h-m)}:i=0,\ldots,d)$, which shows that lemma is valid.\qed

\vskip 10pt {\noindent\bf\em Proof of Theorem \ref{th-suff}}.  Let
$V_{\tau}\,(\tau=1,\ldots,g)$ be the irreducible $\delta$-$\mathcal
{F}$-variety with $(\xi_{\tau_1},\ldots,\xi_{\tau n})$ as its
generic point over $\mathcal {F}$. We will show later that all the
$\delta$-varieties $V_{\tau}$ are the same.

Firstly, we claim that the generic points of $V_{\tau}$ which lie on
$\P_1,\ldots,\P_d$ as well as on $^a\P_0,^a\P'_0,\ldots,$
$^a\P^{(h-1)}_0$ are included in $\{(\xi_{\tau1},\ldots,\xi_{\tau
n}):\, \tau=1,\ldots,g\}$. Without loss of generality, we consider
$V_1$. To prove the claim holds for $V_1$, we need to use the following assertion:

 ($\star$) If $(\eta_{10},\ldots,\eta_{n0},\ldots,$ $\eta_{1h},\ldots,\eta_{nh})$
is a generic point of the algebraic ideal
$\mathbb{I}(V_1)^{<h>}=\mathbb{I}(V_1)\cap \mathcal
{F}[y_1,\ldots,y_n,\ldots,y_1^{(h)},\ldots,$ $y_{n}^{(h)}]$ which lies
on
$^a\P_\sigma^{(k)},^a\P^{(l)}_0\,(\sigma=1,\ldots,d;k=0,\ldots,h;l=0,\ldots,h-1)$,
then there must exist some $\tau$ such that $\eta_{j0}=\xi_{\tau j}$
for $j=1,\ldots,n$.

Assume that ($\star$) is valid. Suppose $\eta=(\eta_1,\ldots,\eta_n)\in\ee^n$ is a generic point of $V_1$ which lie on
$\P_1,\ldots,\P_d$ as well as on $^a\P_0,^a\P'_0,\ldots,$$^a\P^{(h-1)}_0$.
Then the algebraic point $(\eta_1,\ldots,\eta_n,\ldots,\eta_1^{(h)},\ldots,\eta_n^{(h)})$ is a generic point of
the algebraic ideal $\mathbb{I}(V_1)^{<h>}=\mathbb{I}(V_1)\cap \mathcal
{F}[y_1,\ldots,y_n,\ldots,y_1^{(h)},\ldots,$ $y_{n}^{(h)}]$ which lies on
$^a\P_\sigma^{(k)},^a\P^{(l)}_0\,(\sigma=1,\ldots,d;k=0,\ldots,h;l=0,\ldots,h-1)$. By ($\star$),
there must exist some $\tau$ such that $\eta_{j}=\xi_{\tau j}$
for $j=1,\ldots,n$. Since $\eta$ is a differential point of $V_1$, thus $\eta=\xi_\tau$ for some $\tau$.

Now we are going to prove ($\star$). Similar to the proof of Lemma \ref{lm-sfc4},
rewrite $F$ as an algebraic polynomial $f(u_{\sigma j}^{(
k)},u_{0j}^{(l)},$ $u_{00}^{(h)},\ldots,$ $u_{0n}^{(h)})$ and
consider the condition 2) as a pure algebraic factorization.
Let $(\eta_{10},\ldots,\eta_{n0},\ldots,$
$\eta_{1h},\ldots,\eta_{nh})$ be such a generic point of
$\mathbb{I}(V_1)^{<h>}$ other than $(\xi_{11},\ldots,\xi_{1
n},\ldots,\xi_{11}^{(h)},$ $\ldots,\xi_{1 n}^{(h)})$. Then we have
the following $\ff$-isomorphism $\mathcal
{F}(\eta_{10},\ldots,\eta_{n0},\ldots,\eta_{1h},\ldots,\eta_{nh})$
$\cong \mathcal {F}(\xi_{11},\ldots,\xi_{1 n},$
$\ldots,\xi_{11}^{(h)},$ $\ldots,\xi_{1 n}^{(h)})$ which maps
$\eta_{jk}$ to $\xi_{1 j}^{(k)}$ for $j=1,\ldots,n$ and
$k=0,\ldots,h$. By Lemma~\ref{le-algisomorphism extension}, there
exist $w_{\sigma jk}, w_{0jl}\in\ee$ such that the above algebraic
isomorphism can be extended to the isomorphism {\tiny
$$\mathcal{F}(\eta_{10},\ldots,\eta_{n0},\ldots,\eta_{1h},\ldots,\eta_{nh},u_{\sigma
j}^{(k)},u_{0j}^{(l)})\cong
 \mathcal {F}(\xi_{11},\ldots,\xi_{1n},\ldots,\xi_{11}^{(h)},\ldots,\xi_{1 n}^{(h)},w_{\sigma jk}, w_{0jl})$$}
where $\sigma=1,\ldots,d;j=0,\ldots,n;k=0,\ldots,h;l=0,\ldots,h-1$,
and $u_{\sigma j}^{(k)}$ and $u_{0j}^{(l)}$ map to $w_{\sigma jk}$
and $w_{0jl}$ respectively. Since
$(\eta_{10},\ldots,\eta_{n0},\ldots,\eta_{1h},\ldots,\eta_{nh})$
lies on $^a\P_\sigma^{(k)},^a\P_0^{(l)}$, the relation $u_{\sigma
0}^{(k)}+\sum_{j=1}^n \sum_{m=0}^k{k \choose m} u_{\sigma j}^{(m)}$
$\eta_{j,k-m}=0$ implies  that $w_{\sigma0k}+\sum_{j=1}^n
\sum_{m=0}^k{k \choose m}w_{\sigma j m}\xi_{1j}^{(k-m)}$
$=0\,(\sigma=1,\ldots,d)$ and the relation $u_{0
0}^{(l)}+\sum_{j=1}^n \sum_{m=0}^l{l \choose m} u_{0 j}^{(m)}$
$\eta_{j,l-m}=0$ implies that $w_{00l}+\sum_{j=1}^n \sum_{m=0}^l{l
\choose m}w_{0 j m}\xi_{1j}^{(l-m)}$ $=0\, (l=0,\ldots,h-1)$.
Furthermore, if $w_{00h}+\sum_{i=1}^n w_{0 i
h}\xi_{1i}^{(h)}+\sum_{j=1}^n \sum_{m=0}^{h-1}{h \choose m}$ $w_{0 j
m}\xi_{1j}^{(h-m)}=0$ is valid, then from  Lemma \ref{lm-sfc4}, it
follows that the algebraic polynomial $f(w_{\sigma j
k},w_{0jl},w_{00h},\ldots,w_{0nh})$ $=0$. Then, by the Hilbert
Nullstellensatz, when regarded as a
 polynomial in the algebraic indeterminates
$u_{00}^{(h)},\ldots,u_{0n}^{(h)}$, $f(w_{\sigma j k},$
$w_{0jl},u_{00}^{(h)},\ldots,$ $u_{0n}^{(h)})\in
(u_{00}^{(h)}+u_{01}^{(h)}\xi_{11}+\cdots+$ $u_{0n}^{(h)}\xi_{1n}$
$+\sum_{j=1}^n \sum_{m=0}^{h-1}{h \choose m}$ $w_{0 j
m}\xi_{1j}^{(h-m)})$.  Thus,
$u_{00}^{(h)}+u_{01}^{(h)}\xi_{11}+\cdots+$ $u_{0n}^{(h)}\xi_{1n}$
$+\sum_{j=1}^n \sum_{m=0}^{h-1}{h \choose m}$ $w_{0 j
m}\xi_{1j}^{(h-m)}$ divides $f(w_{\sigma j k},$
$w_{0jl},u_{00}^{(h)},\ldots,$ $u_{0n}^{(h)})$. By the above
isomorphism, $f(u_{\sigma j}^{( k)},u_{0j}^{(l)},u_{00}^{(h)},$
$\ldots,u_{0n}^{(h)})$ is divisible by
$u_{00}^{(h)}+u_{01}^{(h)}\eta_{10}+\cdots+u_{0n}^{(h)}\eta_{n0}+\sum_{j=1}^n
\sum_{m=0}^{h-1}{h \choose m}u_{0 j}^{(m)}$ $\eta_{j,h-m}$. The
first factorization expression of condition 2) shows that when
regarded as an algebraic polynomial in the variables
$u_{00}^{(h)},u_{01}^{(h)},\ldots,u_{0n}^{(h)}$, $f(u_{\sigma j}^{(
k)},u_{0j}^{(l)},u_{00}^{(h)},$ $\ldots,u_{0n}^{(h)})$ $=A
\prod_{\tau=1}^g\big(u_{00}^{(h)}+u_{01}^{(h)}\xi_{\tau1}+\cdots+u_{0n}^{(h)}$
$\xi_{\tau n}+\sum_{j=1}^n\sum_{k=0}^{h-1}{h \choose
k}u_{0j}^{(k)}\xi_{\tau j}^{(h-k)}\big)$. Thus, there exists some
$\tau$ such that $\eta_{j0}=\xi_{\tau j}\,(j=1,\ldots,n)$, which
completes the proof of the claim.

Denote the dimension and order of $V_{\tau}$ by $d_\tau$ and
$h_\tau$ respectively. We claim that $d_\tau=d$ and $h_\tau=h$.
Since $V_{\tau}$ meets $\P_1,\ldots,\P_d$ and
$(\xi_{\tau1},\ldots,\xi_{\tau n})$ are such points in their
intersection variety, by Theorem~\ref{th-inter1}, $d_{\tau}\geq d.$
If $d_\tau >d\geq0$, then $V_\tau$ meets $\P_1,\ldots,\P_d,\P_0$.
Let $(\eta_1,\ldots,\eta_n)$ be a generic point of
$\mathbb{V}(\mathbb{I}(V_\tau),\P_1,\ldots,\P_d,\P_0)$. Then by
Lemma~\ref{le-genericpointextension}, $(\eta_1,\ldots,\eta_n)$ is
also a generic point of $V_\tau$. Since $(\eta_1,\ldots,\eta_n)$
lies on $\P_0$, it also lies on
$^a\P_0,^a\P'_0,\ldots,^a\P^{(h-1)}_0$. From the above claim, there
exists some $\tau$ such that
$(\eta_1,\ldots,\eta_n)=(\xi_{\tau1},\ldots,\xi_{\tau n}).$ Thus,
$(\xi_{\tau1},\ldots,\xi_{\tau n})$ lies on $\P_0$, which implies
that $F(\bu_{0},\ldots,\bu_d)$ is a zero $\delta$-polynomial, which
is a contradiction. So $d_\tau=d.$

It remains to show that $h_\tau=h$. We first prove $h_\tau\geq h$.
Suppose the contrary, then $h_\tau \leq h-1$. Similar to  the proof
of Theorem~\ref{th-g}, we can prove that $
\mathbb{V}\big([\mathbb{I}(V_\tau),\P_1,\ldots,\P_d]^{<h>},^a\P_0,\ldots,$ $^a\P^{(h-1)}_0\big)=\emptyset
$. But $(\xi_{\tau_1},\ldots,\xi_{\tau
n})$ is an element of
$\mathbb{V}(\mathbb{I}(V_\tau),\P_1,\ldots,\P_d)$ which also lies on
$^a\P_0,^a\P'_0,\ldots,$ $^a\P^{(h-1)}_0$, which is a contradiction. Now suppose that
$h_\tau >h$, then $h_\tau-1\geq h$. From Theorems~\ref{th-g} and \ref{th-zerochowform}, every point
of $V_\tau$ which lies both on $\P_1,\ldots,\P_d$ and on
$^a\P_0,^a\P'_0,\ldots,^a\P^{(h_\tau-1)}_0$ is a generic point of
$V_\tau$. But the generic points of $V_{\tau}$ which lie on
$\P_1,\ldots,\P_d$ as well as on
$^a\P_0,^a\P'_0,\ldots,^a\P^{(h-1)}_0$ are included in
$\{(\xi_{\tau1},\ldots,\xi_{\tau n}):\, \tau=1,\ldots,g\}$.  So some
$(\xi_{\tau1},\ldots,\xi_{\tau n})$ lies on
$^a\P_0,^a\P'_0,\ldots,^a\P^{(h_\tau-1)}_0$. Since $h_\tau-1\geq h$,
we have $(\xi_{\tau1},\ldots,\xi_{\tau n})$ lies on $^a\P^{(h)}_0$,
which implies $F(\bu_{0},\ldots,\bu_d)\equiv 0$, a contradiction.
Thus, we have proved that $d_\tau=d$ and $h_\tau=h$.

Since the solutions of $V_\tau$ and $\P_1,\ldots,\P_d$ which also
lie on $^a\P_0,^a\P'_0,\ldots,^a\P^{(h-1)}_0$ are  generic points of
$V_\tau$ and these are therefore contained in
$\{(\xi_{\tau1},\ldots,\xi_{\tau n}):\, \tau=1,\ldots,g\}$. Hence,
the differential Chow form of $V_\tau$ is of the form
$$F_\tau(\bu_0,\ldots,\bu_{d})=A_\tau\prod_{\rho=1}^g \big((u_{00}+u_{01}\xi_{\tau1}+\cdots+u_{0n}\xi_{\tau n})^{(h)}\big)^{l_{\tau \rho}},$$ where
$l_{\tau \rho}=1 $ or $0$ according to whether
$(\xi_{\tau1},\ldots,\xi_{\tau n})$ is in $V_\tau$. Since both
$F_\tau$ and $F$ are irreducible, they differ at most by a factor in
$\mathcal {F}$. Therefore, $V_\tau\,(\tau=1,\ldots,g)$ are the same
$\delta$-variety, and $F(\bu_0,\ldots,\bu_{d})$ is their
differential Chow form. \qed


\vskip5pt In order to define $\delta$-Chow varieties in the next
subsection, we will introduce the concept of order-unmixed
$\delta$-varieties. A $\delta$-variety $V$ is called {\em
order-unmixed} if all its components have the same dimension and
order.
Let $V$ be an order-unmixed $\delta$-variety with dimension $d$ and
order $h$ and $V=\bigcup_{i=1}^l V_{i}$ its minimal irreducible
decomposition with $F_{i}(\bu_{0}, \bu_{1}, \ldots, \bu_{d})$ the
Chow form of $V_{i}$. Let
 \begin{equation}\label{eq-sss1}
 F(\bu_{0},\ldots,\bu_{d})=\prod_{i=1}^l F_{i}(\bu_{0}, \bu_{1}, \ldots, \bu_{d})^{s_{i}}
 \end{equation}
  with $s_{i}$ arbitrary nonnegative integers.
Associated with \bref{eq-sss1}, we introduce the concept of {\em
differential algebraic cycle}, or simply $\delta$-cycle, $\VB =
\sum_{i=1}^l s_i \VB_i$ as a generalization of the concept of
algebraic cycle in algebraic geometry \cite{gelfand}, where $s_i$ is
called the multiplicity of $\VB_i$ in $\VB$.
Recall that, we have defined the $\delta$-degree $m$ and leading
$\delta$-degree $g$ for an irreducible $\delta$-variety $V$ in
Definitions \ref{def-ddeg} and \ref{def-ldeg} respectively.
Let $g_i$ and $m_i$ be the leading $\delta$-degree and
$\delta$-degree of $V_i$ respectively. Then the leading
$\delta$-degree and $\delta$-degree of $\VB$ is defined to be
$\sum_{i=1}^l s_ig_i$ and $\sum_{i=1}^l s_im_i$ respectively.

Given a $\delta$-polynomial $G(\bu_0,\ldots,\bu_d)$ with
$\ord(G,u_{00})=h$, it may be reducible over $\mathcal {F}$ such
that some of its irreducible factors are free of $u_{00}^{(h)}$. In
that case, if the product of all such factors is $L$, then we define
the primitive part of $G$ w.r.t. $u_{00}^{(h)}$ to be $G/L$.
Otherwise, its primitive part w.r.t. $u_{00}^{(h)}$  by convention
is defined to be itself. Then we have

\begin{theorem} \label{th-suff2}
Let $F(\bu_{0},\bu_{1},\ldots,\bu_{d})$ be a $\delta$-polynomial and
$\widetilde{F}$ the primitive part of $F$ with respect to the
variable $u_{00}^{(h)}$. If $F$ satisfies the four conditions in
Theorem \ref{th-suff}, then $\widetilde{F}$ is the Chow form for
$\delta$-cycle of dimension $d$ and order $h$.
\end{theorem}

\proof By definition, $F=B\widetilde{F}$, where $\ord(B,\bu_0)<h$.
Since $F$ is $\delta$-homogenous in $\bu_i$ for each $i$,
$\widetilde{F}$ is $\delta$-homogenous in each $\bu_i$ too. And
since $B$ is free of $u_{00}^{(h)}$, i.e. $B$ divides $A$, then
$\widetilde{F}$ satisfies conditions 2) and 3), and moreover the
$(\xi_{\tau 1},\ldots,\xi_{\tau n})$ in the factorization are the
same as that of $F$. By the proof of Theorem~\ref{th-suff}, we have
$\mathbb{I}(\xi_{\tau 1},\ldots,\xi_{\tau n})$ is of dimension $d$
and order $h$ over $\mathcal {F}$. Then similarly as the proof of
Lemma~\ref{lem-d} and Theorem \ref{th-choworder}, we conclude that
$\mathbb{I}(-\sum_{j=1}^n v_{0j}\xi_{\tau j},\ldots,-\sum_{j=1}^n
v_{dj}\xi_{\tau j})$ is of dimension $d$  over $\mathcal {F}\langle
v_{ij}:i=0,\ldots,d;j=1,\ldots,n\rangle$ and its relative order
w.r.t. any parametric set is $h$, where
$v_{ij}\,(i=0,\ldots,d;j=1,\ldots,n)$ are $\delta$-indeterminates
over $\mathcal {F}\langle \xi_{\tau 1},\ldots,\xi_{\tau n}\rangle$.
In particular, $\trdeg\,\mathcal {F}\langle
\zeta_0,\zeta_1,\ldots,\zeta_d\rangle/\mathcal {F}\langle
\zeta_1,\ldots,\zeta_d\rangle=h,$ where $\zeta_i=-\sum_{j=1}^n
v_{ij}\xi_{\tau j}$. Thus  $B(-\sum_{j=1}^n v_{0j}\xi_{\tau
j},\ldots,-\sum_{j=1}^n v_{dj}\xi_{\tau j})\ne0.$ But
 \newline
$F(-\sum_{j=1}^n v_{0j}\xi_{\tau j},\ldots,-\sum_{j=1}^n
v_{dj}\xi_{\tau j})=0,$ so $\widetilde{F}(-\sum_{j=1}^n
v_{0j}\xi_{\tau j},\ldots,-\sum_{j=1}^n$ $ v_{dj}\xi_{\tau j})$
$=0$. It follows that $\widetilde{F}(S^0\xi_\tau,\ldots,S^d
\xi_\tau)=0$, for if we suppose the contrary, then
$B(S^0\xi_\tau,\ldots,S^d \xi_\tau)=0$. But if we specialize
$s^i_{jk}(j<k,j>0)$ to 0 and $s^i_{0k}\,(k>0)$ to $-v_{ik}$, then
$B(-\sum_{j=1}^n v_{0j}\xi_{\tau j},\ldots,-\sum_{j=1}^n
v_{dj}\xi_{\tau j})=0,$ which is a contradiction. Thus,
$\widetilde{F}$ satisfies condition 4).

Now we claim that $\widetilde{F}$ is the Chow form of some
$\delta$-cycle.
Let $V_{\tau}\,(\tau=1,\ldots,g)$ be the irreducible
$\delta$-$\mathcal {F}$-variety  with
$(\xi_{\tau_1},\ldots,\xi_{\tau n})$ as its generic point over
$\mathcal {F}$.
Following the steps in the proof of Theorem~\ref{th-suff} exactly,
we arrive at the conclusion that the Chow form of $V_\tau$ is of the
form
$$F_\tau(\bu_0,\ldots,\bu_{d})=A_\tau\prod_{\rho=1}^g \big((u_{00}+u_{01}\xi_{\tau1}+\cdots+u_{0n}\xi_{\tau n})^{(h)}\big)^{l_{\tau \rho}},$$ where
$l_{\tau \rho}=1 $ or $0$ according to whether
$(\xi_{\tau1},\ldots,\xi_{\tau n})$ is in $V_\tau$.  Since each
$\xi_\tau$ is in at least one of the $\delta$-varieties $V_i$, the
Chow form of $\bigcup_{\tau=1}^g V_\tau$ is of the form
$G(\bu_0,\ldots,\bu_d)=\prod_{\tau=1}^g
(F_\tau(\bu_0,\ldots,\bu_{d}))^{s_\tau}=C(\bu_0,\ldots,\bu_d)\prod_{\rho=1}^g
\big((u_{00}+u_{01}\xi_{\tau1}+\cdots+u_{0n}\xi_{\tau
n})^{(h)}\big)^{\eta_{\tau \rho}}$ with $\eta_{\tau \rho}>0$. Since
$\widetilde{F}$ and $G$ have the same factors
$(u_{00}+u_{01}\xi_{\tau1}+\cdots+u_{0n}\xi_{\tau n})^{(h)}$ and the
primitive factor of $\widetilde{F}$ w.r.t. $u_{00}^{(h)}$ is itself,
thus we can find $\eta_{\tau \rho}$ such that $\widetilde{F}=G$
which completes the proof.   \qed

\subsection{Differential Chow quasi-variety for a differential algebraic cycle}
A $\delta$-cycle $V$ in the $n$ dimensional $\delta$-affine space
with dimension $d$, order $h$, leading $\delta$-degree $g$, and
$\delta$-degree $m$ is said to be of index $(n,d,h,g,m)$.
In this section, we will define the $\delta$-Chow quasi-variety in
certain cases such that each point in this $\delta$-variety
represents a $\delta$-cycle with a given index $(n,d,h,g,m)$.

For a given index $(n,d,h,g,m)$, a $\delta$-polynomial
$F(\bu_0,\ldots,\bu_d)$ which has unknown coefficients
$a_\lambda\,(\lambda=0,\ldots,D)$ and satisfies the following two
conditions is referred to as a {\em $\delta$-polynomial with index
$(n,d,h,g,m)$}.

1) $F$ is a homogenous polynomial of the same degree
$m$ in each set of indeterminates
$\bu_i=(u_{i0},u_{i1},\ldots,u_{in})\,(i=0,\ldots,d)$ and their derivatives.
Furthermore, for each $u_{ij}$, $\ord(F,u_{ij})$ is either $h$ or $-\infty$. In
particular, $\ord(F,u_{00})=h$.

2) As a polynomial in
$u_{00}^{(h)},u_{01}^{(h)},\ldots,u_{0n}^{(h)}$, its total degree is
$g$. In particular, $\deg(F,u_{00}^{(h)})=g$.

We want to determine the necessary and sufficient conditions imposed
on $a_\lambda\,(\lambda=0,\ldots,D)$ in order that $F$ is the Chow
form for a $\delta$-cycle with index $(n,d,h,g,m)$.
Proceeding in this way, if the necessary and sufficient conditions
given in Theorem \ref{th-suff} can be expressed by some
$\delta$-polynomials in $a_\lambda$, then the $\delta$-variety
defined by them is called the $\delta$-Chow (quasi)-variety. More
precisely, we have the following definition.
\begin{definition}
Let $F(\bu_0,\ldots,\bu_d)$ be a $\delta$-polynomial with
$\delta$-$\ff$-indeterminates $a_i\,(i=0,\ldots,D)$ in $\ee$ as
coefficients and with index $(n,d,h,g,m)$.
A quasi-$\delta$-variety $\CV$ in the variables $a_i$ is called the
{\em differential Chow quasi-variety} with index $(n,d,h,g,m)$ if a
point $\bar{a}_i$ is in $\CV$ if and only if $\widetilde{F}$ is the
Chow form for a $\delta$-cycle with index $(n,d,h,g,m_1)$ with
$m_1\le m$, where $\widetilde{F}$ is obtained from $F$  by first
replacing $a_i$ by $\bar{a}_i$  and then taking the primitive part
with respect to the variable $u_{00}^{(h)}$.
\end{definition}

In the case $h=0$, since Theorems~\ref{th-main2} and \ref{th-suff}
become their algebraic counterparts, we can obtain the equations for the algebraic Chow variety
in the same way as in \cite[p.56-57]{hodge}.
So in the following, we only consider the case $h>0$.
For $h>0$, the case $g=1$ is relatively simple. The following result
shows how to determine the defining equations for the $\delta$-Chow
quasi-variety with index $(n,d,h,g,m)$ in the case  $g=1$.

\begin{theorem}
\label{th-cv1} Let $F(\bu_0,\ldots,\bu_d)$ be a $\delta$-polynomial
with $\delta$-$\ff$-indeterminates $a_\nu\,(\nu=0,\ldots,D)$ as
coefficients and with index $(n,d,h,g,m)$ with $g=1$. Let $I_F$ be
the initial of $F$ w.r.t. the elimination ranking $u_{00}\succ
u_{i,j}$ and $a_0,\ldots,a_I$ the coefficients of $I_F$.
Then we can find a set of $\delta$-homogeneous $\delta$-polynomials
$$R_{\omega}(a_0,\ldots,a_D)\,(\omega=1,\ldots,\upsilon)$$
in $a_\nu$ such that $\V(R_{\omega}:
\omega=1,\ldots,\upsilon)\setminus \V(a_0,\ldots,a_I)$ is the
$\delta$-Chow quasi-variety of index $(n,d,h,g,m)$ with $g=1$.
%
\end{theorem}
\proof In order for $F$ to be a differential Chow form, by Theorem
\ref{th-main2}, $F$ must be $\delta$-homogeneous in each $\bu_i$.
Let $\lambda$ be a $\delta$-$\ff$-indeterminate. For each $i$,
replacing $\bu_i$ by $\lambda\bu_i$ in $F$, we should have
$$F(\bu_0,\ldots,\bu_{i-1},\lambda\bu_i,\bu_{i+1},\ldots,\bu_d)=\lambda^{m}F(\bu_0,\ldots,\bu_d).$$
Comparing the coefficients of the power products of $\lambda,
u_{ij}$ and their derivatives, we obtain a system of linearly
homogenous equations
$R_{\omega}(a_0,\ldots,a_D)=0,\,(\omega=1,\ldots,e_1)$ in $a_\nu$,
which are the conditions for $F$ to be $\delta$-homogeneous and with
degree $m$ in each $\bu_i$. So by Gaussian elimination in linear
algebra, we can obtain a basis for the solution space of
$R_{\omega}=0\,(\omega=1,\ldots,e_1)$. More precisely, if the
coefficient matrix of this linear equations is of rank $r$, then $r$
of $\{a_0,\ldots,a_D\}$ are the linear combinations of the other
$D+1-r$ of $a_\nu$. Now substitute these $r$ relations into $F$ and
denote the new $\delta$-polynomial by $F_1$. That is, $F_1$ is a
$\delta$-homogenous $\delta$-polynomial in each $\bu_i$, which only
involves $D+1-r$ independent coefficients $a_\nu$.

Since $g=1$, $F_1$ can be written in the form
$$F_1(\bu_0,\ldots,\bu_d)=A_0
u_{00}^{(h)}+A_1u_{01}^{(h)}+\cdots+A_nu_{0n}^{(h)}+B,$$ where $A_i$
and $B$ are free of $u_{0k}^{(h)}$. Denote
$-(A_1u_{01}^{(h)}+\cdots+A_nu_{0n}^{(h)}+B)/A_0$ by $\gamma$. Then
$u_{00}^{(h)}=\gamma$ is the solution of $F_1$ as an algebraic
polynomial in $u_{00}^{(h)}$. Let $\xi_j=\frac{\partial
F_1}{\partial u_{0 j}^{(h)}}\Big/\frac{\partial F_1}{\partial u_{0
0}^{(h)}}\Big|_{u_{00}^{(h)}=\gamma}=A_j/A_0\big|_{u_{00}^{(h)}=\gamma}$
for $j=1,\ldots,n$. Proceeding as in the proof of
Theorem~\ref{th-fac1}, we define the derivatives of $\xi_j$ to be
$\xi_{j}^{(k)}=(\delta \xi_{j}^{(k-1)})\big|_{u_{00}^{(h)}=\gamma}$.
%
%
It is easy to see that this definition is well defined.
Since $F_1$ is $\delta$-homogeneous in $\bu_0$, by
Theorem~\ref{th-dhomo}, for $r\neq0$
$$\sum_{j=0}^{n}\sum_{k\geq0}{k+r \choose r}u_{0 j} ^{(k)}\frac{\partial
F_1}{\partial u_{0 j}^{(k+r)}} =0.$$ In the case $r=h$, we have
$\sum_{j=0}^{n}u_{0 j}\frac{\partial F_1}{\partial u_{0
j}^{(h)}}=0$.
Set $u_{00}^{(h)}=\gamma$ in the identity $\sum_{j=0}^{n}u_{0
j}\frac{\partial F_1}{\partial u_{0 j}^{(h)}}=0$, then we have
$u_{00}+\sum_{j=1}^{n}u_{0 j}\xi_j=0$ with $u_{00}^{(h)}=\gamma$. So
$(\xi_1,\ldots,\xi_n,\ldots,\xi_{1}^{(h-1)},\ldots,\xi_{n}^{(h-1)})$
is a solution of $^a\P_0,^a\P_0',\ldots,^a\P_0^{(h-1)}$ and
$\gamma=-(\sum_{j=1}^{n}u_{0 j}\xi_j)^{(h)}$. So
$F_1(\bu_0,\ldots,\bu_d)=A_0(u_{00}+\sum_{j=1}^{n}u_{0j}\xi_j)^{(h)}$.
As a consequence, with these $\xi_i$, the second condition and the second part of the
third condition in Theorem \ref{th-suff} are satisfied.

In order for $F_1$ to be the Chow form for some $\delta$-variety, by
Theorem \ref{th-suff}, $(\xi_1,\ldots,\xi_n)$ should satisfy
$\P_{\sigma}=0\,(\sigma=1,\ldots,d)$ and $F_1(S^0\xi,\ldots,S^d\xi)$
$=0$ where $S^i$ are $(n+1)\times (n+1)$-skew symmetric matrices
with elements independent indeterminates and
$\xi=(1,\xi_1,\ldots,\xi_n)^T$.

Firstly, setting $y_j=A_j/A_0$ in $\P_{\sigma}=0$, we get
$u_{\sigma0}A_0+\sum_{j=1}^n u_{\sigma j}A_j=0$. Then we obtain some
equations in $a_\nu$ by equating to zero the coefficients of the
various $\delta$-products of $\bu_0,\ldots,\bu_d$. This gives
$\delta$-polynomials
$R_{\omega}(a_0,\ldots,a_D)\,(\omega=e_1+1,\ldots,e_2)$.

Secondly,  we obtain some $\delta$-equations
$\chi_{\tau}(a_\nu,y_1,$ $\ldots,y_n)$ by equating to zero the
coefficients of all $\delta$-products of the independent
indeterminates $s^i_{jk}\,(j>k)$ in $F_1(S^0Y,\ldots,$ $S^dY)=0$
with $Y=(1,y_1,\ldots,y_n)^T$.
Then setting $y_j^{(k)}=\xi_{j}^{(k)}$ in the above $\chi_{\tau}$
and clearing denominators, we obtain polynomial equations $p_{\mu}$
in $u_{ij}^{(k)}$ and $a_\nu$. Equating to zero the coefficients of
the power products of the $u_{ij}^{(k)}$ in $p_{\mu}$, we finally
obtain $\delta$-polynomials in $a_\nu$:
$R_{\omega}(a_0,\ldots,a_D)\,(\omega=e_2+1,\ldots,\upsilon)$.
We then obtain the defining equations $R_{\omega}=0$ for the Chow variety.

We now show that all the $R_{\omega}$ are $\delta$-homogenous
$\delta$-polynomials. Denote $\ba=(a_0,\ldots,a_D)$. We have known
for $\omega=1,\ldots,e_1$, $R_{\omega}$ are linearly homogenous
$\delta$-polynomials in $\ba$. Since $F_1$ as well as $A_i$ are
linearly homogenous in $\ba$,
$R_{\omega}(\ba)\,(\omega=e_1+1,\ldots,e_2)$ are linearly homogenous
$\delta$-polynomials. To show
$R_{\omega}(\ba)\,(\omega=e_2+1,\ldots,\upsilon)$ are
$\delta$-homogenous $\delta$-polynomials,  by induction on $k$, we first show  that for any
$\delta$-$\ff$-indeterminate $\lambda$ and for each $j$ and $k$,
$\xi_j^{(k)}(\lambda \ba)=\xi_j^{(k)}(\ba)$.
By the expression of $\gamma$, it is clear that
$\gamma(\lambda\ba)=\gamma.$ Since
$\xi_j=A_j/A_0\big|_{u_{00}^{(h)}=\gamma}$ and
$\xi_j(\lambda\ba)=\frac{A_j(\lambda\ba)}{A_0(\lambda\ba)}\big|_{u_{00}^{(h)}=\gamma(\lambda\ba)}=\frac{\lambda
A_j(\ba)}{\lambda A_0(\ba)}\big|_{u_{00}^{(h)}=\gamma}=\xi_j(\ba)$.
So it holds for $k=0$. Suppose it holds for $k-1$, that is
$\xi_j^{(k-1)}(\lambda\ba)=\xi_j^{(k-1)}(\ba)$. Since
$\xi_{j}^{(k)}=(\delta \xi_{j}^{(k-1)})\big|_{u_{00}^{(h)}=\gamma}$,
$\xi_{j}^{(k)}(\lambda\ba)=\Big(\delta\big(
\xi_{j}^{(k-1)}(\lambda\ba)\big)\Big)\Big|_{u_{00}^{(h)}=\gamma(\lambda\ba)}=(\delta
\xi_{j}^{(k-1)}(\ba))\big|_{u_{00}^{(h)}=\gamma(\ba)}=\xi_{j}^{(k)}(\ba)$.
Now we are going to show that
$R_{\omega}(\ba)\,(\omega=e_2+1,\ldots,\upsilon)$ are
$\delta$-homogenous $\delta$-polynomials. Since $F_1$ is linearly
homogenous in $\ba$, each $\chi_{\tau}(a_\nu,y_1,\ldots,y_n)$ is
linearly homogenous in $\ba$. Setting $y_j^{(k)}=\xi_j^{(k)}$ in
$\chi_{\tau}$, it is clear that
$\chi_{\tau}(a_\nu,\xi_1,\ldots,\xi_n)(\lambda
\ba)=\lambda\cdot\chi_{\tau}(a_\nu,\xi_1,\ldots,\xi_n)$ and the
denominator of $\chi_{\tau}(a_\nu,\xi_1,\ldots,\xi_n)$ is a pure
algebraic homogenous polynomial in $\ba$. Thus, $p_\mu$ is
$\delta$-homogenous in $\ba$ and the $\delta$-homogeneity of
$R_{\omega}$ follows.

Let $a_{0},\ldots,a_I$ be the coefficients of $I_{F}$. Then we claim
that the quasi-projective $\delta$-variety
$\mathbb{CV}=\mathbb{V}(R_{\omega}:\omega=1,\ldots,\upsilon)\setminus
\mathbb{V}(a_0,\ldots,a_I)$ is the $\delta$-Chow quasi-variety.
Indeed, for every element $(\bar{a_0},\ldots,\bar{a_D})$ in
$\mathbb{CV}$, following the proof of this theorem, $\overline{F}$
with coefficients $\bar{a_\nu}$ satisfies the four conditions in
Theorem~\ref{th-suff}. And since $g=1$, its primitive part
must be irreducible and satisfies the four conditions too, which
consequently must be the Chow form for some irreducible
$\delta$-variety with index $(n,d,h,1,m_1)$ with $m_1\leq m$. \qed

The following example illustrates the procedure to compute the
$\delta$-Chow quasi-variety in the case of $g=1$.

\begin{example}\label{ex-dc2}
We consider a $\delta$-polynomial which has 16 terms and has index
$(2,1,1,1,2)$ to illustrate the proof of Theorem \ref{th-cv1}:
 $F=a_1 u_{12}^2u_{01}u_{00}'+a_2 u_{ 11}
 u_{ 12}  u_{ 02}  u'_{ 00}$ $ +a_{ 3}  u_{ 01}
 u_{ 02} u_{12}  u'_{ 10}$ $ + a_{ 4}  u_{ 02}^2  u_{ 11}  u'_{10} +a_{ 5}  u_{ 12}
^2 u_{ 00}  u'_{ 01} +a_{ 6}  u_{ 10}  u_{ 12}  u_{ 02} u'_{ 01}$
$ +a_{ 7}  u_{ 00}  u_{ 02}  u_{ 12}  u'_{ 11} +a_{ 8} u_{ 02}
^2 u_{ 10}  u'_{ 11} +a_{ 9}  u_{ 10} u_{ 01}  u_{ 12} u'_{ 02}+a_{
10} u_{ 00}  u_{ 02}  u_{ 11}
 u'_{ 12}$ $+a_{ 11}  u_{ 11}  u_{ 12}  u_{ 00}
 u'_{ 02}$ $  +a_{ 12}  u_{ 01}  u_{ 02}  u_{ 10} u'_{ 12}$ $+
 a_{13} u_{00} u_{11}^2 u_{02}+
 a_{14} u_{00}$ $ u_{01}u_{11} u_{12}+
 a_{15} u_{01} u_{10} u_{02} u_{11}+
 a_{16} u_{10} u_{01}^2  u_{12}$.
We will derive the conditions about the coefficients $a_\nu$ under
which $F$ is a Chow form.
Firstly, in order for $F$ to be $\delta$-homogenous, we have
$R_1=a_5+a_1,R_2=a_8+a_4,R_3=a_9+a_6,R_4=a_{10}+a_7,R_5=a_{11}+a_2,R_6=a_{12}+a_3$.
 Replacing $a_5, a_8, a_9, a_{10}, a_{11}, a_{12}$ by
           $-a_1,-a_4,-a_6,-a_7,-a_2,-a_3$ respectively in $F$ to obtain $F_1$.
%

 For such an $F_1$, $A_0=a_1 u_{12}^2 u_{01}+a_2 u_{11}u_{12}u_{02},\,A_1=-a_1 u_{12}^2 u_{00}+a_6 u_{10}u_{12}u_{02},\,A_2$ $=-a_6 u_{10}u_{01}u_{12}-a_{2}u_{11}u_{12}$ $u_{00},$
 and $B=a_{ 3}  u_{ 01}
 u_{ 02}  u_{ 12}  u'_{ 10} + a_{ 4}  u_{ 02}^2  u_{ 11}  u'_{ 10}+a_{ 7}  u_{ 00}  u_{ 02}  u_{ 12}$ $ u'_{ 11} -a_{ 4} u_{ 02} ^2 u_{
10}  u'_{ 11} -a_{ 7} u_{ 00}  u_{ 02}  u_{ 11}
 u'_{ 12}-a_{ 3}  u_{ 01}  u_{ 02}  u_{ 10}
 u'_{ 12}
 + a_{13} u_{00} u_{11}^2 u_{02}+
 a_{14} u_{00}$ $ u_{01}u_{11}$ $u_{12}+
 a_{15} u_{01} u_{10}$ $u_{02} u_{11}+
 a_{16} u_{10} u_{01}^2  u_{12}.$  Then $\gamma=(F_1-A_0u_{00}')/A_0$,
 $\xi_1=A_1/A_0|_{u_{00}'=\gamma}$, and $\xi_2=A_2/A_0|_{u_{00}'=\gamma}$.
 To confirm that $u_{10}+u_{11}\xi_1+u_{12}\xi_2=0$, we must have
 $R_{7}=a_6-a_1=0,R_8=a_2+a_1=0$.
%
%

In order to satisfy the fourth condition of Theorem \ref{th-suff},
we obtain a set of $\delta$-polynomial equations
$R_\omega(a_0,\ldots,a_{16})=0$ which have more complicated forms.
By simplifying them with $R_{7}=0$ and $R_8=0$, we obtain
 $R_9=a_7a_1(a_{15}+a_{16}),R_{10}=a_7a_1(a_{14}+a_{16}),R_{11}=a_7a_1(a_{13}-a_{16}),R_{12}=a_1(a_4-a_7),R_{13}=a_1(a_3+a_7)$,
 $R_{14}=a_1a_{16}(a_1-a_7)$,
 $R_{15}=a_{15}a_1^2+a_1a_7a_{16}$, $R_{16}=a_{14}a_1^2+a_1a_7a_{16}$, $R_{17}=a_{7}a_1(a_1-a_7)$, $R_{18}=a_1^2a_{13}-a_1a_7a_{16}$, $R_{19}=a_1^3-a_7^2a_1$.
 Thus the Chow quasi-variety is
 $\mathbb{V}(R_1,\ldots,R_{19})/\mathbb{V}(a_1,a_2)
 =\mathbb{V}(a_2+a_1,
 a_3+a_1,
 a_4-a_1,
 a_5+a_1,
 a_6-a_1,
 a_7-a_1,
 a_8+a_1,
 a_9+a_1,
 a_{10}+a_1,
 a_{11}-a_1,
 a_{12}-a_1,
 a_{14}+a_{13},
 a_{15}+a_{13},
 a_{16}-a_{13})/\mathbb{V}(a_1)$.
 From Example \ref{ex-cf1}, it is easy to check that
 each point of this quasi-variety is the coefficients of the Chow form for
 $\mathbb{V}(a_1y_1'+a_{13}y_2)$ for some $a_1,a_{13} \in \mathcal {F}$.
 Note that $a_{13}$ could be zero and the result is still valid.
\end{example}

We are unable to prove the existence of the Chow quasi-variety in
the case of $g>1$. The main difficulty is how to do elimination for
a mixed system consisting of both differential and algebraic
equations. In our case, conditions  2, and the second part of
condition 3 of Theorem \ref{th-suff} generate algebraic equations in
the coefficients of $F$ and $\xi_{ij}$, while condition 1, the first
part of condition 3, and condition 4 of Theorem \ref{th-suff}
generate $\delta$-equations. And we need to eliminate variables
$\xi_{ij}$ from these equations.

The following example shows that the $\delta$-Chow quasi-variety can
be easily defined in the case of $n=1$.

\begin{example}\label{ex-dc1}
If $n=1$, then $d=0$ and every irreducible $\delta$-homogeneous
$\delta$-polynomial in $\bu_0=(u_{00},u_{01})$ is the differential
Chow form for some irreducible  $\delta$-variety.
\end{example}
\proof Let $F(\bu_0)=F(u_{00},u_{01})$ be an irreducible
$\delta$-homogenous $\delta$-polynomial with degree $m$ and order
$h$. Then $F(-\frac{u_{00}}{u_{01}},-1)=(-\frac{1}{u_{01}})^m
F(u_{00},u_{01})$. Let
$g(-\frac{u_{00}}{u_{01}})=F(-\frac{u_{00}}{u_{01}},-1).$ It is easy
to show that $g(y)$ is an irreducible $\delta$-polynomial.
%
By Example \ref{ex-cf1},  the Chow form of the prime $\delta$-ideal
$\sat(g(y))$ is $(-u_{01})^mg(-\frac{u_{00}}{u_{01}}) = (-u_{01})^m
F(-\frac{u_{00}}{u_{01}}, -1)$ $ = F(u_{00}, u_{01})$, and the
result is proved. \qed

As a consequence, the Chow quasi-variety in the case of $n=1$ always
exists.
\begin{example}\label{ex-dc3}
Let $\bu_0=(u_{00},u_{01})$ and $F(\bu_0)$ a homogenous
$\delta$-polynomial with index $(1,0,h,g,m)$ and coefficients
$a_0,\ldots,a_D$. Let $I_F$ be the initial of $F$ w.r.t. the
elimination ranking $u_{00}\succ u_{01}$ and $a_0,\ldots,a_I$ the
coefficients of $I_F$.
Let $R_{\omega}(a_0,\ldots,$ $a_D)\,(\omega=1,\ldots,e_1)$ be the
equations obtained in the proof of Theorem \ref{th-cv1}. Then under
the condition $R_{\omega}(a_0,\ldots,a_D)=0\,(\omega=1,\ldots,e_1)$,
$F$ will become a $\delta$-homogenous $\delta$-polynomial. Then by
Example \ref{ex-dc1},  $\V(R_{\omega}: \omega=1,\ldots,e_1)\setminus
\V(a_0,\ldots,a_I)$ is the $\delta$-Chow quasi-variety with index
$(1,0,h,g,m)$.
\end{example}

\section{Generalized differential Chow form and differential resultant}
\label{sec-gdcf} We mentioned that the differential Chow form can be
obtained by intersecting the $\delta$-variety with generic
$\delta$-hyperplanes. In this section, we show that when
intersecting an irreducible $\delta$-variety of dimension $d$ by
$d+1$ generic $\delta$-hypersurfaces, we can obtain the generalized
Chow form which has similar properties to the Chow form.
As a direct consequence, we can define the differential resultant
and obtain its properties.

\subsection{Generalized differential Chow form}
Let $V\subset\ee^n$ be an irreducible $\delta$-$\ff$-variety with
dimension $d$ and order $h$,  and
\begin{equation}\label{eq-gpol}
 \P_i=u_{i0}+\sum_{j=1}^n\sum_{k=0}^{s_i} u_{ijk}y_{j}^{(k)} +\sum_{
 \begin{array}{c} \alpha \in \mathbb{Z}^{n(s_i+1)}_{\geq 0} \\ 1<|\alpha|\leq
 m_i
 \end{array}}u_{i \alpha}(\Y^{(s_i)})^{\alpha},\, (i=0,\ldots,d)
\end{equation}
a generic $\delta$-polynomial of order $s_i\ge0$ and degree
$m_i\ge1$, where  $u_{i0}, u_{ijk},u_{i\alpha}\,$ $(i=0,  \ldots,
d;j=1, \ldots, n;k=0,\ldots,s_i;\alpha \in
\mathbb{Z}^{n(s_i+1)}_{\geq 0},1<|\alpha|\leq m_i)$  are
$\delta$-$\ff$-indeterminates in $\ee$ and $(\Y^{[s_i]})^{\alpha}$
is a monomial in $\mathcal {F}[\Y^{[s_i]}]$ with exponent vector
$\alpha=(\alpha_{10},\ldots,\alpha_{n0},\alpha_{11},$
$\ldots,\alpha_{n1},\ldots,$ $\alpha_{1s_i},\ldots,\alpha_{n s_i})$,
i.e. $(\Y^{[s_i]})^{\alpha}=\prod \limits_{j=1}^n \prod
\limits_{k=0}^{s_i} (y_{j}^{(k)})^{\alpha_{jk}}$ and
$|\alpha|=\sum_{j=1}^n\sum_{k=0}^{s_i} \alpha_{jk}$.  For
convenience in the rest of the paper, we denote the nonlinear part
of each $\P_i$ by $f_i$, that is,
 $$\P_i=u_{i0}+\sum_{j=1}^n\sum_{k=0}^{s_i} u_{ijk}y_{j}^{(k)}+f_i\,(i=0,\ldots,d)$$

Denote $\bu$ to be the set consisting of all the $u_{ijk}$ and
$u_{i\alpha}$ for $i=0,\ldots,d$. Let $\xi=(\xi_{1}, \ldots,
\xi_{n})\subset\ee^n$  a generic point of $V$, which is free from
$\ff\langle\bu,u_{00},\ldots,u_{d0}\rangle$.
We define $d+1$ elements $\zeta_{0}, \zeta_{1}, \ldots, \zeta_{d}$
of $\mathcal {F}\langle\bu, \xi_{1}, \ldots, \xi_{n} \rangle$:
\begin{eqnarray} \label{equ12}
\zeta_{i}=- \sum_{j=1}^n \sum_{k=0}^{s_i} u_{ijk}\xi_{j}^{(k)} -
f_i(\xi_{1}, \ldots,    \xi_{n})\, (i=0, \ldots,    d)
\end{eqnarray}

Similar to the proof of Lemma \ref{lem-d}, we can prove that if
$d>0$ then any set of $d$ elements of $\zeta_0,\ldots,\zeta_d$ is a
$\delta$-transcendence basis of $\mathcal {F}\langle\bu\rangle
\langle \zeta_{0},    \ldots,    \zeta_{d}\rangle$ over $\mathcal
{F}\langle\bu\rangle$. If $d=0$, $\zeta_0$ is $\delta$-algebraic
over $\mathcal{F}\langle\bu\rangle$.
We thus have
\begin{lemma} \label{lem-dd}
 $\dtrdeg \,\mathcal {F}\langle\bu\rangle \langle
\zeta_{0},    \ldots,    \zeta_{d} \rangle/\mathcal {F}\langle
\bu\rangle=d$.
\end{lemma}
%
%

Let $\I_{\zeta}$ be the prime $\delta$-ideal in $\mathcal{R} =
\mathcal {F}\langle\bu \rangle\{z_{0}, \ldots, z_{d}\}$ having
$\zeta=(\zeta_{0}, \ldots, \zeta_{d})$ as a generic point. By
Lemma~\ref{lem-dd}, the dimension of $\I_{\zeta}$ is $d$. Then, the
characteristic set of
 $\I_{\zeta}$ w.r.t. any ranking  consists of an irreducible $\delta$-polynomial $g(z_{0}, \ldots,  z_{d})$ in $\mathcal {R}$ and
$\I_{\zeta} = \sat(g)$.
Since the coefficients of $g(z_{0}, \ldots,z_{d})$ are elements in
$\mathcal {F}\langle\bu \rangle$, without loss of generality, we
assume that $g(\bu; z_{0}, \ldots,  z_{d})$ is irreducible in
$\mathcal {F}\{ \bu; z_{0}, \ldots, z_{d}\}$.
We shall subsequently replace $z_{0}, \ldots, z_{d}$ by $u_{00},
\ldots, u_{d0}$,   and obtain
\begin{equation}\label{eq-cf1}
G(\bu_{0},\bu_{1}, \ldots, \bu_{d})=g(\bu;u_{00},  \ldots,
u_{d0}),
\end{equation}
where $\bu_{i}=(u_{i0},\ldots, u_{ijk},\ldots,u_{i\alpha},\ldots)$
is the sequence of the coefficients of $\P_i$.

\begin{definition}
The $\delta$-polynomial defined in \bref{eq-cf1} is called the {\em
generalized Chow form} of $V$ or the prime $\delta$-ideal $\I(V)$
with respect to $\P_i(i=0,\ldots,d)$.
\end{definition}

%
%
Similar to Theorem \ref{chowhomogenous}, we can prove that the
generalized Chow form is a $\delta$-homogeneous $\delta$-polynomial
in each set of indeterminates $\bu_{i}$, but in this case the
homogeneous degree for distinct $\bu_{i}$ may be distinct.
The {\em order of the generalized Chow form w.r.t. $\bu_i$}, denoted
by $\ord(G, \bu_{i})$, is defined to be $\max_{u\in\bu_i}\ord(g,u)$.
 Now we will consider the order of the generalized Chow form.
\begin{theorem} \label{generalized-order}
Let $\mathcal {I}$ be a prime $\delta$-$\ff$-ideal with dimension
$d$ and order $h$ defined over $\ff$, and $G(\bu_{0}, \bu_{1},$ $
\ldots, \bu_{d})=g(\bu; u_{00}, u_{10}, \ldots, u_{d0})$ its
generalized Chow form. Then for a fixed $i$ between 0 and $d$,
$\ord(g,u_{i0})=h+s-s_i$ with $s=\sum_{l=0}^d s_l$. Moreover,
$\ord(G,\bu_i)=h+s-s_i$.
\end{theorem}
\proof  Use the notations as above in this section.
%
Let $\mathcal {I}_{d}=[\mathcal
{I},\P_{0},\ldots,\P_{i-1},\P_{i+1},$ $\dots,\P_{d}] \subset
\mathcal {F}\langle
\bu_{0},\ldots,\bu_{i-1},\bu_{i+1},\ldots,\bu_{d}\rangle
\{y_{1},\ldots,y_{n}\}$. By Theorem \ref{th-ord}, $\mathcal {I}_{d}$
is a prime $\delta$-ideal with dimension 0 and  order
$h+s_0+\cdots+s_{i-1}+s_{i+1}+\cdots+s_d=h+s-s_i$, where
$s=\sum_{l=0}^d s_l$.

Let $\I_{\zeta,\xi}$ $=[\mathcal {I},\P_{0},\ldots,\P_{d}]$ $\subset
\mathcal {F}\langle\bu\rangle
\{u_{00},\ldots,u_{d0},y_{1},\ldots,y_{n}\}$
 and
$\mathcal {I}_1=[\mathcal {I},\P_{0},\ldots,$ $\P_{d}]$ $=[\mathcal
{I}_{d},\P_0] \subset \mathcal {F}\langle \widehat{\bu}
\setminus\{u_{i0}\}\rangle \{u_{i0},y_{1},$ $\ldots,y_{n}\}$, where
$\widehat{\bu}=\bu_0\cup\cdots\cup\bu_d$.
%
%
Denote $\ord(G,u_{i0})$ by $h_1$.  Similar to the proof of Lemma
\ref{lm-cf2}, we can show that $\mathcal
{A}=g(\bu;u_{00},\ldots,u_{d0}),$ $\frac{\partial g}{\partial
u_{i0}^{(h_1)}}y_{1}-\frac{\partial g}{\partial u_{i10}^{(h_1)}},
\ldots, \frac{\partial g}{\partial
u_{i0}^{(h_1)}}y_{n}-\frac{\partial g}{\partial u_{in0}^{(h_1)}}$ is
a characteristic set of $\I_{\zeta,\xi}$ w.r.t. the elimination
ranking $u_{00}\prec \cdots u_{i-1,0} \prec u_{i+1,0}\prec
u_{d0}\prec u_{i0}\prec y_{1}\prec\cdots\prec y_{n}$. Clearly,
$\mathcal {I}_1$ is the $\delta$-ideal generated by $\I_{\zeta,\xi}$
in $\mathcal {F}\langle \widehat{\bu}\setminus\{u_{i0}\}\rangle
\{u_{i0},y_{1},\ldots,y_{n}\}$. Since
$\{u_{00},\ldots,u_{i-1,0},u_{i+1,0},\ldots,u_{d0}\}$ is a
parametric set of $\I_{\zeta,\xi}$, $\mathcal {A}$ is also a
characteristic set of  $\mathcal {I}_1$ w.r.t. the elimination
ranking $u_{i0}\prec y_{1} \prec \cdots \prec y_{n}$. Since
$\dim(\mathcal {I}_1)=0$, from Corollary~\ref{cor-0 order}, we have
$\ord(\mathcal {I}_1)=\ord(\mathcal {A})=\ord(g,u_{i0})$.

On the other hand, let $\mathcal {I}_1^{(l)}=[\mathcal
{I}_{d},u_{i0}^{(l)}+\sum_{j=1}^n\sum_{k=0}^{s_{i}}
u_{ijk}y_{j}^{(k)}+f_i] \subset \mathcal
{F}\langle\widehat{\bu}\setminus\{u_{i0}\}\rangle \{u_{i0},$
$y_{1},\ldots,y_{n}\}$ $(l=1,\ldots,s_i)$. Since $\dim(\mathcal
{I}_1^{(l)})=0$, $u_{i0}$ is a leading variable of $\mathcal
{I}_1^{(l)}$ for any ranking. Thus, by Lemma~\ref{th-order+1},
we have $\ord(\mathcal {I}_1^{(l+1)})=\ord(\mathcal {I}_1^{(l)})+1$,
which follows that $\ord(\mathcal {I}_1^{(s_i)})=\ord(\mathcal
{I}_1)+s_i$. And it is easy to see that $\ord(\mathcal
{I}_1^{(s_i)})=\ord(\mathcal {I}_{d})+s_i$. Indeed, let $\mathcal
{A}$ be a characteristic set of $\mathcal {I}_{d}$ w.r.t. some
orderly ranking $\mathscr{R}$, and let $t$ be the pseudo remainder
of $u_{i0}^{(s_i)}+\sum_{j=1}^n\sum_{k=0}^{s_{i}}
u_{ijk}y_{j}^{(k)}+f_i$ w.r.t. $\A$ under the ranking $\mathscr{R}$.
Clearly, $\ord(t,u_{i0})=s_i.$  It is obvious that for some orderly
ranking, $\{\mathcal{A}, t \}$ is a characteristic set of $\mathcal
{I}_{1}^{(s_i)}$ with $\lead(\mathcal{A})$ and $u_{i0}^{(s_i)}$ as
its leaders, so $\ord(\mathcal {I}_1^{(s_i)})=\ord(\mathcal
{I}_{d})+s_i$. Thus, $\ord(\mathcal {I}_1)=\ord(\mathcal
{I}_{d})=h+s-s_i$, and consequently, $\ord(g,u_{i0})=h+s-s_i$.

 It remains to show that $\ord(g,u_{ijk})$
$(j=1,\ldots,n; k=0,1,\ldots,s_i)$ and $\ord(g,u_{i\alpha})$ cannot
exceed $\ord(g,u_{i0})$. If $\ord(g,u_{ijk})=l>\ord(g,u_{i0})$, then
differentiate the identity
$g(\bu;\zeta_0,\ldots,\zeta_{d})=0$ w.r.t.
$u_{ijk}^{(l)}$, we have $\frac{\partial g}{\partial
u_{ijk}^{(l)}}(\bu;\zeta_0,\ldots,\zeta_{d})=0.$
Thus, $\frac{\partial g}{\partial u_{ijk}^{(l)}}$ can be divisible
by $g$, a contradiction.  So $\ord(g,u_{ijk})\leq\ord(g,u_{i0})$.
Similarly, we can prove that
$\ord(g,u_{i\alpha})\leq\ord(g,u_{i0})$. Thus,
$\ord(G,\bu_i)=\ord(g,u_{i0}).$ \qed

\vskip5pt In the following, we consider the factorization of the
generalized Chow form. Denote $h+s-s_i$ by $h_i$ ($i=0,\ldots,d$)
where $s=\sum_{l=0}^d s_l$. Now consider $G$ as a polynomial in
$u_{00}^{(h_0)}$ with coefficients in $\mathcal{F}_0=\mathcal
{F}\langle \widetilde{\bu}\rangle(u_{00},\ldots,u_{00}^{(h_0-1)})$,
where $\widetilde{\bu}=\cup_{i=0}^d\bu_i\backslash\{u_{00}\}$. Then,
in an algebraic extension field of $\mathcal{F}_0$, we have
 $$g=A\prod_{\tau=1}^{t_0} (u_{00}^{(h_0)}-\gamma_{\tau})$$
where $t_0=\deg(g,u_{00}^{(h_0)})$.
 Let $\xi_{\tau
\rho k}=g_{\tau \rho k}/g_{\tau 0}\,(\rho=1,\ldots,n ;\,
k=0,\ldots,s_0)$ and $\xi_{\tau \alpha}=g_{\tau \alpha}/g_{\tau 0}$,
where $g_{\tau \rho k}=\frac{\partial g}{\partial u_{0 \rho
k}^{(h_0)}}\Big| _{u_{00}^{(h_0)}=\gamma_{\tau}}$, $g_{\tau
\alpha}=\frac{\partial g}{\partial u_{0 \alpha}^{(h_0)}}\Big|
_{u_{00}^{(h_0)}=\gamma_{\tau}}$ and $g_{\tau 0}=\frac{\partial
g}{\partial u_{0 0}^{(h_0)}}\Big| _{u_{00}^{(h_0)}=\gamma_{\tau}}$.
Similarly as in Section 4.4, we can uniquely define the derivatives
of $\gamma_\tau$ and $\xi_{\tau\rho0}$ to make them elements in a
differential extension field of $\ff\langle\widetilde{\bu}\rangle$.
From $g(\bu;\zeta_{0},\ldots,\zeta_{d})=0$, if we differentiate this
equality w.r.t. $u_{0\rho k}^{(h_0)}$, then we have

\begin{equation} \label{equ-generalized}
 \overline{\frac{\partial g}{\partial u_{0\rho k}^{(h_0)}}}+
 \overline{\frac{\partial g}{\partial u_{00}^{(h_0)}}}
(-\xi_{\rho}^{(k)})=0
\end{equation}
And if we differentiate $g(\bu;\zeta_{0},\ldots,\zeta_{d})=0$ w.r.t. $u_{0\alpha}^{(h_0)}$, then
\begin{equation} \label{equ-generalized2} \overline{\frac{\partial g}{\partial u_{0\alpha
}^{(h_0)}}}+ \overline{\frac{\partial g}{\partial
u_{00}^{(h_0)}}}(-(\xi^{(s_{0})})^\alpha)=0
\end{equation}
where
$(\xi^{(s_0)})^\alpha=(\Y^{(s_0)})^\alpha|_{(y_1,\ldots,y_n)=(\xi_1,\ldots,\xi_n})$.
And in the above equations, $\overline{\frac{\partial g}{\partial
u_{0\rho k }^{(h_0)}}}$ and $\,\overline{\frac{\partial g}{\partial
u_{0\alpha }^{(h_0)}}}$  represent  $\frac{\partial g}{\partial
u_{0\rho k }^{(h_0)}}$ and $\,\frac{\partial g}{\partial
u_{0\alpha}^{(h_0)}}$  when substituting $u_{i0}$ by $\zeta_i$.
 For each $\rho=1,\ldots,n$ and $ k=0,\ldots,s_0$, multiplying the equations in ~(\ref{equ-generalized}) by $u_{0\rho
 k}$, and for $\alpha \in \mathbb{Z}^{n(s_0+1)}_{\geq 0}, \,1<|\alpha|\leq m_0$,  multiplying the equations in ~(\ref{equ-generalized2}) by $u_{0\alpha}$,
 then adding all of the equations obtained together, we have
\[\zeta_0 \overline{\frac{\partial g}{\partial u_{00}^{(h_0)}}}+\sum\limits_{\rho=1}^n\sum\limits_{k=0}^{s_0} u_{0\rho k}\overline{\frac{\partial g}{\partial u_{0\rho k}^{(h_0)}}}
+\sum_{
 \begin{array}{c} \alpha \in \mathbb{Z}^{n(s_0+1)}_{\geq 0} \\ 1<|\alpha|\leq m_0
 \end{array}}u_{0 \alpha}\overline{\frac{\partial g}{\partial u_{0\alpha}^{(h_0)}}}=0.\] Thus, the
 $\delta$-polynomial
$u_{00}\frac{\partial g}{\partial
u_{00}^{(h_0)}}+\sum\limits_{\rho=1}^n\sum\limits_{k=0}^{s_0}
u_{0\rho k}\frac{\partial g}{\partial u_{0\rho
k}^{(h_0)}}+\sum\limits_{
 \begin{array}{c} \alpha \in \mathbb{Z}^{n(s_0+1)}_{\geq 0} \\ 1<|\alpha|\leq
 m_0
 \end{array}}u_{0 \alpha}\frac{\partial g}{\partial
u_{0\alpha}^{(h_0)}}$ \, vanishes at $
(u_{00},\ldots,u_{d0})=(\zeta_{0},\ldots,\zeta_{d}).$ Since it is at
most of the same order as $g$, it must be divisible by $g$. And
since it has the same degree as $g$, there exists some $a\in
\mathcal {F}$ such that $$u_{00}\frac{\partial g}{\partial
u_{00}^{(h_0)}}+\sum\limits_{\rho=1}^n\sum\limits_{k=0}^{s_0}
u_{0\rho k}\frac{\partial g}{\partial u_{0\rho
k}^{(h_0)}}+\sum\limits_{
 \begin{array}{c} \alpha \in \mathbb{Z}^{n(s_0+1)}_{\geq 0} \\ 1<|\alpha|\leq
 m_0
 \end{array}}u_{0 \alpha}\frac{\partial g}{\partial
u_{0\alpha}^{(h_0)}}=ag.$$
Setting $u_{00}^{(h_0)}=\gamma_{\tau}$ in both sides of the above
equation, we have
$$u_{00}g_{\tau 0}+\sum_{\rho=1}^n\sum_{k=0}^{s_0} u_{0\rho k}g_{\tau \rho k}+\sum_{
 \begin{array}{c} \alpha \in \mathbb{Z}^{n(s_0+1)}_{\geq 0} \\ 1<|\alpha|\leq m_0
 \end{array}}u_{0 \alpha}g_{\tau \alpha}=0.$$ Or,
\[u_{00}+\sum_{\rho=1}^n\sum_{k=0}^{s_0} u_{0\rho k}\xi_{\tau \rho k}+\sum_{
 \begin{array}{l} \alpha \in \mathbb{Z}^{n(s_0+1)}_{\geq 0} \\ 1<|\alpha|\leq m_0
 \end{array}}u_{0 \alpha}\xi_{\tau \alpha}=0. \]
Then, we have
\[\begin{array}{c}(u_{00}+\sum\limits_{\rho=1}^n\sum\limits_{k=0}^{s_0} u_{0\rho
k}\xi_{\tau \rho k}+\sum\limits_{
 \begin{array}{l} \alpha \in \mathbb{Z}^{n(s_0+1)}_{\geq 0} \\ 1<|\alpha|\leq m_0
 \end{array}}u_{0 \alpha}\xi_{\tau \alpha})^{(h_0)}\\=\gamma_{\tau}+(\sum\limits_{\rho=1}^n\sum\limits_{k=0}^{s_0} u_{0\rho k}\xi_{\tau \rho k}+\sum\limits_{
 \begin{array}{l} \alpha \in \mathbb{Z}^{n(s_0+1)}_{\geq 0} \\ 1<|\alpha|\leq m_0
 \end{array}}u_{0 \alpha}\xi_{\tau \alpha})^{(h_0)}=0.\end{array}.\]

We have the following theorem
\begin{theorem}
Let $G(\bu_{0},\bu_{1},\ldots,\bu_{d})$ be the generalized Chow form
of a $\delta$-$\ff$-variety of dimension $d$ and order $h$. Then,
there exist $\xi_{\tau \rho}(\rho=1,\ldots,n;\tau=1,\ldots,t_0)$ in
a $\delta$-extension field of $\mathcal
{F}\langle\widetilde{\bu}\rangle$ such that
\begin{eqnarray}\label{eq-generafactor}
G(\bu_{0},\bu_{1},\ldots,\bu_{d})=A(\bu_{0},\bu_{1},\ldots,\bu_{d})\prod^{t_0}_{\tau=1}
\P_0(\xi_{\tau1},\ldots,\xi_{\tau n})^{(h_0)}
%
\end{eqnarray}
where  $A(\bu_{0},\bu_{1},\ldots,\bu_{d})$ is in $\mathcal
{F}\{\bu_{0},\bu_{1},\ldots,\bu_{d}\}$ and
$t_0=\deg(G,u_{00}^{(h_0)})$.
\end{theorem}
\proof From what we have proved,
{\tiny $$G(\bu_{0},\bu_{1},\ldots,\bu_{d})=A(\bu_{0},\bu_{1},\ldots,\bu_{d})\prod^{t_0}_{\tau=1}\Big(u_{00}+\sum_{\rho=1}^n\sum_{k=0}^{s_0}
u_{0\rho k}\xi_{\tau \rho k}+\sum_{
 \begin{array}{l} \alpha \in \mathbb{Z}^{n(s_0+1)}_{\geq 0} \\ 1<|\alpha|\leq
 m_0
 \end{array}}u_{0 \alpha}\xi_{\tau \alpha}\Big)^{(h_0)}.$$}

Denote $\xi_{\tau \rho0}$ by $\xi_{\tau \rho}$. To complete the
proof, it remains to show that $\xi_{\tau \rho k}=(\xi_{\tau \rho
0})^{(k)}(k=1,\ldots,s_0)$ and $\xi_{\tau \alpha}=\prod
\limits_{\rho=1}^n \prod \limits_{j=0}^{s_0} \Big( \big(\xi_{\tau
\rho 0}\big)^{(j)}\Big)^{\alpha_{\rho j}}$. From
equation~(\ref{equ-generalized}) and
equation~(\ref{equ-generalized2}) , we have
$\xi_{\rho}^{(k)}=\overline{\frac{\partial g}{\partial u_{0\rho k
}^{(h_0)}}}\Big/ \overline{\frac{\partial g}{\partial
u_{00}^{(h_0)}}}
%
\,(k=0,\ldots,s_0)$ and
$(\xi^{(s_0)})^\alpha=\overline{\frac{\partial g}{\partial
u_{0\alpha }^{(h_0)}}}\Big/\overline{\frac{\partial g}{\partial
u_{00}^{(h_0)}}}$. So we have  the equalities:
$\Big(\overline{\frac{\partial g}{\partial u_{0\rho 0
}^{(h_0)}}}\Big/\overline{\frac{\partial g}{\partial
u_{00}^{(h_0)}}} \Big)^{(k)}=\overline{\frac{\partial g}{\partial
u_{0 \rho k}^{(h_0)}}}\Big/\overline{\frac{\partial g}{\partial
u_{00}^{(h_0)}}}$ and $\overline{\frac{\partial g}{\partial
u_{0\alpha }^{(h_0)}}}\Big/\overline{\frac{\partial g}{\partial
u_{00}^{(h_0)}}}=\prod \limits_{\rho=1}^n \prod \limits_{j=0}^{s_0}
\bigg(\Big(\overline{\frac{\partial g}{\partial u_{0 \rho0
}^{(h_0)}}}\bigg/\overline{\frac{\partial g}{\partial
u_{00}^{(h_0)}}}\Big)^{(j)}\bigg)^{\alpha_{\rho j}}$. Thus,
$\Big(\frac{\partial g}{\partial u_{0\rho 0
}^{(h_0)}}\Big/\frac{\partial g}{\partial
u_{00}^{(h_0)}}\Big)^{(k)}-\frac{\partial g}{\partial u_{0\rho k
}^{(h_0)}}\Big/\frac{\partial g}{\partial u_{00}^{(h_0)}}$ and
$\frac{\partial g}{\partial u_{0\alpha }^{(h_0)}}\Big/\frac{\partial
g}{\partial u_{00}^{(h_0)}}-\prod \limits_{\rho=1}^n \prod
\limits_{j=0}^{s_0} \bigg(\Big(\frac{\partial g}{\partial u_{0 \rho0
}^{(h_0)}}\bigg/\frac{\partial g}{\partial
u_{00}^{(h_0)}}\Big)^{(j)}\Big)^{\alpha_{\rho j}}$ vanish at
$(u_{00},$ $\ldots,u_{d0})=(\zeta_{0},\ldots,\zeta_{d})$. Similarly
as in the proof of Theorem~\ref{th-zerochowform}, we can see that both of the differential polynomials
$\Big(\frac{\partial g}{\partial u_{0\rho 0
}^{(h_0)}}\Big/\frac{\partial g}{\partial
u_{00}^{(h_0)}}\Big)^{(k)}-\frac{\partial g}{\partial u_{0\rho k
}^{(h_0)}}\Big/\frac{\partial g}{\partial u_{00}^{(h_0)}}$ and
$\frac{\partial g}{\partial u_{0\alpha }^{(h_0)}}\Big/\frac{\partial
g}{\partial u_{00}^{(h_0)}}-\prod \limits_{\rho=1}^n \prod
\limits_{j=0}^{s_0} \bigg(\Big(\frac{\partial g}{\partial u_{0 \rho0
}^{(h_0)}}\bigg/\frac{\partial g}{\partial
u_{00}^{(h_0)}}\Big)^{(j)}\bigg)^{\alpha_{\rho j}}$ vanish at
$u_{00}^{(h_0+j)}=\gamma_{\tau}^{(j)}(j  \geq 0)$.  Thus, $\xi_{\tau
\rho 0}^{(k)}=\xi_{\tau \rho k}$ and $\xi_{\tau \alpha}-\prod
\limits_{\rho=1}^n \prod \limits_{j=0}^{s_0} \Big( \big(\xi_{\tau
\rho 0}\big)^{(j)}\Big)^{\alpha_{\rho j}}=0.$ The proof is
completed. \qed

\begin{theorem}
The points $(\xi_{\tau1},\ldots,\xi_{\tau n})\,(\tau=1,\ldots,t_0)$
in \bref{eq-generafactor} are  generic points of the
$\delta$-$\ff$-variety $V$, and satisfy the equations
\[\P_\sigma(y_1,\ldots,y_n)=u_{\sigma
0}+\sum_{\rho=1}^n\sum_{k=0}^{s_{\sigma}} u_{\sigma \rho
k}y_{\rho}^{(k)}+f_{\sigma} =0\,(\sigma=1,\ldots,d)\]
\end{theorem}
\proof The proof is similar to that  of Theorem~\ref{th-zerochowform}. \qed

\begin{theorem} \label{th-generalized-sf}
Let $G(\bu_{0},\ldots,\bu_{d})$ be the generalized  Chow form of $V$
and $S_{G}=\frac{\partial G}{\partial u_{00}^{(h_0)}}$ with
$\ord(G,u_{00})=h_0$.
Suppose that $\bu_i(i=0,\ldots,d)$ are $\delta$-specialized to sets
$\bv_i$ of specific elements in  $\mathcal{E}$ and
$\overline{\P}_{i}\,(i=0,\ldots,d)$ are obtained by substituting
$\bu_i$ by $\bv_i$ in $\P_i$.
If\, $\overline{\P}_i=0(i=0,\ldots,d)$ meet $V$, then
$G(\bv_{0},\ldots,\bv_{d})$ $=0$. Furthermore, if
$G(\bv_{0},\ldots,\bv_{d})=0$ and $S_{G}(\bv_{0},\ldots,\bv_{d})\neq
0$, then the $d+1$ $\delta$-hypersurfaces $\overline{\P}_{i}=0$
$(i=0,\ldots,d)$ meet $V$.
%
\end{theorem}
\proof The proof is similar to that  of Theorem~\ref{th-sf}. \qed

\subsection{Differential resultant of multivariate differential polynomials}

As an application of the generalized Chow form, we can define the
differential resultant of $n+1$ generic $\delta$-polynomials in $n$
variables.
Let $\mathcal{I}=[0]$ be the $\delta$-ideal generated by $0$ in
$\mathcal{F}\{\Y\}$. Then $\dim(\mathcal{I})=n$. Let
$G(\bu_0,\bu_1,\ldots,\bu_n)$ be the generalized Chow form for
$\mathcal{I}$. Then we will define $G(\bu_0,\bu_1,\ldots,\bu_n)$ to
be the differential resultant for the $n+1$ generic
$\delta$-polynomials given in \bref{eq-gpol}.

\begin{definition}
The {\em differential resultant} for the $n+1$ generic
$\delta$-polynomials $\P_i$ in \bref{eq-gpol} is defined to be the
generalized Chow form of $\mathcal{I}=[0]$ associated with these
$\P_i$, and will be denoted by
$R(\bu_0,\ldots,\bu_n)=G(\bu_0,\ldots,\bu_n)$.
\end{definition}
\begin{theorem}\label{th-resc}
Let $R(\bu_{0},\ldots,\bu_{n})$ be the differential resultant of the
$n+1$ $\delta$-polynomials $\P_{0},\ldots,\P_{n}$ given in
\bref{eq-gpol} with $\ord(\P_i)=s_i$ and $\deg(\P_i)=m_i$, where
$\bu_{i}=(u_{i0},\ldots,$
$u_{ijk},\ldots,u_{i\alpha},\ldots)\,(i=0,\ldots,n).$ Denote
$s=\sum\limits_{i=0}^n s_{i}$, $D=\max_{i=0}^n\{m_i\}$ and
$\bu=\cup_{i=0}^n\bu_{i}\backslash\{u_{i0}\}$. Then there exist
$h_{jk}\in \mathcal {F}\langle
\bu\rangle[y_1,\ldots,y_n,\ldots,y_{1}^{(s)},\ldots,y_{n}^{(s)}]$
such that
$$R(\bu_{0},\ldots,\bu_{n})=\sum_{j=0}^n\sum_{k=0}^{s-s_j}h_{jk}\delta^{k}\P_{j}.$$
Moreover, the degree of $h_{jk}$ in $\Y$ is bounded by $(sn+n)^2
D^{sn+n}+D(sn+n)$.
\end{theorem}
\proof  Let $\mathcal {J}$ be the $\delta$-ideal generated by
$\P_{0},\ldots,\P_{n}$ in  $\mathcal
{F}\{\bu_0,\ldots,\bu_{n},y_1,\ldots,y_n\}$. Let $\mathscr{R}$ be
the elimination ranking $\bu\prec y_n \prec \cdots \prec y_1 \prec
u_{n0}\prec \cdots \prec u_{00}$ with arbitrary ranking endowed on
$\Theta(\bu)=(\theta u: u\in \bu; \theta \in \Theta)$. Clearly,
$\mathcal {J}$ is a prime $\delta$-ideal with $\P_{0},\ldots,\P_{n}$
as its characteristic set w.r.t. $\mathscr{R}$. Thus, $\bu \cup
\{y_1,\ldots,y_n\}$ is a parametric set of $\mathcal {J}$. From the
definition of $R$, $R\in \mathcal {J}.$  In
$R(\bu_{0},\ldots,\bu_{n})=G(\bu_{0},\ldots,\bu_{n})=g(\bu;u_{00},\ldots,u_{n0})$,
let $u_{i0}\,(i=0,\ldots,n)$ be replaced respectively by $$
\P_i-\sum_{j=1}^n\sum_{k=0}^{s_i} u_{ijk}y_{j}^{(k)} -\sum_{
 \begin{array}{c} \alpha \in \mathbb{Z}^{n(s_i+1)}_{\geq 0} \\ 1<|\alpha|\leq
 m_i
 \end{array}}u_{i \alpha}(\Y^{(s_i)})^{\alpha},\, (i=0,\ldots,d),$$
 and let $R$ be expanded as a polynomial in $\P_{0},\ldots,\P_{n}$
 and their derivatives. The term not involving $\P_{0},\ldots,\P_{n}$
 or their derivatives will be a $\delta$-polynomial only involving $\bu \cup
\{y_1,\ldots,y_n\}$ which also belong to $\mathcal {J}$. Since
$\mathcal {J}\bigcap \mathcal {F}\{\bu,y_1,\ldots,y_n\}=\{0\}$, such
term will be identically zero. So $R$ is a linear combinations of
$\P_0,\ldots,\P_n$ and some of their derivatives. Since
$\ord(R,u_{i0})=s-s_i$, the above expansion for $R$ involves $\P_i$
only up to the order $s-s_i$ and the coefficients in the linear
combination are $\delta$-polynomials in $\mathcal
{F}\{\bu\}[y_1,\ldots,y_n,\ldots,y_{1}^{(s)},\ldots,y_{n}^{(s)}]$.
Denote $\mathcal {R}=\mathcal
{F}\langle\bu\rangle[y_1,\ldots,y_n,\ldots,y_{1}^{(s)},\ldots,y_{n}^{(s)}]$.
Thus, $R\in (\delta^{s-s_0}\P_0,\ldots,\delta
\P_0,\P_0,\ldots,\delta^{s-s_n}\P_n,\ldots,\delta \P_n,\P_n)
\subseteq \mathcal {R}$,
which implies that $(\delta^{s-s_0}\P_0,\ldots,\delta
\P_0,\P_0,\ldots,\delta^{s-s_n}\P_n,\ldots,\delta \P_n,\P_n)$ in
$\mathcal {R}$ is the unit ideal.
By \cite[Theorem 1]{ritt0}, there exist $A_{jk}\in\mathcal {R}$ with $\deg(A_{jk})\leq (sn+n)^2 D^{sn+n}+D(sn+n)$
such that
$$1=\sum_{j=0}^n\sum_{k=0}^{s-s_j}A_{jk} \delta^{k}\P_{j},$$ where $D=\max\{m_0,m_1,\ldots,m_n\}$.
If we multiply the above equation by $R$ and denote $A_{jk}R$ by
$h_{jk}$, we complete the proof. \qed

As a consequence of the above five theorems proved in this section,
the properties of the differential resultant listed in Theorem
\ref{th-main3} are proved.

Let $R(\bu_0,\ldots,\bu_n)$ be the differential resultant for the
$n+1$ generic $\delta$-polynomials $\P_i$ in \bref{eq-gpol}.
When each $\bu_i$ are specialized to specific elements
$\bv_i\in\ff^{n+1}$, $R(\bv_0,\ldots,$ $\bv_n)$ is called the {\em
differential resultant} of $\overline{\P}_i(i=0,\ldots,n)$ which are
obtained by replacing $\bu_i$ by $\bv_i$ in $\P_i$.
By Theorem \ref{th-resc}, vanishing of the differential resultant
for $n+1$ $\delta$-polynomials in  $\ff\{\Y\}$ is a necessary
condition for them to have a common solution.

\begin{remark}\label{rm-res}
It is easy to see that if $s_i=0$, then the differential resultant
of $\P_i(i=0,\ldots,n)$ becomes the Macaulay resultant for $n+1$
polynomials in $n$ variables \cite{cox2,joun1}. From Theorem
\ref{th-main3}, we see that the differential resultant has similar
properties to that of the Macaulay resultant. Special attention
should be payed to the second property which is a differential
analog to the so-called Poisson type formulas for algebraic
resultants \cite{Pedersen}. Also note that many properties of the
Macaulay resultant are yet to be extended to the differential case.
The most significant one might be to find a matrix representation
for the differential resultant similar to the one given in
\cite[p.102]{cox2}. Note that such a formula was claimed to be given
in \cite{dres1,dres2}, which is not correct as we mentioned in
Section 1 of this paper.
As a latest development,  we defined the differential sparse
resultant and proposed a single exponential algorithm to compute it
\cite{sdres}.

\end{remark}

Similar to the differential Chow form, the differential resultant
can be computed with the differential elimination algorithms
\cite{boulier2010,ardm1,ritt,sit,alexey}.
%

\begin{example}
The simplest nonlinear differential resultant is the case
$n=1,d_0=d_1=2,s_0=0,s_1=1$. Denote $y_1$ by $y$.
Let $\P_0 = u_{00} + u_{01} y + u_{02} y^2$, $\P_1 = u_{10} + u_{11}
y + u_{12} y' + u_{13} y^2 + u_{14} yy' +  u_{15} (y')^2$.
Then the differential resultant for $\P_0$ and $\P_1$ is a
$\delta$-polynomial $R(\bu_0,\bu_1)$ such that $\ord(R,\bu_0)=1$,
$\ord(R,\bu_1)=0$ and $R$ is $\delta$-homogenous of degree $8$ in
$\bu_0$ and degree 2 in $\bu_1$ respectively. Totally, $R$ has 206
terms. Moreover, $R$ has a matrix representation which is a factor
of the determinant of the coefficient matrix of
$\P_0,y'\P_0,y^2\P_0,yy'\P_0,y'^2\P_0,\P_0',y\P_0',y'\P_0',yy'\P_0',y'^2\P_0,\P_1,$
$y\P_1,y'\P_0,yy'\P_1$ w.r.t. the monomials
$\{y^{l_0}(y')^{l_1}|0\leq l_0\leq 4,0\leq l_1< 4,l_0+l_1\leq4\}$.
\end{example}

\section{Conclusion}
\label{sec-conc}

In this paper, an intersection theory for generic differential
polynomials is presented by giving the explicit formulas for the
dimension and order of the intersection of an irreducible
differential variety with a generic differential hypersurface. As a
consequence, we show that the differential dimension conjecture is true for
generic differential polynomials.

The  Chow form for an irreducible differential variety is defined.
Most of the properties of the algebraic Chow form
have been extended to its differential counterpart. In particular,
we introduce the concept of differential Chow quasi-variety for a
special class of differential algebraic cycles.
Furthermore, the generalized Chow form for an irreducible
differential variety is also defined and its properties are proved.
As an application of the generalized differential Chow form, we can
give a rigorous definition  for the differential resultant and
establish its properties which are similar to that of the Sylvester
resultant of two univariate polynomials and the Macaulay resultant
of multivariate polynomials.

The results given in this paper enrich the field of differential
algebraic geometry. Further, many new problems can be raised
naturally. Some of them are already mentioned in Remarks and
\ref{re-sf} and \ref{rm-res}. We mentioned in Section 1 that the
algebraic Chow form has many important applications. It is very
interesting to see whether some of these applications can be
extended to the differential case.

As we mentioned in Section 5, the theory of differential Chow quasi-varieties is
not fully developed and the main difficulty is to develop an elimination theory
for mixed systems with both algebraic and differential equations.
%
%
%
%

In this paper, we only consider Chow forms for affine differential
varieties. It is not difficult to extend most of the results in this
paper to differential Chow form of differentially projective
varieties. Note that differentially projective varieties were
defined by Kolchin in \cite{kol51}. It is expected that Theorems
\ref{th-sf},  \ref{th-var1}, and \ref{th-cv1} could improved for
differentially projective varieties.
%
%

\newpage

\end{document}